\documentclass[10pt]{amsart}
\usepackage[all]{xy}
\usepackage{graphicx}
\usepackage[colorlinks=true, urlcolor=blue,bookmarks=true,bookmarksopen=true,
citecolor=blue]{hyperref}

\usepackage{amssymb,amsmath,amsthm,amscd}

\newcommand{\Z}{\mathbb{Z}}
\newcommand{\R}{\mathbb{R}}
\newcommand{\C}{\mathbb{C}}
\renewcommand{\i}{\,{i}\,}
\renewcommand{\L}{\mathrm{L}}
\newcommand{\RegClass}{\mathrm{C}}
\newcommand{\ind}{\mathrm{ind}}
\newcommand{\im}{\mathrm{im}}
\newcommand{\PathSpace}{\mathcal{P}}
\newcommand{\Id}{\mathrm{Id}}
\newcommand{\A}{\mathcal{A}}
\newcommand{\Mp}{{\mathcal{M}'}}
\newcommand{\M}{\mathcal{M}}
\newcommand{\dbar}{\bar\partial}
\newcommand{\dbartld}{\tilde\partial}
\newcommand{\cutoff}[2]{\chi^{#1}_{#2}}
\ifx\llbracket\undefined%
\newcommand{\strip}[1]{#1\!#1}%
\else
\newcommand{\strip}[1]{%
\ifx#1[\llbracket\else%
\ifx#1]\rrbracket\else%
\errmessage{bad usage of \strip}\fi\fi}
\fi

\newtheorem{theo}{Theorem}[section]
\newtheorem{cor}[theo]{Corollary}
\newtheorem{prop}[theo]{Proposition}
\newtheorem{rem}[theo]{Remark}
\newtheorem{lem}[theo]{Lemma}
\newtheorem{defi}[theo]{Definition}

\begin{document}

\title{Lagrangian Intersections and the Serre Spectral Sequence}

\author[J.-F. Barraud,\ O. Cornea]{Jean-Fran\c{c}ois Barraud and Octav Cornea$^{1}$}
\date{November 23, 2003; Revised January 27, 2006}

\address{J-F.B.: UFR de Math\'ematiques\newline
\indent Universit\'e de Lille 1\newline
\indent 59655 Villeneuve
d'Ascq\newline
\indent France}
\email{barraud@agat.univ-lille1.fr}

\address{O.C.: University of Montr\'eal\newline
\indent Department of Mathematics and Statistics
\newline \indent CP 6128 Succ. Centre Ville
\newline \indent Montr\'eal, QC H3C 3J7
\newline\indent Canada}
\email{cornea@dms.umontreal.ca}

\begin{abstract}\hskip5pt
For a transversal pair of closed Lagrangian submanifolds $L, L'$ of a symplectic manifold $M$ so that $\pi_{1}(L)=\pi_{1}(L')=0=c_{1}|_{\pi_{2}(M)}=\omega|_{\pi_{2}(M)}$ and a generic almost complex structure $J$ we construct an invariant with a high
homotopical content which consists in the pages of order $\geq 2$ of
a spectral sequence whose differentials provide an algebraic measure of
the high-dimensional moduli spaces of pseudo-holomorpic strips of finite energy that join $L$ and $L'$. When $L$ and $L'$ are hamiltonian isotopic, we show that the pages of the spectral sequence coincide (up to a horizontal translation) with the terms of the Serre-spectral sequence of the path-loop fibration $\Omega L\to PL\to L$ and we deduce some applications.
\end{abstract}
\maketitle
 \addtocounter{footnote}{+1}
\footnotetext{Partially supported by an NSERC Discovery grant and by a FQRNT group research grant.}
\tableofcontents

\section{Introduction.}\label{sec:intro}
Consider a symplectic manifold $(M,\omega)$ which is convex at
infinity together with two closed (compact, connected, without boundary) Lagrangian submanifolds $L$, $L'$
in general position. We fix from now on the dimension of $M$ to be $2n$.
Unless otherwise stated we assume in this introduction that
\begin{equation}\label{eq:connectivity}
\pi_{1}(L)=\pi_{1}(L')=0=c_{1}|_{\pi_{2}(M)}=\omega|_{\pi_{2}(M)}
\end{equation}
and we shall keep this assumption in most of the paper.

One of the main tools in symplectic
topology is Floer's machinery (see \cite{Sal} for a recent exposition)
which, once a generic almost complex
structure compatible with $\omega$ is fixed on $M$, gives rise to a
Morse-type chain complex $(CF_{\ast}(L,L'), d_{F})$ such that
$CF_{\ast}(L,L')$ is the free $\Z/2$-vector space generated by
(certain) intersection points in $L\cap L'$ and $d_{F}$ counts the
number of connecting orbits (also called ``Floer trajectories" - in this case they are
pseudo-holomorphic strips) joining intersection points of relative
(Maslov) index equal to $1$ (elements of Floer's construction are
recalled in \S\ref{sec:spec}). In this construction are only involved
$1$ and $2$-dimensional moduli spaces of connecting trajectories.

\

The present paper is motivated by the following problem: {\em extract
out of the structure of higher dimensional moduli spaces of Floer trajectories
useful homotopical type data which is not limited to Floer
homology (or cohomology)}.

\

This question is natural because the properties of Floer trajectories
parallel those of negative gradient flow lines of a Morse function
(defined with respect to a generic riemannian metric) and the
information encoded in the Morse-Smale negative-gradient flow of such
a function is much richer than only the homology of the ambient
manifold. Indeed, in a series of papers on ``Homotopical Dynamics"
\cite{Co1}\cite{Co2}\cite{Co3}\cite{Co4} the second author has
described a number of techniques which provide ways to ``quantify"
algebraically the homotopical information carried by a flow.
In particular, in \cite{Co2} and \cite{Co4} it is shown how to
estimate the moduli spaces arising in the Morse-Smale context
when the critical points involved are consecutive in the
sense that they are not joined by any ``broken" flow line.
However, the natural problem of finding a computable algebraic
method to ``measure" \emph{general, high dimensional moduli
spaces} of connecting orbits has remained open till now
even in this simplest Morse-Smale case.
Of course, in the Floer case, a significant additional difficulty is
that there is no ``ambient" space with a meaningful topology.

\

We provide a solution to this problem in the present paper.
The key new idea
can be summarized as follows: \emph{in ideal conditions,
the ring of coefficients used to define a Morse type complex can be enriched so
that the resulting chain complex contains information about high
dimensional moduli spaces of connecting orbits}.

\

Roughly, this ``enrichment" of the coefficients is achieved by viewing the
relevant connecting orbits as loops in an appropriate space $\tilde{L}$ in which
the finite number of possible ends of the orbits are naturally
identified to a single point. The ``enriched" ring is then provided by the (cubical) chains of the pointed (Moore) loop space of $\tilde{L}$.
This ring turns out to be sufficiently rich algebraically such as to encode reasonably well the geometrical complexity  of the combinatorics of the higher dimensional moduli spaces.  Operating with the new chain complex is no more difficult than using the usual Morse complex. In particular, there is a natural filtration
of this complex and the pages of order higher than $2$ of the
associated spectral sequence (together with the respective differentials) provide
our invariant. Moreover, these pages are computable purely algebraically in
certain important cases.

\

This technique is quite powerful and is general enough so that
each manifestation of a Morse type complex in the literature offers a potential
application. From this point of view, our construction is certainly just a
first - and, we hope, convincing - step.

\subsection{The main result.} Fix a path-connected
component $\mathcal{P}_{\eta}(L,L')$ of the space
$\mathcal{P}(L,L')=\{\gamma\in C^{\infty}([0,1], M) : \gamma(0)\in
L$, $\gamma(1)\in L'\}$. The construction of Floer homology
depends on the choice of such a component. We denote the corresponding
Floer complex by $CF_{\ast}(L,L';\eta)$ and the resulting homology by
$HF_{\ast}(L,L';\eta)$. In case $L'=\phi_{1}(L)$ with $\phi_{1}$ the
time $1$-map of a hamiltonian isotopy $\phi :M\times [0,1]\to M$
(such a $\phi_{1}$ is called a hamiltonian diffeomorphism) we denote
by $\mathcal{P}(L,L';\eta_{0})$ the path-component of $\mathcal{P}(L,L')$ such that
$[\phi_{t}^{-1} (\gamma(t)]=0\in\pi_{1}(M,L)$ for some (and thus all)
$\gamma \in \eta_{0}$. We omit $\eta_{0}$ in the notation for the
Floer complex and Floer homology in this case. Given two spectral
sequences $(E^{r}_{p,q}, d^{r})$ and $(G^{r}_{p,q},d^{r})$ we say
that they are {\em isomorphic up to translation} if there exists an
integer $k$ and an isomorphism of chain complexes
$(E^{r}_{\ast+k,\star},d^{r})\approx (G^{r}_{\ast,\star},d^{r})$ for
all $r$. Recall that the path-loop fibration $\Omega L\to PL \to L$
of base $L$ has as total space the space of based paths in $L$ and as
fibre the space of based loops. Given two points $x,y\in L\cap L'$ we
denote by $\mu(x,y)$ their relative Maslov index and by
$\mathcal{M}(x,y)$ the non-parametrized moduli space of Floer
trajectories connecting $x$ to $y$ (see \S\ref{sec:spec} for the
relevant definitions). We denote by $\mathcal{M}$ the disjoint union
of all the $\mathcal{M}(x,y)$'s. We denote by $\mathcal{M}'$ the
space of all parametrized pseudo-holomorphic strips. All homology
groups below have $\Z/2$-coefficients.

\begin{theo}\label{theo:main}
Under the assumptions above there exists a spectral sequence
$$EF(L,L';\eta)=(EF^{r}_{p,q}(L,L';\eta), d_{F}^{r}) \ , r\geq 1$$
with the following properties:
\begin{itemize}
\item[a.] If $\phi:M\times [0,1]\to M$ is a hamiltonian isotopy, then
$(EF^{r}_{p,q}(L,L';\eta), d^{r})$ and
$(EF^{r}_{p,q}(L,\phi_{1}L';\phi_{1}\eta),d^{r})$ are isomorphic up
to translation for $r\geq 2$ (here $\phi_{1}\eta$ is the component
represented by $\phi_{t}(\gamma(t))$ for $\gamma\in \eta$).
\item[b.]
$EF^{1}_{p,q}(L,L';\eta)\approx CF_{p}(L,L';\eta)\otimes H_{q}(\Omega
L)$, $EF^{2}_{p,q}(L,L';\eta)\approx HF_{p}(L,L';\eta)\otimes
H_{q}(\Omega L)$.
\item[c.] If $d_{F}^{r}\not=0$, then there exist
points $x,y\in L\cap L'$ such that $\mu(x,y)=r$ and
$\mathcal{M}(x,y)\not=\emptyset$. \item[d.] If $L'=\phi'L$ with $\phi'$ a
Hamiltonian diffeomorphism, then for $r\geq 2$ the spectral sequence
$(EF^{r}(L,L'),d_{F}^{r})$ is isomorphic up to translation to the
$\Z/2$-Serre spectral sequence of the path loop fibration $\Omega
L\to PL\to L$.
\end{itemize}
\end{theo}

\subsection{Comments on the main result.}
We survey here the main features of the theorem.
\subsubsection{Geometric interpretation of the spectral
sequence.} \label{subsubsec:geom_int}
The differentials appearing in
the spectral sequence $EF(L,L';\eta)$ provide an algebraic measure of
the Gromov compactifications $\overline{\mathcal{M}}(x,y)$ of the
moduli spaces $\mathcal{M}(x,y)$ in - roughly - the following sense.
Let $\tilde{L}$ be the quotient topological space obtained by
contracting to a point a path in $L$ which passes through each point
in $L\cap L'$ and is homeomorphic to $[0,1]$. Let $\tilde{M}$ be the
space obtained from $M$ by contracting to a point the same path.
Clearly, $L$ and $\tilde{L}$ (as well as $M$ and $\tilde{M}$) have
the same homotopy type. Each point $u\in \mathcal{M}(x,y)$ is
represented by a pseudo-holomorphic strip $u: \R\times [0,1]\to M$
with $u(\R,0) \subset L$, $u(\R,1)\subset L'$ and such that
$\lim_{s\to-\infty}u(s,t)=x, \lim_{s\to +\infty}u(s,t)=y,\ \forall
t\in [0,1]$. Clearly, to such a $u$ we may asociate the path in $L$
given by $s \to u(s,0)$ which joins $x$ to $y$. Geometrically, by
projecting onto $\tilde{L}$, this associates to $u$ an element of
$\Omega \tilde{L}\simeq \Omega L$. The action functional can be used to
reparametrize uniformly the loops obtained in this way so that
 the resulting application extends in a continuous manner to the whole of $\overline{\mathcal{M}}(x,y)$ thus
producing a continuous map $\Phi_{x,y}:\overline{\mathcal{M}}(x,y)\to
\Omega L$. The space $\overline{\mathcal{M}}(x,y)$ has the structure
of a manifold with boundary with corners (see \S\ref{sec:spec} and
\S\ref{sec:appendix}) which is compatible with the maps $\Phi_{x,y}$.
If it would happen that
$\partial\overline{\mathcal{M}}(x,y)=\emptyset$ one could measure
$\overline{\mathcal{M}}(x,y)$ by the image in $H_{\ast}(\Omega L)$ of
its fundamental class via the map $\Phi_{x,y}$. This boundary is
almost never empty so this elementary idea fails. However, somewhat
miraculously, the differential $d^{\mu(x,y)}_{F}$ of $EF(L,L';\eta)$
reflects homologically what is left of
$\Phi_{x,y}((\overline{\mathcal{M}}(x,y))$ after ``killing" the
boundary $\partial\overline{\mathcal{M}}(x,y)$.

>From this perspective, it is clear that it is not so important where
the spectral sequence $EF(L,L';\eta)$ converges but rather whether it
contains many non-trivial differentials.

\subsubsection{Role of the Serre spectral sequence.} Clearly,
point a. of the theorem shows that the pages of order higher than $1$
of the spectral sequence together with all their differentials are
invariant (up to translation) with respect to hamiltonian isotopy.
Moreover, b. implies that Floer homology is isomorphic to
$EF^{2}_{\ast,0}(L,L';\eta)$ and so our invariant extends Floer
homology. It is therefore natural to expect to be able to estimate
the invariant $EF(L,L';\eta)$ when $L'$ is hamiltonian-isotopic to
$L$ (and $\eta=\eta_{0}$) in terms of some algebraic-topological
invariant of $L$. The fact that this invariant is precisely the Serre
spectral sequence of $\Omega L\to PL \to L$ is remarkable because,
due to the fact that $PL$ is contractible, this last spectral
sequence always contains non-trivial differentials. As we shall see
this trivial algebraic-topological observation together with the
geometric interpretation of the differentials discussed in
\S\ref{subsubsec:geom_int} leads to interesting applications.

\subsection{Some applications.}
Here is an overview of some of the consequences discussed in the
paper. It should be pointed out that we focus in this paper only on
the applications which follow rather rapidly from the main result. We
intend to discuss others that are less immediate in later papers.

We shall only mention in this subsection applications that take place
in the case when $L$ and $L'$ are hamiltonian isotopic and so we make
here this assumption.

\subsubsection{Algebraic consequences.} Under the assumption at (\ref{eq:connectivity}), a first
consequence of the theorem is that, if $\mathcal{K}=
\bigcup_{x,y}\{\Phi_{x,y}(\overline{\mathcal{M}}(x,y))\}\subset
\Omega L$ and let $\widehat{\mathcal{K}}$ be the closure of $K$ with respect to
concatenation of loops, then the inclusion
$\widehat{\mathcal{K}}\stackrel{k}{\hookrightarrow} \Omega L$ is surjective in homology. An immediate consequence of this is as follows. Notice first
that the space $\mathcal{M}'$ maps injectively onto
a subspace $\tilde{\mathcal{M}}$ of $\mathcal{P}(L,L')$ via the map that associates to each pseudo-holomorphic strip $u:\R\times [0,1]\to M$ the path
$u(0,-)$. Let $e:\tilde{\mathcal{M}}\to L$ be defined by $e(u)=u(0,0)$.
We show that
\begin{equation}\label{eq:loop_surj}
 H_{\ast}(\Omega e):H_{\ast}(\Omega
\tilde{\mathcal{M}}; \Z/2)\to H_{\ast}(\Omega L; \Z/2)\ {\rm is\
surjective}~.~
\end{equation}
This complements a result obtained by Hofer \cite{Ho} and
independently by Floer \cite{Fl} which claims that $H_{\ast}(e)$ is also surjective.

Another easy consequence is
that for a generic class of choices of $L'$, the image of the group
homomorphism $\Pi=\omega| : \pi_{2}(M,L\cup L')\to \R$ verifies
\begin{equation}\label{eq:rk_pi} {\rm rk\/}(Im(\Pi))\geq
\sum_{i} {\rm dim\/}_{\Z/2}H_{i}(L;\Z/2)-1~.~
\end{equation}

\subsubsection{Existence of pseudo-holomorphic ``strips".}\label{subsubsec:intro_exist} A rather
immediate consequence of the construction of $EF(L,L')$ is that
through each point in $L\backslash L'$ passes at least one strip
$u\in\mathcal{M}'$ of Maslov index at most $n$. By appropriately
refining this argument we shall see that we may even bound the energy of these
strips which ``fill" $L$
by the energy of a hamiltonian diffeomorphism that carries $L$ to $L'$.
More precisely, denote by $||\phi||_{H}$ the Hofer norm (or energy - see \cite{Ho2} and equation (\ref{eq:hofer_norm})) of a Hamiltonian
diffeomorphism $\phi$. We put (as in \cite{Chek} and \cite{Oh}):
$$\nabla(L,L')=\inf_{\psi\in\mathcal{H}, \psi(L)=L'}||\psi||_{H}$$ where
$\mathcal{H}$ is the group of compactly supported hamiltonian
diffeomorphisms.
We prove that \emph{through each point of $L\backslash L'$ passes a pseudo-holomorphic
strip which is of Maslov index at most $n$ and whose symplectic area is at most
$\nabla(L,L')$}. This fact has many interesting geometric consequences.
We describe a few in the next paragraph.

\subsubsection{Non-squeezing and Hofer's energy.}\label{subsubsec:non_squeeze}
Consider on $M$ the riemannian metric induced by some fixed generic almost
complex structure which tames $\omega$. The areas below are defined with
respect to this metric. For two points $x,y\in L\cap L'$ let

\begin{equation}\label{eq:strips}
\begin{matrix}\mathcal{S}(x,y)=\{ u\in C^{\infty}([0,1]\times
[0,1], M) \ :& u([0,1],0)\subset L, u([0,1],1)\subset L', \\
\hspace{55pt} u(0,[0,1])=x ,u(1,[0,1])=y\}~.~ &
\end{matrix}\end{equation}
Fix the notation: $$a_{L,L'}(x,y)=\inf\{area(u) : u\in
\mathcal{S}(x,y)\}~.~$$
Let $a_{k}(L,L')=\min\{a_{L,L'}(x,y) : x,y\in L\cap L' \ , \ \mu(x,y)=k\}$ and,
similarly, let $A_{k}(L,L')$ be the maximum of all $a_{L,L'}(x,y)$ where $x,y\in L\cap L'$ verify $\mu(x,y)=k$.

We prove that:
$$
a_{n}(L,L') \leq \nabla(L,L')\ ~.~
$$

For $x\in L\backslash L'$ let $\delta(x)\in [0,\infty)$ be the
maximal radius $r$ of a standard symplectic ball $B(r)$ such that
there is a symplectic embedding $e_{x,r}:B(r)\to M$ with
$e_{x,r}(0)=x$, $e_{x,r}^{-1}(L)= B(r)\cap \R^{n}$ and
$e_{x,r}(B(r))\cap L'=\emptyset$. We thank Fran\c{c}ois Lalonde who
noticed that, as we shall see, $\delta_{x}$ does not depend of $x$.
Therefore, we introduce the \emph{ball separation energy} between $L$ and $L'$
by
$$\delta(L,L')=\delta_{x}~.~$$ We show a second inequality

\begin{equation}\label{eq:squeeze}
\frac{\pi}{2}\ \delta(L,L')^{2}\leq A_{n}(L,L')
\end{equation}
and also that:

\begin{equation}\label{eq:squee_Hofer}
\frac{\pi}{2}\ \delta(L,L')^{2}\leq \nabla(L,L')~.~
\end{equation}

The results summarized in \S\ref{subsubsec:intro_exist} as well as the inequalitities
(\ref{eq:squeeze}) and (\ref{eq:squee_Hofer}) are first proved under the assumption
at (\ref{eq:connectivity}). However, we then show that our spectral sequence
may also be constructed (with minor modifications) when $L$ and $L'$ are hamiltonian isotopic under the single additional assumption $\omega|_{\pi_{2}(M,L)}=0$
and as a consequence these three results also remain true in this setting.

The inequality (\ref{eq:squee_Hofer}) is quite powerful:
it implies that $\nabla(-,-)$ (which is easily seen to be symmetric and to satisfy the triangle inequality) is also non-degenerate thus reproving - when $\omega|_{\pi_{2}(M,L)}=0$ - a result of Chekanov \cite{Chek}.
The same inequality is of course reminiscent the known
displacement-energy estimate in
\cite{Lalonde_McDuff} and, indeed, this estimate easily follows from
(\ref{eq:squee_Hofer})  (of course, under the assumption $\omega|_{\pi_{2}(M)}=0$) by applying this inequality to the diagonal embedding $M\to M\times M$.

\subsection{The structure of the paper.} In \S \ref{sec:spec} we start
by recalling the basic notations and conventions used in the paper
as well as the elements from Floer's theory that we shall need. We
then pass to the main task of the section which is to present the
construction of $EF(L,L';\eta)$. A key technical ingredient in this
construction is the fact that the compactifications of the moduli
spaces of Floer trajectories, $\overline{\mathcal{M}}(x,y)$, have a
structure of manifolds with corners. This property is closely related
to the gluing properties proven by Floer in his classical paper
\cite{Fl2} and is quite similar to more recent results proven by
Sikorav in \cite{Siko1}. In fact, this same property also appears to
be a feature of the Kuranishi structures used by Fukaya and Ono in
\cite{Fuk}. For the sake of completeness we include a complete proof
of the existence of the manifold-with-corners structure in the
Appendix. We then verify the points a., b., c. of Theorem
\ref{theo:main}. In \S \ref{sec:serre} we prove point d. of Theorem
\ref{theo:main}. This proof is based on one hand on the classical
method of comparing the Floer complex to a Morse complex of a Morse
function on $L$ and, on the other hand, on a new Morse theoretic result
which shows that if in the construction of $EF(L,L')$ the moduli
spaces of pseudo-holomorphic curves are replaced with moduli spaces
of negative gradient flow lines, then the resulting spectral sequence
is the Serre spectral sequence of the statement. The whole
construction of $EF(L,L')$ has been inspired by precisely this
Morse theoretic result which, in its turn, is a natural but
non-trivial extension of some ideas described in \cite{Co2} and
\cite{Co4}.

Finally, \S \ref{sec:appli} contains the applications mentioned above
as well as various other comments.

\subsection*{Acknowledgements}
We thank Fran\c{c}ois Lalonde for useful discussions.
The second author is grateful to the organizers of the IAS/Park City
summer-school in 1997 and of the Fields/CRM semester in the Spring of 2001
for encouraging his presence at these meetings and thus easing
his introduction to symplectic topology. We both thank the organizers of
the Oberwolfach Arbeitsgemeinschaft in October 2001 during which this
project has started. We thank Katrin Wehrheim for pointing out some
imprecisions in the appendix.

\section{The spectral sequence.}\label{sec:spec}

It turns out that it is more natural to construct a richer invariant than
the one appearing in Theorem \ref{theo:main}. The spectral sequence of the
theorem will be deduced as a particular case of this construction.

As before let $L,L'$ be closed lagrangian submanifolds of the fixed
symplectic manifold $(M,\omega)$. In this section we assume that their intersection
is transversal and that
$\omega|_{\pi_{2}(M)}=c_{1}|_{\pi_{2}(M)}
=0=\pi_{1}(L)=\pi_{1}(L')$. As $\pi_{2}(M)\to \pi_{2}(M,L)$ is surjective
(and similarly for $L'$) we deduce $\omega|_{\pi_{2}(M,L)}=\omega|_{\pi_{2}(M,L')}=0$.

\subsection{Recalls and Notation.}\label{subsec:recall}
We start by recalling some elements from Floer's construction.
This machinery has now been described in detail in various sources (for example,\cite{Fl2},\cite{Oh2})
so that we shall only give here a very brief presentation.

We fix a path $\eta\in
\mathcal{P}(L,L')=\{\gamma\in C^{\infty}([0,1], M) : \gamma(0)\in L$,
$\gamma(1)\in L'\}$ and let $\mathcal{P}_{\eta}(L,L')$ be the
path-component of $\mathcal{P}(L,L')$ containing $\eta$. We also fix
an almost complex structure $J$ on $M$ that tames $\omega$ in the
sense that the bilinear form $X,Y\to \omega(X,JY)=\alpha (X,Y)$ is a
Riemannian metric. The set of all the almost complex structures on $M$
that tame $\omega$ will be denoted by $\mathcal{J}_{\omega}$.
Moreover, we also consider a smooth
Hamiltonian $H:[0,1]\times M\to \R$ and its associated
family of hamiltonian vector fields $X_{H}$ determined
by the equation $$\omega(X^{t}_{H},Y)=-dH_{t}(Y) \ , \ \forall Y$$
as well as the hamiltonian isotopy $\phi^{H}_{t}$ given by
\begin{equation}\label{eq:hamiltonian}
\frac{d}{dt}\phi^{H}_{t}=X_{H}^{t}\circ \phi^{H}_{t} \ , \phi^{H}_{0}=id~.~
\end{equation}
The gradient of $H$, $\nabla H$, is computed with respect to $\alpha$ and it
verifies $J\nabla H=X_{H}$.

We shall also assume that $\phi^{H}_{1}(L)$ intersects transversely $L'$.
Moreover, $H$ and all the hamiltonians considered in this paper are assumed
to be constant outside of a compact set.

\subsubsection{The action functional and pseudo-holomorphic strips.}
\label{subsubsec:action_funct} The idea behind the whole construction
is to consider the action functional
\begin{equation}\label{eq:action}
\mathcal{A}_{L,L',H}:\mathcal{P}_{\eta}(L,L')\to \R \ , \ x\to
-\int \overline{x}^{\ast}\omega +\int_{0}^{1}H(t,x(t))dt
\end{equation}
where $\overline{x}(s,t):[0,1]\times [0,1]\to M$ is such that
$\overline{x}(0,t)=\eta(t)$, $\overline{x}(1,t)=x(t)$, $\forall t\in
[0,1]$, $x([0,1],0)\subset L$, $x([0,1],1)\subset L'$. The fact that
$L$ and $L'$ are simply connected Lagrangians and $\omega$ vanishes
on $\pi_{2}(M)$ implies that $\mathcal{A}_{L,L',H}$ is well-defined.
To shorten notation we neglect the subscripts $L,L',H$ in case no confusion
is possible. We shall also assume $\mathcal{A}(\eta)=0$ (this is of course
not restrictive).

Given a vector field $\xi$ tangent to $TM$ along $x\in
\mathcal{P}(L,L')$ we derive $\mathcal{A}$ along $\xi$ thus
getting
\begin{eqnarray}\label{eq:derivative}d\mathcal{A}(\xi)&=& -\int_{0}^{1}
\omega(\xi,\frac{dx}{dt})dt +\int_{0}^{1}dH_{t}(\xi)(x(t))dt=\\
&=& \int_{0}^{1}\alpha(\xi,J\frac{dx}{dt}+\nabla
H(t,x))dt~.~\nonumber
\end{eqnarray}
This means that the critical points of
$\mathcal{A}_{L,L'}$ are precisely the orbits of $X_{H}$ that start on $L$,
end on $L'$ and which belong to $\mathcal{P}_{\eta}(L,L')$. Obviously,
these orbits are in bijection with a subset of $\phi_{1}^{H}(L)\cap L'$ so they are finite
in number. A particular important case is when $H$ is constant. Then these orbits coincide
with the intersection points of $L$ and $L'$ which are in the class of
$\eta$. We denote the set of these orbits by $I(L,L'; \eta, H)$. In case $H$ is constant
we shall also use the more intuitive notation $L\cap_{\eta}L'$.

The putative
associated equation for  the negative $L^{2}$-gradient flow lines of
$\mathcal{A}$ has been at the center of Floer's work and is:
\begin{equation}\label{eq:floer}
\frac{\partial u}{\partial s}+J(u)\frac{\partial u}{\partial t}+\nabla H(t,u)=0
\end{equation}
with
 $$ u(s,t):\R\times [0,1]\to M \ , u(\R,0)
 \subset L \ , \ u(\R,1)\subset L' ~.~$$
When $H$ is constant, the solutions of (\ref{eq:floer}) are called
\emph{pseudo-holomorphic strips}. They coincide with the zeros of the
operator $\overline{\partial}_{J}=\frac{1}{2}(d+ J\circ d\circ
\mathbf{i})$. It is
well known that (\ref{eq:floer}) does not define a flow in any
convenient sense.

Let $\mathcal{S}(L,L')=\{u\in C^{\infty}(\R\times [0,1],M) :
u(\R,0)\subset L \ , \ u(\R,1)\subset L' \}$ and for
$u\in\mathcal{S}(L,L')$ consider the energy
\begin{equation}\label{eq:energy}
E_{L,L',H}(u)=\frac{1}{2}\int_{\R\times [0,1]} ||\frac{\partial
u}{\partial s}||^{2}+ ||\frac{\partial u}{\partial t}-X^{t}_{H}(u)||^{2}\ ds\ dt
~.~
\end{equation}
The key point of the whole theory is that, for a generic choice of
$J$, the solutions $u$ of (\ref{eq:floer}) which {\em are of finite
energy}, $E_{L,L',H}(u)<\infty$, do behave very much like (negative)
flow lines of a Morse-Smale function when viewed as elements in
$C^{\infty}(\R,\mathcal{P}_{\eta}(L,L'))$ (in particular,
$\mathcal{A}$ is decreasing along such solutions) . The type
of genericity needed here sometimes require that $J$ be
time-dependent. In other words $J=J_{t}$, $t\in [0,1]$ is a
one-parameter family of almost complex structures each taming
$\omega$. In this case the equation (\ref{eq:floer}) is understood as
$\partial u/\partial s +J_{t}(u)\partial u/\partial t+\nabla H(t,u)=0$. We now put
\begin{equation}\label{eq:param_moduli}
\mathcal{M}'=\{u\in\mathcal{S}(L,L') : u \ {\rm verifies
(\ref{eq:floer}) } \ , \ E_{L,L',H}(u)<\infty \}~.~
\end{equation}
The translation $u(s,t)\to u(s+k,t)$ obviously induces an $\R$ action
on $\mathcal{M}'$ and we let $\mathcal{M}$ be the quotient space. An
important feature of $\mathcal{M}'$ is that for each $u\in
\mathcal{M}'$ there exist $x,y\in I(L,L';\eta, H)$ such that the (uniform)
limits verify
\begin{equation}\label{eq:ends}
\lim_{s\to-\infty}u(s,t)=x(t) \ , \ \lim_{s\to+\infty}u(s,t)=y(t)~.~
\end{equation}
We let $\mathcal{M}'(x,y)=\{u\in \mathcal{M}' : u \ {\rm verifies} \
(\ref{eq:ends})\}$ and $\mathcal{M}(x,y)=\mathcal{M}'(x,y)/\R$.
Therefore, $\mathcal{M}=\bigcup_{x,y} \mathcal{M}(x,y)$. In case we need
to explicitely indicate to which pair of Lagrangians, to what Hamiltonian
and to what almost complex structure are associated these moduli spaces we shall add $L$ and $L'$, $H$, $J$ as subscripts
(for example, we may write $\mathcal{M}_{L,L',H,J}(x,y)$).

\subsubsection{Dimension of $\mathcal{M}(x,y)$ and the Maslov index.}
\label{subsubsec:Maslov_ind} Let $\mathcal{L}(n)$ be the set of
lagrangian subspaces in $(\R^{2n},\omega_{0})$. It is well-known that
$H^{1}(\mathcal{L}(n);\Z)\approx \Z$ has a generator given by a
morphism called the Maslov index $\mu:\mathcal{L}(n)\to S^{1}$
(geometrically it is given as the class dual to the Maslov cycle
constituted by the lagrangian subspaces non-transversal to the
vertical lagrangian). For $x,y\in I(L,L';\eta,H)$ we let
(as in (\ref{eq:strips}))
$$\begin{matrix}\mathcal{S}(x,y)=\{ u\in C^{\infty}([0,1]\times
[0,1], M) \ :& u([0,1],0)\subset L, u([0,1],1)\subset L', \\
\hspace{55pt} u(0,t)=x(t) ,u(1,t)=y(t)\}~.~ &
\end{matrix}$$
 and suppose that $u\in\mathcal{S}(x,y)$. Following
the work of Viterbo \cite{Viterbo}, the Maslov index of $u$,
$\mu(u)$, is given as the degree of the map $S^{1}=\partial
([0,1]\times
[0,1])\stackrel{\gamma}{\longrightarrow}\mathcal{L}(n)\stackrel{\mu}{\longrightarrow}
S^{1}$ with the loop $\gamma$ defined as follows. First notice that
$u^{\ast}TM$ is a trivial symplectic bundle (and all trivializations
are homotopic). We fix such a trivialization. This allows the
identification of each space $T_{x}L\subset T_{x}M$ to an element of
$\mathcal{L}(n)$ (and similarly for $T_{x}L'$). We then define the
loop $\gamma :S^{1}\to \mathcal{L}(n)$ as follows.
We let $\gamma_{0}$ be the path of Lagrangians $(\phi_{t}^{H})^{-1}_{\ast}T_{x(1)}L'$
and we let $\gamma_{1}$ be the path $(\phi_{t}^{H})^{-1}_{\ast}T_{y(1)}L'$.
We then join $(\phi^{H}_{1})^{-1}_{\ast}T_{x(1)}L'$ to
$(\phi^{H}_{1})^{-1}_{\ast}T_{y(1)}L'$ by a path of lagrangian subspaces $\gamma'(t)\subset
T_{u(t,0)}M$ such that for each $t$, $\gamma'(t)$ is transversal to
$T_{u(t,0)}L$ and let $\gamma=\gamma_{0}\ast\gamma'\ast\gamma_{1}^{-1}\ast\gamma''$ where
$\gamma''(t)$ is the path $t\to T_{u(1-t,1)}L'$. It is easy to see
that such a path $\gamma'$ does exist and that the degree of the
composition $\mu\circ\gamma$ is independent of the choice of
$\gamma'$ as well as of that of the trivialization. Given that $L$
and $L'$ are simply connected and $c_{1}|_{\pi_{2}(M)}=0$ we see
that for any $u,v\in \mathcal{S}(x,y)$ we have
$\mu(u)=\mu(v)$. Therefore, for any $x,y\in
I(L,L';\eta,H)$ we may define
$$\mu(x,y)=\mu(u)\ , \ u\in \mathcal{S}(x,y)~.~$$
This implies that, in this case, for any three points $x,y,z\in
I(L,L';\eta,H)$ we have
\begin{equation}\label{eq:transit}
\mu(x,z)=\mu(x,y)+\mu(y,z)~.~
\end{equation}
The fundamental role of the Maslov index in relation to the
properties of the action functional is provided by the fact that the
linearized operator $D^{H,J}_{u}$ associated to
the operator $\overline{\partial}_{J}+(1/2)\nabla H$ at $u$ is Fredholm
of index $\mu(u)$. In case $J$ is such that  $D^{H,J}_{u}$ is
surjective for all $u\in \mathcal{M}'(x,y)$ and all $x,y\in
I(L,L';\eta,H)$ (see \ref{sec:appendix}), it follows that the spaces
$\mathcal{M}'(x,y)$ are smooth manifolds (generally non-compact) of
dimension $\mu(x,y)$. Under certain circumstances the theory works in
the same way even if $L$ and $L'$ are non-transversal (for example if
$L=L'$) but in that case the choice of $H$ needs to be generic.
In all cases, we shall call a pair $(H,J)$ regular if the
surjectivity condition mentioned above is satisfied. In our setting it is
easy to see that for any $x\in I(L,L';\eta, H)$,
the space $\mathcal{M}'(x,x)$ is reduced to the constant
solution equal to $x$. Because of that we will always assume in the paper
that in the writing $\mathcal{M}(x,y)$ we have $x\not=y$. Thus,
$\mathcal{M}(x,y)$ is also a smooth manifold whose dimension is $\mu(x,y)-1$.
The set of regular $(H,J)$'s is generic and we assume below that we are
using such a pair.

\subsubsection{Naturality of Floer's equation.}\label{subsubsec:naturality}
Let $L''=(\phi_{1}^{H})^{-1}(L')$. Consider the map
$b_{H}: \mathcal{P}(L,L'')\to \mathcal{P}(L,L')$ defined by
$(b_{H}(x))(t)=\phi_{t}^{H}(x(t))$. Let $\eta'\in\mathcal{P}(L,L'')$ be
such that $\eta=b_{H}(\eta')$. Clearly, $b_{H}$ restricts to a map
between $\mathcal{P}_{\eta'}(L,L'')$ and $\mathcal{P}_{\eta}(L,L')$ and,
moreover, by our assumption on $\phi$, the intersection of $L$
and $L''$ is transverse and the same map restricts to a
bijection $L\cap_{\eta'} L''=I(L,L'';\eta',0)\to I(L,L';\eta, H)$.

We also have $$\mathcal{A}_{L,L',H}(b_{H}(x))=\mathcal{A}_{L,L'',0}(x)~.~$$
Indeed,
let $\overline{x}(s,t):[0,1]\times [0,1]\to M$ be such that
$\overline{x}(0,t)=\eta'(t)$, $\overline{x}(1,t)=x(t)$, $\forall t\in
[0,1]$, $x([0,1],0)\subset L$, $x([0,1],1)\subset L''$ and let $\tilde{x}(s,t)=\phi_{t}^{H}(x(s,t))$.
We then have (by using
(\ref{eq:hamiltonian}) and leting $\phi=\phi^{H}$):
\begin{eqnarray*}\int_{[0,1]\times [0,1]}\tilde{x}^{\ast}\omega &= &
\int_{[0,1]\times [0,1]} \overline{x}^{\ast}(\phi^{\ast}\omega)+
\int_{0}^{1}\int_{0}^{1}
\omega(\frac{\partial \tilde{x}}{\partial s},\frac{\partial \phi}{\partial t})dsdt =\\
=\int_{[0,1]\times [0,1]}\overline{x}^{\ast}\omega &+&
\int_{0}^{1}(\int_{0}^{1}dH(\frac{\partial\tilde{x}}{\partial s})ds)dt =
-\mathcal{A}_{L,L'',0}(x)+\int_{0}^{1}H(b_{H}(x)(t))dt ~.~
\end{eqnarray*}

Moreover, the map $b_{H}$ does identify the geometry of the two
action functionals. This is due to the fact that for $u:\R\times
[0,1]\to M$ with $u(\R,0)\subset L$, $u(\R,1)\subset L''$,
$\tilde{u}(s,t)=\phi_{t}(u(s,t))$, $\tilde{J}=\phi^{\ast}J$ we have
$$\phi_{\ast}(\frac{\partial u}{\partial s}+ \tilde{J} \frac{\partial u}{\partial t}) =
\frac{\partial \tilde{u}}{\partial s}+J (\frac{\partial
\tilde{u}}{\partial t}-X_{H})~.~$$ Therefore, the map $b_{H}$ induces
diffeomorphisms (that we shall denote by the same symbol):
$$b_{H}:\mathcal{M}_{L,L'',\tilde{J},0}(x,y)\to \mathcal{M}_{L,L',J,H}(x,y)$$
where we have identified $x,y\in L\cap_{\eta'} L''$ with their orbits $\phi^{H}_{t}(x)$
and $\phi^{H}_{t}(y)$.

\subsubsection{Gromov compactification of $\mathcal{M}(x,y)$.}\label{subsubsec:gromov_comp}
The non-compactness of $\mathcal{M}(x,y)$ for $x,y\in I(L,L';\eta,H)$ is only due
to the fact that, as in the Morse-Smale case, a sequence of strips $u_{n}\in \mathcal{M}(x,y)$
might ``converge" to a broken strip. For example, it might converge to an element of
$\mathcal{M}(x,z)\times\mathcal{M}(z,y)$ for some other $z\in I(L,L';\eta, H)$. The type of
convergence used here has been studied extensively and it is called Gromov convergence.
Moreover, there are natural compactifications of the moduli spaces $\mathcal{M}(x,y)$
called Gromov compactification and denoted by $\overline{\mathcal{M}}(x,y)$ so that each of
the spaces $\overline{\mathcal{M}}(x,y)$ is a manifold with boundary and there is
a homeomorphism:
\begin{equation}\label{eq:Grom_comp}
\partial\overline{\mathcal{M}}(x,y)=\bigcup_{z\in I(L,L';\eta,H)}\ \overline{\mathcal{M}}(x,z)
\times \overline{\mathcal{M}}(z,y)~.~
\end{equation}

It is shown in the Appendix \ref{sec:appendix}, that
the manifolds $\overline{\mathcal{M}}(x,y)$ are manifolds with corners.
We insist there mainly on the homogenous case, when $H=0$. However, as
the maps $b_{H}$ constructed in \S \ref{subsubsec:naturality} are compatible with
equation (\ref{eq:Grom_comp}) this result is also true for a general $H$.

\subsection{Construction of the spectral sequence.}\label{subsec:constr_SS}

\subsubsection{Deformed pseudo-holomorphic strips viewed as paths.}
\label{subsubsec:defor_pseudo}
To each element $u\in \mathcal{M}'(x,y)$ we associate a continuous
path
\begin{equation}\label{eq:def_path}
\gamma_{u}:[0, \mathcal{A}(x)-\mathcal{A}(y)]\to
\mathcal{P}_{\eta}(L,L')
\end{equation} in a rather obvious way:
$\gamma_{u}(\mathcal{A}(x)-\mathcal{A}(y))=y$,
$\gamma_{u}(0)=x$ and  for $\tau \in (0,
\mathcal{A}(x)-\mathcal{A}(y))$,
$\gamma_{u}(\tau)=u(h_{u}(-\tau),[0,1])$ where
$$h_{u}:(\mathcal{A}(y)-\mathcal{A}(x),0)\to \R$$ is
defined by
$\mathcal{A}(u(h_{u}(\tau),[0,1]))=\tau+\mathcal{A}(x)$.
In short, $\gamma_{u}$ associates to $\tau$ the unique element of
$\mathcal{P}(L,L')$ which is of the form $u(\xi,-):[0,1]\to M$
for some $\xi\in\R$ and on which $\mathcal{A}$ has the value
$\mathcal{A}(x)-\tau$. The
function $h_{u}$ is well defined because $\mathcal{A}$ is
strictly decreasing along $u$ and it is easy to see that $\gamma_{u}$
is continuous (we shall use here the compact-open $C^{0}$-topology on
$\mathcal{P}(L,L')$). Obviously, $\gamma_{u}$ only depends of the
class of $u$ in $\mathcal{M}(x,y)$ and thus we have a map:
$$\gamma_{x,y}:\mathcal{M}(x,y)\to
C^{0}([0, \mathcal{A}(x)-\mathcal{A}(y)],
\mathcal{P}_{\eta}(L,L'))\ , \ \gamma_{x,y}(u)=\gamma_{u} ~.~$$ To
simplify notation let $$C_{x,y}\mathcal{P}=
C^{0}([0,\mathcal{A}(x)-\mathcal{A}(y)],
\mathcal{P}_{\eta}(L,L'))$$ which is taken to be void in case
$\mathcal{A}(x)\leq\mathcal{A}(y)$. The map
$\gamma_{x,y}$ is easily seen to be continuous in view of the
description of the charts of $\mathcal{M}(x,y)$. Moreover, in view of
the definition of Gromov compactness (or by using the description of
the small neighbourhoods of broken Floer orbits given in the Appendix
\ref{sec:appendix}) we see that this map extends to a continuous map
$$\overline{\gamma}_{x,y}:\overline{\mathcal{M}}(x,y)\to
C_{x,y}\mathcal{P}~.~$$ Notice that there exists an obvious
continuous composition map given by concatenation of paths
\begin{equation}\label{eq:concat}
\#: C_{x,y}\mathcal{P}\times C_{y,z}\mathcal{P}\to C_{x,z}\mathcal{P}
\end{equation}
which is associative in the obvious sense. As an immediate
consequence of the proof of (\ref{eq:Grom_comp}) we also see that for
each element $u=(u_{1},u_{2},\ldots u_{k})\in
\mathcal{M}(x,z_{1})\times\mathcal{M}(z_{1},z_{2})
\times\ldots\times\mathcal{M}(z_{k-1},y)
\subset\partial\overline{\mathcal{M}}(x,y)$ we have:
\begin{equation}\label{eq:comm_eval_concat}
\overline{\gamma}_{x,y}(u)=
\gamma_{x,z_{1}}(u_{1})\#\gamma_{z_{1},z_{2}}(u_{2})\#\ldots\#
\gamma_{z_{k-1},y}(u_{k})
\end{equation}

\subsubsection{Some additional path spaces.}\label{subsubsec:add_path_sp}
We fix here some more notation. Let $w$ be a path (homeomorphic to [0,1])
 embedded in $L$ that
joins all points $\{x(0) : x\in I(L,L';\eta, H) \}$ and let $\tilde{M},\tilde{L}$ be
respectively the quotient topological spaces obtained by contracting
$w$ to a point. Obviously, the quotient maps $M\to \tilde{M}$, $L\to
\tilde{L}$ are homotopy equivalences. We also have homotopy
equivalences $\mathcal{P}(L,L')\to \mathcal{P}(\tilde{L},L')$,
$\mathcal{P}_{\eta}(L,L')\to \mathcal{P}_{\eta}(\tilde{L},L')$.  We
shall denote any of these quotient maps by $q$.
We shall also need the obvious map
$l:\mathcal{P}_{\eta}(\tilde{L},L')\to \tilde{L}, \ \  l(\gamma)=\gamma(0)$.
Notice that the spaces
$\tilde{L},\tilde{M},\mathcal{P}(\tilde{L},L')$ have a distinguished
base point, $\ast$, given by the class of the path $w$ and
$(l\circ q)(I(L,L';\eta,H))=\ast$.

For any pointed topological space $X$ we recall
that $\Omega X$ is the space of continuous loops in $X$ that are
based at the distinguished point of $X$ and are parametrized by the
interval $[0,1]$. This space is homotopy equivalent to the space of
Moore loops on $X$, $\Omega' X$, which consists of the continuous
loops in $X$ that are parametrized by arbitrary intervals $[0,a]$,
$a\in [0,\infty)$ (and, again, are based at the distinguished point
of $X$).

The compositions $l\circ q$
induce maps
$$\mathcal{Q}_{x,y}:C_{x,y}\mathcal{P}\to
\Omega'\tilde{L}, \
(\mathcal{Q}_{x,y}(a))(\tau)=(l\circ q)(a(\tau))~.~$$

Concatenation
of loops gives Moore loops the structure of topological monoid. This
operation, denoted by $\cdot$, commutes in an obvious way with the
maps $\mathcal{Q}_{-,-}$ and the operation $\#$ of (\ref{eq:concat}).

Fix also the notation
\begin{equation}\label{eq:eval}
\Phi_{x,y}=\mathcal{Q}_{x,y}\circ
\overline{\gamma}_{x,y}~.~
\end{equation}

For further use, notice that the
space $\mathcal{P}(L,L')$ (and therefore also
$\mathcal{P}(\tilde{L},L')$) is homotopy equivalent to the homotopy
pull-back of the two inclusions $L\hookrightarrow M$ and
$L'\hookrightarrow M$.

\subsubsection{An algebraic construction.}\label{subsubsec:alg_constr}
For a topological space $X$ let $S_{\ast}(X)$ be the $\Z/2$-cubical (normalized)
 chain complex of $X$ . We use cubical chains - that is chains whose domains are unit cubes (see \cite{Massey} for definitions) - instead of singular chains because in this case, for two spaces $X,Y$,
we have an obvious map $S_{k}(X)\times S_{q}(Y)\to
S_{k+q}(X\times Y)$ defined by $(\sigma\times \sigma')(x,y)=(\sigma(x),\sigma'(y))$.
Moreover, the multiplication $\cdot$ directly induces a natural
multiplication denoted again by $\cdot:S_{k}(\Omega' X)\otimes S_{l}(\Omega' X)\to S_{k+l}(\Omega' X)$
defined by $(\sigma\cdot \sigma')(x,y)=\sigma(x)\cdot\sigma(y)$ where $x\in [0,1]^{k}$,
$y\in [0,1]^{l}$.

In particular, this turns $S_{\ast}(\Omega' \tilde{L})$
into a differential ring that we shall denote from now on by  $\mathcal{R}_{\ast}$.

\begin{defi}\label{def:repres} A representing
chain system for the moduli spaces associated to $L,L',J,H,\eta$ is
a family $\{s_{xy}\in S_{\mu(x,y)-1}(\overline{\mathcal{M}}(x,y)) : x,y\in I(L,L';\eta,H)\}$
such that:
\begin{itemize}
\item[i.] The image of $s_{xy}$ in $S_{\ast}
(\overline{\mathcal{M}}(x,y),\partial \overline{\mathcal{M}}(x,y))$ is a cycle representing
the fundamental class.
\item[ii.] With the identifications given by equation (\ref{eq:Grom_comp}) we have
 $\partial s_{xy}=\sum_{z} s_{xz}\times s_{zy} \in S_{\ast}(\overline{\mathcal{M}}(x,y))$.
\end{itemize}
\end{defi}

\begin{lem}\label{lem:repres_chain_syst}
With the assumptions and notations above, there exists a representing
chain system for the moduli spaces
$\mathcal{M}_{L,L',J,H,\eta}(-,-)$.
\end{lem}

\begin{proof}We construct the $s_{xy}$'s by induction. Assume the construction
accomplished for $\mu(x,y)-1<k$. Consider now a pair $x,y$ with
$\mu(x,y)-1=k$. We may assume that $\mathcal{M}(x,y)$ is connected
(if not we apply the the argument below one component at a time).
Using the identifications in (\ref{eq:Grom_comp}) consider the chain
$c_{xy}=\sum_{z} s_{xz}\times s_{zy}\in S_{k-1}(\partial
\overline{\mathcal{M}}(x,y))$. We denote the differential in
$S_{\ast}(-)$ by $\partial$ and we compute $\partial
c_{xy}=\sum_{z}\partial s_{xz}\times s_{yz} +\sum s_{xz}\times
\partial s_{zy}= \sum_{z,k}(s_{xt}\times s_{tz})\times s_{yz} +
\sum_{z,j}s_{yz}\times (s_{zj}\times s_{jy})= 2(\sum_{s,r}
s_{xs}\times s_{sr}\times s_{rz})=0$. The homology class represented
by $c_{xy}$ is the fundamental class of $\partial
\overline{\mathcal{M}}(x,y)$. This is because the image of this class
in any one of $H_{k-1}(\overline{\mathcal{M}}(x,z)\times
\overline{\mathcal{M}}(z,y), \partial
(\overline{\mathcal{M}}(x,z)\times \overline{\mathcal{M}}(z,y)))$
coincides with the class represented by $s_{xz}\times s_{zy}$ which
is the fundamental class. Therefore, $c_{xy}\in Im(\partial :S_{k}(
\overline{\mathcal{M}}(x,y)) \to
S_{k-1}(\overline{\mathcal{M}}(x,y)))$. Let $s_{xy}$ be such that
$\partial s_{xy}=c_{xy}$. By construction, property ii. of a
representing system is then satisfied. The first property is also
satisfied because the image of $s_{xy}$ is a cycle in
$S_{k}(\overline{\mathcal{M}}(x,y)),\partial
\overline{\mathcal{M}}(x,y))$ and the homology connectant $\delta $
of the pair $(\overline{\mathcal{M}}(x,y)),\partial
\overline{\mathcal{M}}(x,y))$ is an isomorphism in dimension $k$ and
it verifies $\delta ([s_{xy}])=[c_{xy}]$.
\end{proof}

\begin{rem}\label{rem:caution_triang}{\rm
Representing chain systems appear naturally when the
compactified moduli spaces $\overline{\mathcal{M}}(x,y)$ are triangulated
(or rather ``cubulated") in a way compatible with formula (\ref{eq:Grom_comp}):
the $s_{xy}$'s may then be taken to be the sum of the top dimensional cubes.
However, the existence of such a triangulation is not obvious.
The most direct approach to constructing such a triangulation is to proceed
by induction. Assuming that a triangulation
of $\partial \overline{\mathcal{M}}(x,y)$ is constructed the induction
step is then to extend this triangulation to the whole of $\overline{\mathcal{M}}(x,y)$. For this extension to exist one needs to check that the Kirby-Siebenmann obstruction vanishes - fact which is not apriori clear.}
\end{rem}

We now fix  a representing chain system $\zeta=\{s_{xy}\}$ and we define
\begin{equation}\label{eq:coefficents}a_{xy}\in \mathcal{R}_{\mu(x,y)-1}\ , \
a_{xy}=\Phi_{x,y}(s_{xy})~.~
\end{equation}

Let $m$ be the number of elements of the set $I(L,L';\eta,H)$. Fix one
point $z_{0}\in I(L,L';\eta,H)$  and for each
$x\in I(L,L';\eta,H)$ let $\mu(x)=\mu(x,z_{0})$. In view of
(\ref{eq:transit}) the function $\mu(-)$ so defined only depends of $z_{0}$ up to
a translation by a constant. Let a strict ordering $\succ$ of the
set $I(L,L';\eta,H)$  be such that we have
$\mu(x)>\mu(y) \Rightarrow x\succ y$.

The main algebraic object that we shall be using is the matrix
\begin{equation}\label{eq:matrix}
A=(a_{xy})_{\{x,y\in I(L,L';\eta,H)\}}\in M_{m,m}(\mathcal{R}_{\ast})~.~
\end{equation}

\begin{rem}\label{rem:non-inv_matr} {\rm Of course, despite of our short notation
for $A$, this matrix depends
on $L$, $L'$, $H$, $\eta$, the choice of $J$ and of $\zeta$.}
\end{rem}

If $C$ is a matrix with coefficients in $\mathcal{R}_{\ast}$, then we let $\partial C$ be the matrix
whose coefficients are obtained by applying the differential $\partial$ of $\mathcal{R}_{\ast}$
to the coefficients of $C$.

The key property of $A$ is as follows.

\begin{prop}\label{prop:matrix_sq}
Under the assumptions above we have:
$$A^{2}=\partial A~.~$$
\end{prop}

\begin{proof} This is immediate from the construction of $A$ and from (\ref{eq:Grom_comp})
and (\ref{eq:comm_eval_concat}).
Indeed, we have the following sequence of equalities
\[\begin{array}{rl}
\partial a_{xy}\ =&\partial \Phi_{x,y}(s_{xy})=\Phi_{x,y}(\partial  s_{xy})=
\Phi_{x,y}(\sum_{z}s_{xz}\times s_{zy})=\\
=&\mathcal{Q}_{x,y}\circ\overline{\gamma}_{x,y} (\sum_{z}s_{xz}\times s_{zy})=
\mathcal{Q}_{x,y}\circ (\sum_{z} \overline{\gamma}_{x,z}(s_{xz})
\#\overline{\gamma}_{z,y}(s_{zy}))=\\
=&\sum_{z}(\mathcal{Q}_{x,z}\circ\overline{\gamma}_{x,z})(s_{xz})\cdot
(\mathcal{Q}_{z,y}\circ\overline{\gamma}_{z,y})(s_{zy})=\\
=&\sum_{z}\Phi_{x,z}(s_{xz})\cdot\Phi_{z,y}(s_{zy})=\\
=&\sum_{z}a_{xz}\cdot a_{zy}
\end{array}\]
which is valid for any $x,y\in I(L,L';\eta,H)$.
\end{proof}

\subsubsection{The spectral sequence.}\label{subsubsec: ss_def}
We first use the matrix $A$ to define an $\mathcal{R}_{\ast}$- chain complex
$$\mathcal{C}^{\eta,J,\zeta}(L,L';H)=
(\mathcal{C}_{\ast},d)$$
which should be thought of as an {\em extended Morse type chain complex}
as discussed in the introduction.

We consider the graded $\Z/2$-vector space $\Z/2<I(L,L';\eta,H)>$
where the grading is given by $|x|=\mu(x)\ , \ \forall x\in
I(L,L';\eta,H)$ (recall that the ``absolute" Maslov index function
$\mu : I(L,L';\eta,H)\to \Z$ from \S\ref{subsubsec:alg_constr}
depends on our choice of a fixed point $z_{0}\in I(L,L';\eta,H)$
only up to translation by an integer constant).

Now let $\mathcal{C}_{\ast}$ be equal to the left
$\mathcal{R}_{\ast}$-module $\mathcal{R}_{\ast}\otimes
\Z/2<I(L,L';\eta,H)>$. The module operation is so that for
$c\in\mathcal{R}_{\ast}$ and $a\otimes b\in \mathcal{C}$ we have
$c\cdot(a\otimes b)=(c\cdot a)\otimes b$.  The differential
$d:\mathcal{C}_{\ast}\to \mathcal{C}_{\ast -1}$ is the unique
$\mathcal{R}_{\ast}$-module derivation (in the sense that $d(a\otimes
b)=\partial a\otimes b+a\cdot db$) such that
$$d(x)=\sum_{y} a_{xy}\otimes y \ , \ \forall x\in I(L,L';\eta,H)~.~$$

\begin{cor}\label{cor:differential}
For $\mathcal{C}^{\eta,J,\zeta}(L,L'; H)=(\mathcal{C}_{\ast},d)$ defined
as above we have $d^{2}=0$.
\end{cor}

\begin{proof} For any $x\in I(L,L';\eta,H)$ we have: $d(d(x))= d(\sum_{y} a_{xy}\otimes y)= \\ =
\sum_{y}\partial a_{xy}\otimes y + \sum_{z,y}a_{xy}\cdot a_{yz}\otimes z=
\sum_{t}(\partial a_{xt} + \sum_{s}a_{xs}\cdot a_{st})\otimes t$\\
and all these last terms vanish in view of the equality in Proposition
\ref{prop:matrix_sq} (and because we work over $\Z/2$).
\end{proof}

Consider the spectral sequence which is associated to the filtration
of the complex $\mathcal{C}^{\eta,J,\zeta}(L,L';H)$ defined by:
$$F^{k}\mathcal{C}=\mathcal{R}_{\ast}\otimes \Z/2<x\in
I(L,L';\eta,H) : \mu(x)\leq k> ~.~$$ Clearly, this is a differential
filtration and thus it does indeed induce a spectral sequence which
we shall denote by $EF(L,L';\eta, H, J,
\zeta)=(EF_{pq}^{r}(L,L';\eta,J,H,\zeta),d^{r}_{F})$. We fix the
notation such that an element of bi-degree $(p,q)$ in the spectral
sequence is a class coming from an element in $\mathcal{R}_{q}\otimes
\Z/2< x : \mu(x)=p>$ (this last vector space being isomorphic to
$EF_{pq}^{0}(L,L'; H)$). We shall sometimes omit $\eta$, $J$, $\zeta$
in the notation for the spectral sequence.

We denote by $CF_{\ast}(L,L';H)$ the Floer chain complex associated to $\mathcal{A}_{L,L',H}$
and by $HF_{\ast}(L,L';H)$ the respective Floer homology. The relation of these
with our spectral sequence is as follows.

\begin{prop}\label{prop:simple_theo}
For the spectral sequence defined above we have
\begin{itemize}
\item[a.] $EF^{1}(L,L';H)\simeq CF_{\ast}(L,L';H)\otimes
H_{\ast}(\Omega L)$.
 \item[b.] $EF^{2}(L,L';H)\simeq HF_{\ast}(L,L';H)\otimes
H_{\ast}(\Omega L)$.
 \item[c.] If $d^{r}_{F}\not=0$, then there exist
$x,y\in I(L,L';\eta, H)$, $\mu(x,y)=r$, such that $\mathcal{M}(x,y)\not=0$.
\item[d.] For $r\geq 1$, $(EF^{r}_{pq}(L,L';H), d^{r}_{F})$ is a spectral
sequence of $H_{\ast}(\Omega L)$-modules.
\end{itemize}
\end{prop}

\begin{proof}
The only part of $d$ that counts for the first point is the
internal differential in $S_{\ast}(\Omega' \tilde{L})$. This expresses the $E^{1}$
term as desired.
The differential $d^{1}$ is horizontal and is generated
by the part of $d$ that connects orbits of relative
Maslov index equal to $1$. This is precisely the Floer (classical)
differential and thus implies the second point. The third point is
obvious as $d^{r}_{F}\not=0$ implies that there are some $x,y\in
I(L,L';\eta,H)$ such that $a_{xy}\not=0$ and $\mu(x,y)=r$. The
fourth point is a direct consequence of the fact that the
differential $d$ of $\mathcal{C}_{\ast}$ verifies $d(a\otimes b)=
\partial a\otimes b+a\cdot db$
\end{proof}

\begin{rem}\label{rem:transl}{\rm
 Notice that a different choice for $z_{0}$ only modifies the
resulting spectral sequence by a translation.}
\end{rem}

The spectral sequence of Theorem \ref{theo:main} consists of the
terms of order greater or equal than $1$ of $EF(L,L')=EF(L,L';0)$. In
particular, Proposition \ref{prop:simple_theo} implies the points b.
and c. of this theorem. We still need to prove the rest of the
theorem.

\begin{rem}\label{rem:simple}{\rm It is possible to modify the construction above in
such a way as to replace the ring $\mathcal{R}_{\ast}$ with the richer ring
$S_{\ast}(\Omega' \mathcal{P}_{\eta}(L,L'))$.  However, as $\mathcal{R}_{\ast}$
is sufficient for the applications discussed in this paper
we shall not pursue this extension here.}
\end{rem}

\subsection{Proof of the main theorem.  I: Invariance of the spectral sequence.}\label{subsubsec:inv} Our next aim is to prove the point a.
of Theorem \ref{theo:main}. As we shall see this point will follow
rapidly from the main result of this subsection which is shown in
\S\ref{subsubsec:var_ham} below.

\subsubsection{Variation of the Hamiltonian.}\label{subsubsec:var_ham} Assume
that with $L,L',\eta,H, J,\zeta$ as above we additionally have a
Hamiltonian $H':[0,1]\times M\to \R$ which is also constant outside of
a compact set.  We consider an almost complex structure $J'$
so that the pair $(H',J')$ is regular and so
$EF(L,L';\eta,J',H',\zeta')$ is defined with $\zeta'$ a representing
system of cochains for the moduli spaces associated to $L,L',J' H',\eta$. Let
\begin{equation}\label{eq:ham_distance}
\epsilon(L,L'; H, J)=\min\{ E_{L,L',H}(u) :
u\in\mathcal{M}'_{L,L',J,H}\} \end{equation}
(where $E_{L,L',H}$ is the energy as defined in (\ref{eq:energy})).

\begin{theo}\label{theo:comparison}
Under the assumptions above:
\begin{itemize}
\item[a.] There exists a chain morphism
$$\mathcal{V}:\mathcal{C}^{\eta,J,\zeta}(L,L';H)\to
\mathcal{C}^{\eta,J',\zeta'}(L,L';H')$$
of possibly non-zero degree which induces an isomorphism up to
translation between $EF^{r}(L,L';H)$ and $EF^{r}(L,L'; H')$ for $r\geq 2$.
\item[b.] If $||H'-H||_{0}<\epsilon(L,L';H,J)/4$, then there exists a morphism $\mathcal{V}$
as before which admits a retract.
\end{itemize}
\end{theo}

\begin{rem}\label{rem:rig} {\rm  A morphism of chain complexes $f:C_{\ast}\to D_{\ast+k}$ is said
to admit a retract if there exists another morphism $g:D_{\ast}\to C_{\ast-k}$ such that
$g\circ f=id_{C}$.
Clearly, if $\mathcal{V}$ admits a retract, then the same is true for the morphism induced
by $\mathcal{V}$ on each page of the spectral sequence.
Therefore, point b. of Theorem \ref{theo:comparison} shows, in particular,
that $EF^{r}(L,L';\eta,J,H,\zeta)$ does not depend on $J$ (or $\zeta$)  already for $r\geq 1$.}
\end{rem}

The idea  for the proof of Theorem \ref{theo:comparison} is classical
in Floer's theory : we adapt the previous construction to the case of
the moduli spaces of solutions of an equation similar to
(\ref{eq:floer}) but such as to allow for deformations from the
hamiltonian $H$  to the hamiltonian $H'$.

\begin{proof}
To shorten notation let $I=I(L,L';\eta, H)$, $I'=(L,L'; \eta, H')$.

We start with some recalls on Floer's comparison method. Take a
smooth homotopy $H^{01}:\R\times [0,1]\times M\to \R$ and
a homotopy $J^{01}: \R\times M\to End(TM)$,
$J^{01}_{s}\in\mathcal{J}, \forall s\in\R$  (here $\mathcal{J}$ is
the set of almost complex structures on $M$) such that there exists
$R>0$ with the property that, for $s\geq R$, we have
$(H^{01}_{s}(x),J^{01}_{s}(x))=(H(x),J(x))$ and for $s\leq -R$,
$(H^{01}_{s},J^{01}_{s})=(H'(x),J'(x))$, $\forall x\in M$. Moreover,
we assume that there exists a compact set such that for all $s\in
\R$, $H^{01}_{s}$ is constant outside this compact set. Consider the
equation:
\begin{equation}\label{eq:comp_fl}
\frac{\partial u}{\partial s}+J^{01}(s,u)\frac{\partial u}{\partial
t}+\nabla^{s}_{x}H^{01}(s,t,u)=0
\end{equation}
where $\nabla^{s}_{x}H^{01}(s,t,-)$ is the gradient of the function
$H^{01}(s,t,-)$ with respect to the riemannian metric induced by
$J^{01}_{s}$ and $u:\R\times [0,1]\to M$ with $u(\R,0)\subset L$ and
$u(\R, 1)\subset L'$. We may define the energy $E_{L,L',H^{01}}$ by
replacing $H$ in formula (\ref{eq:energy}) by $H^{01}$.

The finite energy solutions of (\ref{eq:comp_fl}) have properties
that are very similar to those of (\ref{eq:floer}). In particular,
for each such solution $u$ there exist $x\in I$, $y\in I'$ such that
\begin{equation}\label{eq:assymptotes}
\lim_{s\to-\infty}u(s,-)=x \ , \ \lim_{s\to\infty}u(s,-)=y~.~
\end{equation}
 If the linearized operator asociated to (\ref{eq:comp_fl}),
$D^{H^{01},J^{01}}_{u}$, is surjective for each finite energy
solution $u$ we say that the pair $(H^{01}, J^{01})$ is regular.
There is again a generic set of choices of regular such pairs. Again,
to insure genericity of regularity it might be needed to assume that
$J^{01}$ is also time dependent. We shall assume
from now on that $(H^{01},J^{01})$ is regular. We denote by
$\mathcal{M}_{H^{01},J^{01}}(x,y)$ the finite energy solutions of
(\ref{eq:comp_fl}) that satisfy (\ref{eq:assymptotes}). These spaces
are smooth manifolds of dimension $\mu(x,y)$ (the relative Maslov
index in this case being defined by a straightforward adaptation of
the definition  in \S\ref{subsubsec:Maslov_ind}). Gromov
compactifications also exist in this context and we shall denote them
by $\overline{\mathcal{M}}_{H^{01},J^{01}}(x,y)$. They are manifolds
with boundary and they verify:
\begin{equation}\label{eq:boundary_comp}
\partial\overline{\mathcal{M}}_{H^{01}}(x,y)=\bigcup_{z\in
I}\overline{\mathcal{M}}_{H}(x,z)\times
\overline{\mathcal{M}}_{H^{01}}(z,y)\cup \bigcup_{z'\in
I'}\overline{\mathcal{M}}_{H^{01}}(x,z')\times\overline{\mathcal{M}}_{H'}(z',y)
\end{equation}
Moreover, in the same way as the one described in the Appendix \ref{sec:appendix} it is
possible to show that these manifolds are manifolds with
corners.

Another useful remark concerns the functional
$\mathcal{A}_{H^{01}}(s,x):\R\times \mathcal{P}_{\eta}(L,L')\to \R$
which is defined by the action functional formula (\ref{eq:action})
but by using $H^{01}$ instead of $H$. This is clearly a homotopy
between $\mathcal{A}_{H}$ and $\mathcal{A}_{H'}$. Assume now that $H^{01}$ is
a monotone homotopy in the sense that $\frac{\partial
H^{01}}{\partial s}(s,t,y)\leq 0, \ \forall \  s,t,y\in \R\times
[0,1]\times M$.  In this case, if we put
$a_{H^{01}}(s)= \mathcal{A}_{H^{01}}(s,u(s,-))$
for $u$ a solution of (\ref{eq:comp_fl}), then
\begin{equation}\label{eq:monotone}\frac{d
a_{H^{01}}}{d s}=d \mathcal{A}^{s}_{H^{01}}(\frac{\partial
u}{\partial s})+\int_{0}^{1}\frac{\partial H^{01}}{\partial
s}(s,t,u(s,t))dt
\end{equation} and, by (\ref{eq:derivative}), the first term of
the sum is negative and the second is negative or null due to
monotonicity. In other words, monotone homotopies which have been
introduced in the symplectic setting by Floer and Hofer in
\cite{FlHo}, insure that the relevant action functionals decrease
along solutions of (\ref{eq:comp_fl}). Since both $H$ and $H'$ are constant
outside of a compact set we see that after possibly adding some
positive constant to $H$ we may assume that  $H(t,x)>H'(t,x)$ for all $t,x\in [0,1]\times M$.
As adding a constant to $H$ does not modify its Hamiltonian flow and only
changes $\mathcal{A}_{H}$ by the addition of the same constant, we may assume that
monotone homotopies as above always exist and we fix one such homotopy $H^{01}$
for the rest of this proof. To each element $u\in\mathcal{M}_{H^{01}}(x,y)$ we associate
a path $\gamma_{u}:[0,\mathcal{A}_{H}(x)-\mathcal{A}_{H'}(y)]\to \mathcal{P}_{\eta}(L,L')$
defined by the same formula as that used for (\ref{eq:def_path}) but with $\mathcal{A}_{H^{01}}$
instead of $\mathcal{A}_{H}$. We continue the construction in prefect analogy to that described in
\S \ref{subsubsec:defor_pseudo} and we thus get continuous maps
$$\overline{\gamma}_{x,y}:\overline{\mathcal{M}}_{H^{01}}(x,y)\to C_{x,y}\mathcal{P}$$
which are coherent with the maps constructed in
\S\ref{subsubsec:defor_pseudo} in the sense that an obvious analogue
of (\ref{eq:comm_eval_concat}) is verified as implied by
(\ref{eq:boundary_comp}). To pursue the construction along the lines
in \S \ref{subsubsec:add_path_sp} we first need to impose an
additional restriction on the path $w$ used to construct $\tilde{L}$:
we shall assume that $\{y(0): y\in I(L,L';\eta,H')\}\subset w$. With
this non-restrictive assumption and for any $x,y\in I(L,L';H)\cup
I(L,L';H')$, we define $\Phi_{x,y}=\mathcal{Q}_{x,y}\circ
\overline{\gamma}_{x,y}$ as in (\ref{eq:eval}). We pursue the
construction with the step described in \S
\ref{subsubsec:alg_constr}. This construction involves the choice of
$z_{0}\in I(L,L';H)$. We shall also need a similar choice $z'_{0}\in
I(L,L';H')$. To insure the compatibility of these choices we take
$z_{0}$ and $z'_{0}$ so that $\mu(z_{0},z'_{0})=0$ (it is easy to see
that such a couple necessarily exists). With these choices, the
construction described in \S \ref{subsubsec:alg_constr} applied to
$H$ and to $H'$ produces, respectively, matrices $A=(a_{xy})$ and
$A'=(a'_{xy})$ and chain complexes $\mathcal{C}(L,L';H)$,
$\mathcal{C}(L,L';H')$. There is an obvious analogue
$\{\tilde{s}_{xy}\}$ of the representing system of chains for the
moduli spaces $\mathcal{M}_{H^{01}}(x,y)$ so that this system is
compatible with both $\zeta=\{s_{xy}\}$ and with
$\zeta'=\{s'_{x'y'}\}$. The condition ii. in Definition
\ref{def:repres} is replaced by
$\partial\tilde{s}_{xy'}=\sum_{z}s_{xz}\times \tilde{s}_{zy'} +
\sum_{z'}\tilde{s}_{xz'}\times s_{z'y'}$ which reflects equation
(\ref{eq:boundary_comp}). The existence of such representing chain
systems for $H^{01}$ compatible with $\zeta$ and $\zeta'$ then
follows as in Lemma \ref{lem:repres_chain_syst}. Pursuing the
construction we obtain a matrix $B=(b_{xy})\in
M_{m,m'}(\mathcal{R}_{\ast})$ where, as in \S
\ref{subsubsec:alg_constr}, $m$ is the number of elements of
$I(L,L';H)$ and $m'$ is the number of elements in $I(L,L';H')$ and
$b_{xy}=\Phi_{x,y}(\tilde{s}_{xy})$. As in Proposition
\ref{prop:matrix_sq} we see that
\begin{equation}\label{eq:comp_matrix}
\partial B= A\cdot B + B\cdot A'~.~
\end{equation}
For $x\in I(L,L';H)$ we now define $$\mathcal{V}(x)=\sum_{y\in I(L,L';H')} b_{xy}\otimes y$$
and extend this to an $\mathcal{R}_{\ast}$-morphism. We then have
\begin{eqnarray*}
\mathcal{V}(d x)&=&
\mathcal{V}(\sum_{y'}a_{xy'}\otimes y')=\\=\sum_{y'}a_{xy'}\otimes \mathcal{V}(y')&=&
\sum_{y',z}a_{xy'}\cdot b_{y'z}\otimes z =\\= \sum_{z}(\sum_{y'}a_{xy'}\cdot b_{y'z})\otimes z
&=& \sum_{z}(\partial b_{xz}+\sum_{v}b_{xv}\cdot a'_{vz})\otimes z =\\=
\sum_{z}\partial b_{xz}\otimes z  + \sum_{v}(\sum_{z}b_{xv}\cdot a'_{vz}\otimes z)
&=& \sum_{z}\partial b_{xz}\otimes z + \sum_{v} b_{xv}\cdot d v =\\
= d (\sum_{v} b_{xv}\otimes v)&=& d\mathcal{V}(x)~.~
\end{eqnarray*}
Therefore, the map $\mathcal{V}$ so defined is a morphism of chain
complexes which we shall sometimes also denote by
$\mathcal{V}_{H^{01}}$ to emphasize the monotone homotopy to which it
is associated. If the choices of $z_{0}$ and $z'_{0}$ are compatible,
as above, then this morphism is of degree $0$. If $z_{0}$ and
$z'_{0}$ are independent, then this morphism could have a non-zero
degree. Assuming for now the compatible choices from above it is obvious that this morphism preserves filtrations and so it induces a morphism of spectral sequences. Moreover, by the definition of the
Floer comparison morphism $V_{H^{01}}:CF_{\ast}(L,L';J,H)\to
CF_{\ast}(L,L';J',H')$ (induced by the same monotone homotopy) we see
that the morphism induced by $\mathcal{V}_{H^{01}}$ at the $E^{1}$
term of our spectral sequences is the $H_{\ast}(\Omega L)$-module
morphism induced by $V$. But $H_{\ast}(V_{H^{01}})$ is an isomorphism
so $E^{2}(\mathcal{V})$ is also an isomorphism and so
$E^{r}(\mathcal{V})$ is an isomorphism for all $r\geq 2$. Obviously,
in case the choices for $z_{0}$ and $z'_{0}$ are not compatible, then
this is still an isomorphism up to translation and this proves point
a. of the theorem.

For the point b. notice that, for $x,y\in I(L,L';H)$ and $u\in \mathcal{M}_{H}(x,y)$
we have
\begin{equation}\label{eq:energy_action}
\mathcal{A}_{H}(x)-\mathcal{A}_{H}(y)=E_{L,L',H}(u)~.~
\end{equation} Therefore,
$\epsilon(L,L';H, J)=\min\{\mathcal{A}_{H}(x)-\mathcal{A}_{H}(y) : \mathcal{M}_{H,J}(x,y)\not=\emptyset\}$.
It has been proven by the second author together with Andrew Ranicki in \S 2.1 of \cite{CoRa}
that under the assumptions of the theorem and for the case of periodic orbits,
the Floer comparison morphism admits a retract. More precisely,
there exist  monotone homotopies $H^{01}$ and $G^{01}$ so that
$V_{G^{01}}\circ V_{H^{01}}$ is an isomorphism whose matrix is upper triangular
with $1$'s on the diagonal.
The exact same argument applies also here: the only difference with respect to the
proof of Theorem 2.1 in \cite{CoRa} is that we deal with orbits
starting in $L$ and ending in $L'$ instead of periodic orbits, everything else remains
the same. The fact that the matrix for $V_{G^{01}}\circ V_{H^{01}}$ is as above implies
that the matrix for $\mathcal{V}_{G^{01}}\circ \mathcal{V}_{H^{01}}$ is also upper
triangular with $1$'s on the diagonal. Therefore, $\mathcal{V}_{G^{01}}\circ \mathcal{V}_{H^{01}}$
is an isomorphism and this proves the claim.
\end{proof}

\subsubsection{Proof of Theorem \ref{theo:main} a.}\label{subsubsec:proof_a}
Point a. of Theorem \ref{theo:main} is a simple consequence of Theorem \ref{theo:comparison}
and of the naturality property recalled in  \S \ref{subsubsec:naturality}.

In fact, we can as easily prove slightly more. For this we let
$\epsilon(L,L')=\epsilon(L,L';0, J)$ and we recall the setting: $L$,
$L'$ are as before and we have also the Lagrangian $L''$ which is
transversal to $L$ and the almost complex structure $J'$ so that the
complexes $\mathcal{C}^{J}(L,L')=\mathcal{C}^{\eta,J,\zeta}(L,L';0)$,
$\mathcal{C}^{J'}(L,L'')=\mathcal{C}^{\eta',J',\zeta'}(L,L'';0)$ are
defined as well as the associated spectral sequences $EF(L,L')$ and
$EF(L,L'')$. Assume also that we have a hamiltonian diffeomorphism
$\phi$ such that
$$\phi(L'')=L'\ , \ \eta(t)=\phi(\eta'(t)) \ ,
\forall t\in [0,1]~.~$$ We shall assume here that $\phi$ has a
compact support. This is not restrictive for our purposes because
$L$, $L'$ are compact. Denote by $\mathcal{T}$ the set of
$1$-periodic hamiltonians on $M$ which are constant outside some compact set and recall the Hofer norm (or energy) \cite{Ho2} of a compactly
supported hamiltonian diffeomorphism:
\begin{equation}\label{eq:hofer_norm}
||\phi||_{H}=\inf_{H\in\mathcal{T}\ ,\ \phi^{H}_{1}=\phi}\
(\sup_{x,t} H(t,x)-\inf_{x,t}H(t,x))
\end{equation}

\begin{cor}\label{cor:comp_lagran} Under the assumptions above:
\begin{itemize}
\item[a.] There exists a morphism of chain complexes, possibly of non-zero degree
$$\mathcal{W}:\mathcal{C}^{J}(L,L')\to \mathcal{C}^{J'}(L,L'')$$
which induces an isomorphism up to translation between the spectral
sequences $(EF^{r}(L,L'), d^{r})$ and $(EF^{r}(L,L''),d^{r})$ for
$r\geq 2$. \item[b.] If $||\phi||_{H}<\epsilon(L,L')/4$, then
$\mathcal{W}$ admits a retract.
\end{itemize}
\end{cor}

\begin{rem}\label{rem:complexity}{\rm
Point a. of Theorem \ref{theo:main}  is clearly the same as point a.
of Corollary \ref{cor:comp_lagran}.  In
view of the moduli-spaces interpretation of the differentials in $\mathcal{C}^{J}(L,L')$ we may interpret point b. of the corollary as saying that a small enough hamiltonian
isotopy of $L'$ can only increase the algebraic complexity of the
moduli spaces of pseudo-holomorphic strips. A different useful
formulation is that, if $\mathcal{C}^{J}(L,L')$
is not a retract of $\mathcal{C}^{J'}(L,L'')$ (for example if the number of intersection points in $L\cap L''$ is smaller than the
number of intersection points in $L\cap L'$), then at
least as much energy as $\epsilon(L,L')/4$ is needed to deform $L'$
into $L''$.}
\end{rem}

\begin{proof} Let $H\in\mathcal{T}$ be such that $\phi^{H}_{1}=\phi$.
Let $J_{\ast}$ be the almost complex structure on $M$ which satisfies
$\phi^{\ast}(J_{\ast})=J'$. Recall from \S \ref{subsubsec:naturality}
the map $b_{H}:\mathcal{M}_{L,L'',J',0}(x,y)\to
\mathcal{M}_{L,L',J_{\ast},H}(x,y)$ which is defined by
$(b_{H}(u))(s,t)=\phi_{t}^{H}(u(s,t))$ and is a homeomorphism
respecting the various compactifications. Obviously, this map is also
compatible with the maps $\Phi_{x,y}$ and so $b_{H}$ induces an
identification of the two chain complexes (in the sense that it gives
a base preserving isomorphism of chain complexes):
\begin{equation}\label{eq:natural_id}
\overline{b}_{H}:\mathcal{C}^{\eta',J',\zeta'}(L,L'';0)\to
\mathcal{C}^{\eta,J_{\ast},\zeta''}(L,L';H) \end{equation}
where $\zeta''$ is the image of $\zeta'$ by $b_{H}$.
 Clearly,
$\overline{b}_{H}$ induces an isomorphism up to translation between
the respective spectral sequences and as, by Theorem
\ref{theo:comparison} a., we also have a morphism
$$\mathcal{V}:\mathcal{C}^{\eta,J,\zeta}(L,L';0)
\to \mathcal{C}^{\eta,J_{\ast},\zeta''}(L,L';H)$$ which induces an
isomorphism at the level of the spectral sequences we conclude that
the composition $\mathcal{W}=\mathcal{V}\circ
(\overline{b}_{H})^{-1}$ verifies point a.

Point b. of Theorem \ref{theo:comparison} shows that if
$\sup_{x,t}|H(t,x(t))|\leq \epsilon(L,L')/4$
for all $x\in \mathcal{P}_{\eta}(L,L')$, $t\in [0,1]$, then the
conclusion at point b. of the corollary holds. We pick a
hamiltonian $H\in\mathcal{T}$ such that $\phi_{1}^{H}(L'')=L'$ and
$\sup_{x,t}H(t,x)-\inf_{x,t}H(t,x)=||\phi||_{H}+\delta$
where $\delta$ verifies $||\phi||_{H}+\delta\leq \epsilon(L,L')/4$.
By adding an appropriate constant to $H$ we may assume $\inf_{x,t}H(t,x)=0$
and this proves the second point of the corollary.
\end{proof}
\subsection{Proof of the main theorem. II: Relation to the Serre spectral sequence.}\label{sec:serre}
The purpose of this subsection is to show point d. of Theorem
\ref{theo:main}.

\subsubsection{Elements of classical Morse theory}\label{subsubsec:morse}
We shall fix here a Morse function $f:L\to \R$ and
we also fix  a Riemannian metric $\alpha$ on $L$ such that the pair
$(f,\alpha)$ is Morse-Smale.
The Morse-Smale condition means that, if we denote by $\gamma$ the
flow induced by the negative $\alpha$-gradient of $f$, $-\nabla f$,
then the unstable manifolds
$$W^{u}(P)=\{x\in L : \lim_{t\to-\infty}\gamma_{t}(x)=P\}$$ and the stable
manifolds $$W^{s}(Q)=\{x\in L : \lim_{t\to+\infty}\gamma_{t}(x)=Q\}$$
intersect transversely for any two critical points $P,Q\in Crit(f)$.
If the index of the critical points $P$ is equal to $p$, then
$W^{u}(P)$ is diffeomorphic to an open $p$-disk and $W^{s}(P)$ is
diffeomorphic to an open $(n-p)$-disk. It is easy to see that if
$\alpha\in \R$ is a regular value of $f$ such that $f(P)>\alpha>
f(Q)$, then the space of $\gamma$-flow lines that join $P$ to $Q$ is
parametrized by the intersection $W^{u}(P)\cap f^{-1}(\alpha)\cap
W^{s}(Q)$ which, due to the transversality assumption, is seen to be
a manifold of dimension $ind(P)-ind(Q)-1$. This moduli space of
negative gradient flow lines will be denoted by $M_{f,\alpha}(P,Q)$
and the space of all the points situated on elements of
$M_{f,\alpha}(P,Q)$ will be denoted by $M'_{f,\alpha}(P,Q)$ (to
shorten notation we shall sometimes omit the symbol $\alpha$). These
moduli spaces $M_{f}(-,-)$ have properties that parallel those of the
moduli spaces $\mathcal{M}_{L,L',H}(-,-)$ as described in \S
\ref{subsubsec:Maslov_ind} and in \S \ref{subsubsec:gromov_comp} but
with the set $I(L,L';H)$ replaced by the set of critical points of
$f$, $Crit(f)$, and with the difference of Morse indexes
$ind(P)-ind(Q)$ used instead of the Maslov index $\mu(x,y)$. These
properties are much easier to prove for negative-gradient flow lines
than for pseudo-holomorphic strips and, in fact, historically the
Morse case has preceded and inspired Floer's machinery. From an
analytic point of view, the study of the moduli spaces
$M_{f,\alpha}(-,-)$ is clearly a simpler version of the study of
$\mathcal{M}_{L,L',H}(-,-)$ because negative gradient flow lines are
solutions $v:\R\to L$ of the equation
$$\frac{dv}{ds}+\nabla f (v)=0$$
which may be treated as a simplified version of equation
(\ref{eq:floer}). This approach has been developed in detail in
\cite{Sch}.

\subsubsection{Morse flow lines and pseudo-holomorphic strips.}
\label{subsubsec:Morse_fl_strips} There exists another deeper
relation between the moduli spaces of Morse trajectories and the
moduli spaces of pseudo-holomorphic strips which has been established
by Floer \cite{Fl3} and which we now recall. Recall that there exists a
neighbourhood of $L$ in $M$ which is symplectically equivalent to the
total space of a disk bundle associated to the cotangent bundle $T^{\ast}L$.
We shall denote this neighbourhood by $DT^{\ast}L$ and we consider the Hamiltonian
$\overline{f}:DT^{\ast}L\to \R$, $\overline{f}=-f\circ \pi$ where
$\pi:DT^{\ast}L\to L$ is the projection. Notice that if
$L_{f}=\phi_{1}^{\overline{f}}(L)$, then $L_{f}$ is precisely the image of
$-df$ and $L\cap L_{f}$ coincides with the set of critical points of
$f$ (we assume here that $f$ is small enough so that the image of $df$
is contained in $DT^{\ast}L$). The fact that $f$ is a Morse function
is equivalent to the
transversality of $L_{f}$ and $L$. For any $x,y\in L\cap L_{f}$, it
is natural to define a map $c_{f}:M_{f,\alpha}'(x,y)\to
C^{\infty}(\R\times [0,1], M)$ by
$(c_{f}(v))(s,t)=\phi_{t}^{\overline{f}}(v(s))$. Floer's result is
that, if $f$ is sufficiently small in $C^{2}$-norm, then there exists
a (time-dependent) almost complex structure $J^{f}$ such that the
image of this map belongs to $\mathcal{M}'_{L,L_{f},J^{f},0}(x,y)$
and, moreover, the resulting application $c_{f}:M'_{f,\alpha}(x,y)\to
\mathcal{M}'_{L,L_{f},J^{f},0}(x,y)$
 is a diffeomorphism. The fact that $c_{f}$ is surjective is in itself
 highly non-trivial as, apriori, $\mathcal{M}'_{L,L_{f},J^{f},0}(x,y)$ could contain
 some ``long" Floer trajectories which do not belong to $DT^{\ast}L$, however,
 Gromov compactness together with our assumptions on the lack of bubbling
 imply that by making $f$ sufficiently small (for example
 by replacing it with $\lambda f$ with $\lambda>0$ and small) this does not happen.
 Obviously, this application induces a diffeomorphism
$$l_{f}:M_{f,\alpha}(x,y)\to \mathcal{M}_{L,L_{f},J^{f},0}(x,y)$$ and it is
clear that this is compatible with the compactifications and the
stratifications on the two sides.

\subsubsection{The Morse spectral
sequence.}\label{subsubsec:Morse_ss} We now let $w$ be a path in $L$
which is embedded and joins all critical points of $f$. We then
define the quotient map $q:L\to \tilde{L}$ as in \S\ref{subsubsec:add_path_sp}. Following the scheme in \S\ref{subsec:constr_SS} it is easy to see how to build a
spectral sequence asociated to the Morse-index filtration of the
$\mathcal{R}_{\ast}$-chain complex $C^{f,\alpha}= (C_{\ast}, d)$
which is defined by $C_{k}=\bigoplus_{q+p=k}\mathcal{R}_{q}\otimes
\Z/2<Crit_{k}(f)>$ (where $Crit_{k}(f)$ are the critical points of
$f$ which are of Morse index equal to $k$) and
$$dx=\sum_{y\in Crit(f)}m_{xy}\otimes y~.~$$
As in formula (\ref{eq:coefficents}), the coefficients $\{m_{xy}\}$
are defined as images of a representing chain system for the moduli
spaces $M_{f,\alpha}(x,y)$ by the map $v\in M_{f}(x,y)\to
q\circ s_{v}\in\Omega'\tilde{L}$ where $$s_{v}:[0,f(x)-f(y)]\to L$$ is
a reparametrization of $v$ such that $s_{v}(t)=z \Leftrightarrow
f(z)=f(x)-t$. Further, as in \S\ref{subsubsec: ss_def}, the
filtration $F^{k}C= \mathcal{R}_{\ast}\otimes \Z/2<Crit_{j}(f): j\leq
k>$ induces a spectral sequence which we shall denote by
$E(f,\alpha)=(E^{r}_{pq}(f,\alpha), d^{r})$ (again, sometimes we
shall omit $\alpha$ to shorten notation). A result similar to
Proposition \ref{prop:simple_theo} is true after replacing the Floer
complex with the Morse complex and Floer homology with the usual
homology of $L$.

\subsubsection{Reduction to the Morse
case.}\label{subsubsec:red_Morse} We now assume that $f$ is
sufficiently $C^{2}$-small so that Floer's result mentioned above
applies. Clearly, we may extend both $\overline{f}$ and $J^{f}$ to a
hamiltonian and, respectively, an almost complex structure defined on
all of $M$ which shall be denoted by the same respective symbols.

If we let $\eta_{0}$ coincide with $z_{0}$ and both be equal to a
minimum of $f$, then we see that the map $l_{f}$ of
\S\ref{subsubsec:Morse_fl_strips} induces an identification of chain
complexes $\overline{l}_{f}:C^{f,\alpha}\to
\mathcal{C}^{\eta_{0},J^{f}}(L,L_{f};0)$. This obviously preserves
filtrations and identifies the spectral sequences $E(f,\alpha)$ and
$EF(L,L_{f};\eta_{0},J^{f},0)$ .

We now turn to the setting of Theorem \ref{theo:main} d. Therefore,
$L'$ is hamiltonian isotopic to $L$. By Corollary
\ref{cor:comp_lagran}, we then have that
$(EF^{r}(L,L';\eta,J,0),d^{r})$ is isomorphic up to translation to
$(EF^{r}(L,L_{f};\eta_{0},J^{f},0),d^{r})$ for $r\geq 2$. At the same
time, as discussed above, this last spectral sequence is isomorphic
to $E(f,\alpha)$. Thus, to prove Theorem \ref{theo:main} d., it
suffices to show that $E^{r}(f,\alpha)$ is isomorphic to the Serre
spectral sequence of $\Omega L\to PL \to L$ for $r\geq 2$.

\subsubsection{The Morse and Serre spectral
sequences.}\label{subsubsec:Morse_Serre_ss}

The purpose of this sub-subsection is to conclude the proof of
Theorem  \ref{theo:main} by showing:

\begin{theo} \label{theo:Morse_comp}
Assume that $f:L\to \R$ is a Morse function and $\alpha$ is a
riemannian metric on $L$ so that the spectral sequence
$E(f,\alpha)=(E^{r}_{pq}(f,\alpha),d^{r})$ is defined as in
\S\ref{subsubsec:Morse_ss}. For $r\geq 2$ there exists an isomorphism
of spectral sequences between $E(f,\alpha)$ and the Serre spectral
sequence $E(L)=(E^{r}_{pq},d^{r})$ of the path loop fibration of base
$L$.
\end{theo}

\begin{proof} We may assume that the function $f$
has just one minimum that we shall denote by $B$. We also assume that
$f(B)=0$. It is not
restrictive to also suppose that $f$ is self-indexed which means that
for any critical point $x$ of $f$ we have that $f(x)=ind_{f}(x)$.
Take $\epsilon$ to be a very small positive constant and let
$L_{k}=f^{-1}(-\infty,k+\epsilon]$. Of course, by classical Morse
theory, $L_{k}$ is homotopy equivalent to a $k$-th dimensional
skeleton of $L$. Consider the  path-loop fibration $\Omega L\to PL
\to L$ and let $\Omega L\to E_{k}\to L_{k}$ be the pull-back of this
fibration over the inclusion $L_{k}\hookrightarrow L$. We consider
the filtration $\Omega L=E_{0}\hookrightarrow \ldots
E_{k}\hookrightarrow E_{k+1}\hookrightarrow \ldots PL$ and the
resulting filtration of the cubical chain complex $S_{\ast}(PL)$
which is given by  the $S_{\ast}(E_{k})$'s. The spectral sequence
associated to this filtration is, by definition, the Serre spectral
sequence of the statement \cite{Spanier}. The proof of the theorem
consists of the following two steps:
\begin{itemize}
\item[i.] there exists a morphism of chain complexes $\xi:
C^{f,\alpha}\to S_{\ast}(PL)$ so that $\xi(F^{k}C)\subset
S_{\ast}(E_{k})$. Such a $\xi$ induces a morphism of spectral
sequences denoted by $E(\xi):E(f,\alpha)\to E(L)$. \item[ii.] with
$\xi$ as above the morphism $E^{2}(\xi)$ is an isomorphism.
\end{itemize}

Before proceeding with the proof we need to make a few adjustments.
First, notice that instead of using unit paths in the definition of the path-loop
fibration we may as well use Moore paths - these are paths parametrized by arbitrary
intervals $[0,a]$. The resulting fibration is denoted by $\Omega'L\to P'L\to L$
and the associated filtration is $\{E'_{k}\}$. Moreover, as $q:L \to \tilde{L}$
is a homotopy equivalence we may replace the spaces $L_{k}$, $E'_{k}$ by their respective image $\tilde{L}_{k}$ and $\tilde{E}_{k}\subset P'\tilde{L}$ in the latter case via the induced  map $P'q:P'L\to P'\tilde{L}$ (the two induced spectral sequences being obviously isomorphic). For further use, notice also that there is an obvious action $\cdot: \Omega'\tilde{L}\times P'\tilde{L}\to P'\tilde{L}$
which induces $\mathcal{R}_{k}\otimes S_{q}(P'\tilde{L})\to S_{k+q}(P'\tilde{L})$.

\subsubsection{Blow-up of unstable manifolds.}
\label{subsubec:blow_unst}

The first step is based on a geometric construction which, as we shall see,
is of independent interest. This construction provides an efficient
geometric description for the compactification of the unstable manifolds
of $f$.

We fix $x\in Crit(f)$.
Notice that for each element  $v\in \overline{M}_{f}(x,B)$
there exists some $k\geq 0$ such that $v=(v_{1},v_{2},\ldots , v_{k})$ with
$v_{1}\in M_{f}(x,x_{1}) \ , \ \ldots \ , v_{i}\in M_{f}(x_{i-1},x_{i}) \ , \ldots \ , \
v_{k}\in M_{f}(x_{k-1},B)$. This writing is of course unique. We recall the
parametrizations $s_{v}$ for the flow lines represented by $v\in M_{f}(x,B)$
as defined in \S\ref{subsubsec:Morse_ss}. Clearly,
this parametrization extends in an obvious way to the elements
$v=(v_{1},v_{2},\ldots, v_{k})\in \partial\overline{M}_{f}(x,B)$ and we shall continue to denote the parametrization of these elements by $s_{v}$.

 We consider the space
$\widehat{M}(x)$ which is defined as the topological quotient of the space
$\overline{M}_{f}(x,B)\times [0,f(x)]$ by the equivalence relation induced by:
$$((v_{1},\ldots , v_{k}),t)\sim ((v'_{1},\ldots , v'_{k}),t) \ \ {\rm if} \ \ v_{i}=v'_{i} \ \ \forall i \ {\rm with}\ \ f(x_{i-1})
> t \ , \ v_{i}\in M_{f}(x_{i-1},x_{i})~.~$$
In short, two couples $(v,t),(v',t)\in \overline{M}_{f}(x,B)$ are identified
in $\widehat{M}_{f}(x,B)$ if the (possibly broken) negative gradient trajectories
of $v$ and $v'$ coincide above level $t$.
Notice that if $(l_{n},t_{n})\sim
(l'_{n},t_{n})$, where $l_{n},l'_{n}\in \overline{M}_{f}(x,B)$ with
$l_{n}\to l\in \overline{M}_{f}(x,B)$, $l'_{n}\to l'\in \overline{M}_{f}(x,B)$,
 $t_{n}\to t$, then $(l,t)\sim (l',t)$ and $\widehat{M}(x)$ is Hausdorff.

It is useful to introduce the map $S:\overline{M}_{f}(x,B)\times [0,f(x)]\to L$ defined by $S(v,\tau)=s_{v}(f(x)-\tau)$ ($S(v,\tau)$ is thus simply the intersection
of the trajectory $v$ with $f^{-1}(\tau)$). This map factors as
$$S:\overline{M}_{f}(x,B)\times [0,f(x)]\stackrel{k}{\longrightarrow}
\widehat{M}(x)
\stackrel{o}{\longrightarrow}L$$ where $k$ is the quotient map.

We call the space $\widehat{M}(x)$ the {\em blow-up} of the unstable manifold
$W^{u}(x)$. As we shall see this is justified by a number of remarkable
properties of this space. We start with the most immediate. First, the image of $o$ is included and is onto the closure of $W^{u}(x)$.  Secondly,
all the points in $\overline{M}_{f}(x,B)\times \{f(x)\}$ belong to a unique equivalence class which we shall denote by $\ast$ . Furthermore, define paths $s'_{v}:[0,f(x)]\to \widehat{M}(x)$
by the formula $s'_{v}(\tau)=k(v,f(x)-\tau)$. Obviously, $s_{v}'(0)=\ast$ and $s_{v}=o\circ s'_{v}$. Moreover, for each $y\in\widehat{M}(x)$
there exists a unique $t\in [0,f(x)]$
and a unique path $\overline{y}:[0,t]\to \widehat{M}(x)$ such that
$\overline{y}(t)=y$ and $\overline{y}(\tau)=s'_{v}(\tau), \ \forall \tau\in [0,t]$
for some $v\in\overline{M}_{f}(x,B)$. It is easy to see that the map
$$\beta:\widehat{M}(x)\to P'(\widehat{M}(x)),\ \beta(y)=\overline{y}$$
is continuous.
As we also have that $\overline{y}(0)=\ast$ this shows that $\widehat{M}(x)$ is contractible by a contraction that pushes each $y\in \widehat{M}(x)$ along the
path $\overline{y}$ till it reaches $\ast$. We formulate a stronger property next.
For this first notice that for all $y\in Crit(f)\cap \overline{W^{u}(x)}$
there is a natural inclusion $M_{f}(x,y)\times \widehat{M}(y)\subset
\widehat{M}(x)$ which is induced by the product of inclusions
$(M_{f}(x,y)\times \overline{M}_{f}(y,B))\times [0,f(y)]\hookrightarrow \overline{M}_{f}(x,B)\times [0,f(x)]$.

\begin{lem}\label{lem:blow_struct}
The space $\widehat{M}(x)$ is homeomorphic to a closed disk of
dimension equal to $ind_{f}(x)$. Moreover,
$$\partial\widehat{M}(x)=\bigcup_{y}M_{f}(x,y)\times \widehat{M}(y)~.~$$
\end{lem}

\begin{rem}\label{rem:CW_Morse} {\rm a. As we shall see below, the actual proof of Theorem \ref{theo:Morse_comp} only uses that $\widehat{M}(x)$ is a topological manifold with a boundary described as in the statement of the lemma and that $\ast$ has a neighbourhood homeomorphic to a disk. Of course, the fact that $\widehat{M}(x)$ is a topological manifold is not surprising: this space is obviously homeomorphic to the space of all (appropriately parametrized) possibly broken gradient flow lines that join $x$ to points in $L$.

b. While the definition of $\widehat{M}(x)$ based on the equivalence relation $\sim$ is new, the space of all geometric, possibly broken, flow lines ending
in points of $L$ and originating in $x\in Crit(f)$ has appeared before in the Morse theoretic literature, for example, in \cite{Hutchings1}\cite{Hutchings2}. In this last paper, in the proof of Lemma 3.1, it is also claimed that $\widehat{M}(x)$
is homeomorphic to a disk. This fact is of independent interest as it immediately implies that the union of the closures of the unstable manifolds of a self-indexed Morse-Smale function has a natural $CW$-complex structure - the attaching map corresponding to the cell associated to $x$ being simply $o|_{\partial\widehat{M}(x)}$. However, no explicit proof
of the existence of a homeomorphism between $\widehat{M}(x)$ and a closed disk appears to exist in the literature and so we provide one here.}
\end{rem}

{\it Proof of the Lemma.} We fix $i=ind_{f}(x)$ and we recall that
$f$ is self-indexed.

We start by verifying explicitly that $\widehat{M}(x)$ is a topological manifold whose boundary has the description of the statement. We first notice that the restriction of $k$ to  $M_{f}(x,B)\times (0,f(x))$ is a homeomorphism onto its image.
Moreover, the definition
of the equivalence relation $\sim$ directly implies that
the restriction
\begin{equation}\label{eq:local_hom}o|:k(\overline{M}_{f}(x,B)\times [(i-1)+\delta,f(x)])\to W^{u}(x)\cap f^{-1}[(i-1)+\delta,+\infty)
\end{equation}
is a homeomorphism for any small positive $\delta$  where $o$ is, as before, the factor of the map $S$
(for further use, notice also that $W^{u}(x)\cap f^{-1}[(i-1)+\delta, +\infty)$ is homeomorphic to an $i$-disk).

Consider a point $(v,t)\in \overline{M}_{f}(x,B)\times [0,f(x)]$ such that
$v=(v_{1},\ldots,v_{k})\in M_{f}(x,x_{1})\times\ldots \times M_{f}(x_{k-1},B)$
and $t> f(x_{1})$. We want to notice that $k(v,t)$ has a neighbourhood homeomorphic
to an $i$-disk. Indeed, for $\lambda$ sufficiently close to $f(x)$, the point $k(v,\lambda )$ does have such a neighbourhood $V$ because of the homeomorphism at (\ref{eq:local_hom}). This neighbourhood $V$ verifies $V\subset \bigcup_{y\in C(x_{1})}k(M_{f}(x,y)\times \overline{M}_{f}(y,B)\times [0,f(x)])$ where
$C(x_{1})=\{ y\in Crit(f) : x_{1}\in \overline{W^{u}(x)}\cap \overline{W^{s}(y)}\}$.
But this means that, if $V$
is sufficiently small, we may isotope it by sliding it along the paths $s'_{r}, r\in V$ till we get a neighbourhood of $k(v,t)$. ``Sliding" along
the paths $s'_{v}$ is given by $$h((k(r,t'),\tau)=s'_{r}(t'+\tau)$$ and is
well defined and an isotopy when restricted to $k(M_{f}(x,y)\times\overline{M}_{f}(y,B)\times [s,f(x)])$ as long as $\tau+s>f(y)$.
As $y\in C(x_{1})$ we have that $f(y)\leq f(x_{1})$ and thus sliding is indeed possible.

\subsubsection*{a. First look at boundary points.} Next, to continue the proof of the Lemma, we need to show that each point belonging to some $M_{f}(x,y)\times \widehat{M}(y)$ has a neighbourhood homeomorphic to a semi-disk. Let $z=k(v,t)$ with $v=(v_{1},\ldots, v_{k})\in M_{f}(x,x_{1})\times M_{f}(x_{1},x_{2}) \times \ldots M_{f}(x_{k-1},B)$ so that $f(x_{j-1})> t > f(x_{j})$. Because we are only interested in a neighbourhood of $z$ we may assume that the interval $(f(x_{j}),f(x_{j-1}))$ is regular and, in particular,
$t$ is a regular value of $f$.  Recall that
the point $o(z)$ is the intersection with $f^{-1}(t)$ of the broken negative gradient flow line of $f$ represented by $v$.

Let $\tilde{M}_{t}(x)=\{z\in C^{0}([0,f(x)-t],M) : \exists v\in
\overline{M}_{f}(x,B), \ z=s_{v}|_{[0,f(x)-t]} \}$. In short,
a path in $\tilde{M}_{t}(x)$ joins the point $x$ to some point in $f^{-1}(t)$
and it coincides geometrically to the part of a
negative-gradient (possibly broken) flow line of $f$ which is above (and on) level $t$.
Clearly, the spaces $\tilde{M}_{t}$ for $t$ such that $f(x_{j-1})> t \geq f(x_{j})$
are canonically identififed with $\tilde{M}_{j}=\tilde{M}_{f(x_{j})}$. Obviously, for  our fixed point $z=k(v,t)$ there exists a unique point $z'\in \tilde{M}_{j}(x)$ such that $o(z)=z'(f(x)-t)$ (the parametrization used for $z'$ is similar to that used for the paths $s_{v}$). In fact, in view of the definition of $\sim$  it is immediate to see that
the application $z\to (z',t)$ is a {\emph local}
homeomorphism defined on a neighbourhood of $z \in \widehat{M}(x)$ and with values in $\tilde{M}_{j}(x)\times (f(x_{j}),f(x_{j-1}))$. Now notice that
$\tilde{M}_{j}(x)$ is a compact topological manifold whose boundary consists
as usual of broken trajectories. This means that in case $x_{j-1}\not=x$
the trajectory $v$ is broken at $x_{j-1}$ and thus $z$ is mapped by this local homeomorphism to a point in $\partial\tilde{M}_{j}(x)\times (f(x_{j}),f(x_{j-1}))$.
Therefore, $z$ has a semi-disk
neighbourhood in $\widehat{M}(x)$.

\subsubsection*{b. Local study around breaking points.}
A slightly different argument is needed for the
points $k(v,t)$ with $v=(v_{1},\ldots, v_{k})$ as before but with
$t=f(x_{j-1})$. The first such case corresponds to $j=2$. The key observation
is that $x_{1}$ has inside $W^{u}(x_{1})$ a neighbourhood $U$ which is homeomorphic to a disk (of dimension $ind_{f}(x_{1})$). The element $v_{1}$ also has
a neighbourhood $V$ in $M_{f}(x,x_{1})$ which is homeomorphic to a disk.
Together with the definition of $\sim$ this shows that $k(v,t)$ has a
neighbourhood in $\widehat{M}(x)$ which is homeomorphic to the
product $U\times V\times [0,1)$.

 To see this we study the problem locally in a neighbourhood of $x_{1}$ in $L$. We may assume that $f$ is in normal form around $x_{1}$. Let $a=f(x_{1})$ and let $\epsilon, \delta$ be very small positive constants. Let $W$ be a neighbourhood of $x_{1}$ which consists of all the points $x\in f^{-1}[a-\epsilon,a+\epsilon]$ that are situated on flow lines of $f$ whose intersection with $f^{-1}(a)$ is at distance less than $\delta$ from $x_{1}$.
We remark that $D=W^{s}(x_{1})\cap W$ is homeomorphic to a disk of dimension
$n-ind_{f}(x_{1})$, similarly $D'=W^{u}(x_{1})\cap W$ is a disk of dimension
$ind_{f}(x_{1})$. We let $S^{s}=\partial D$ and $S^{u}=\partial D'$, $A^{s}=W\cap
f^{-1}(a+\epsilon)$, $A^{u}=W\cap f^{-1}(a-\epsilon)$. Notice that $A^{s}$, $A^{u}$ are respectively tubular neighbourhoods of
$S^{s}$ inside $f^{-1}(a+\epsilon)$ and of $S^{u}$ inside $f^{-1}(a-\epsilon)$. Therefore, $A^{s}=S^{s}\times D''$, $A^{u}= D'''\times S^{u}$ with $D''$ a disk of dimension $ind(x_{1})$ and $D'''$ a disk of dimension $n-ind(x_{1})$. Moreover,
 the flow provides a homeomorphism between $A'=A^{s}\backslash (S^{s}\times \{0\})$
and $A''=A^{u}\backslash (\{0\}\times S^{u})$. In view of this we may identify
both $A'$ and $A''$ with $S^{s}\times S^{u}\times (0,\delta)$. The set of all
paths in $W$ which join $A^{s}$ to $A^{u}$, which are parametrized by the values of $f$ (similarly to the $s_{v}$'s) and which coincide geometrically to portions of possibly broken flow lines of $f$ is identified with $S^{s}\times S^{u}\times [0,\delta)$
(the broken flow lines correspond to $S^{s}\times S^{u}\times \{0\}$). We now consider
the space $K(x_{1})=(S^{s}\times S^{u}\times [0,\delta))\times [a-\epsilon,a+\epsilon]/\sim'$ where $\sim'$ is the analogue of $\sim$ for our
paths in $W$. It is easy to see that the existence of our
semi-disk neighbourhood  of $k(v,t)$ inside $\widehat{M}(x)$ follows
if we show that any point of type $[(x,y,0),a]$,
has a similar semi-disk neighbourhood inside $K(x_{1})$.
We have $K(x_{1})\approx S^{s}\times (S^{u}\times [0,\delta)\times [a-\epsilon,a+\epsilon]/\sim'')$ where $\sim''$ is the equivalence relation
 induced by $(x,0,t)\sim''(x',0,t)$ if $t\geq a$. This means that we reduced the
 problem to studying the  space $K'(x_{1})= (S^{u}\times [0,\delta)\times [a-\epsilon,a+\epsilon])/\sim''$. Recall that $S^{u}=\partial D'$. It is easy to check now that $K'(x_{1})$ is homeomorphic to the cylinder $D'\times [a-\epsilon,a+\epsilon]$ from which has been eliminated the interior and the base
of a circular cone of hight
$[a-\epsilon,a)$, whose base lies in the interior of $D'\times \{a-\epsilon\}$ and
whose vertex corresponds to $(y,0,a)$. This shows our claim.

An immediate adaptation of this argument
also works when $t=f(x_{j-1})$ even for $j>2$ and this shows that $\widehat{M}(x)$
is indeed a compact topological manifold with boundary.

\subsubsection*{c. Homeomorphism to a disk.}To end the proof of the lemma we still need to show that $\widehat{M}(x)$ is homeomorphic to a disk. The idea is to construct a copy $\partial' \widehat{M}(x)$ of $\partial\widehat{M}(X)$ such that $\partial' \widehat{M}(x)$ is
contained in $\widehat{M}(x)$, it is transverse to the paths $s'_{v}$ and it bounds a topological manifold $\widehat{M}'(x)$ which contains $\ast$ and is
homeomorphic to $\widehat{M}(x)$.
Recall that a neighbourhood $U$ of $\ast$ as
in (\ref{eq:local_hom}) has a boundary that is also transverse to the paths $s'_{v}$. By sliding along these paths it follows that $\widehat{M}'(x)$
is homeomorphic to a disk and as this manifold is homeomorphic to $\widehat{M}(x)$ this concludes the proof.
Before starting this construction we make explicit the notion of transversality used here: given a separating hypersurface $V$ of a topological manifold $N$ and a path
$g:[-a,a]\to N$ such that $g(0)\in V$ we say that $g$ is transversal to $V$ if
for some neighbourhood $U$ of $V$ such that $U\backslash V=U_{0}\bigsqcup U_{1}$
there exists $\epsilon>0$ and $i\in\{0,1\}$ such that $\forall t\leq\epsilon$ we have
$g(-t)\in U_{i}$, $g(t)\in U_{1-i}$.

To construct $\widehat{M}'(x))$ we first fix the notation
$\partial(y)=\overline{M}_{f}(x,y)\times\overline{M}_{f}(y,B)$,
$D(y)=\partial(y)\times [0,f(x)]$
and we let $s''_{v}$ be the path $v\times [0,f(x)]$ in $\overline{M}_{f}(x,B)\times [0,f(x)]$.
We now intend to construct for each $y\in Crit(f)\cap \overline{W^{u}(x)}, y\not=B$
a map
 $$f_{y}: D(y)\to \overline{M}_{f}(x,B)\times [0,f(x)]$$
 which is a homeomorphism onto its image - we shall denote this image by $D'(y)$ -  and has the following additional properties: $f_{y}(v,t)=(v,t)$ if $t\geq f(y)$; $s''_{f_{y}(v,t)}$ is transverse to $D'(y)$ at the point $f_{y}(v,t)$ whenever $t< f(y)$;
$D'(y)$ together with $\partial(\overline{M}_{f}(x,B)\times [0,f(x)])\backslash D_{y}$ bound a topological manifold with boundary $M'_{y}\subset \overline{M}_{f}(x,B)\times [0,f(x)]$ which is homeomorphic to
$\overline{M}_{f}(x,B)\times [0,f(x)]$ and contains $\ast$; if $(v,t)\sim (v',t)$, then $f_{y}(v,t)=f_{y}(v',t)$.
The construction of this auxiliary application is as follows. As $\overline{M}_{f}(x,B)$ is a manifold with corners
and $\partial(y)$
is a part of the boundary of $\overline{M}_{f}(x,B)$
there exists a collar neighbourhood $U(y)$
of $\partial(y)$
inside $\overline{M}_{f}(x,B)$. In particular,
there exists a homeomorphism $f':\partial(y)\times [0,\epsilon)\to U(y)$ so that $f'((v,w),\tau)=v\#_{\tau}w$
where $v\#_{\tau}w$ is the flow line obtained by gluing $v$ to $w$ at $y$ with gluing parameter $\tau$. Of course, for this we need to choose a particular gluing formula
(we may  do this  as discussed in the Appendix
\ref{sec:appendix} in the obviously harder Floer case) and we choose the gluing
parameter in such a way that $v\#_{0}w$ coincides
with $(v,w)$. More generally, for $\tau$ small
enough and $v=(v_{1},\ldots,v_{i})\in M_{f}(x,x_{1})\times\ldots\times M_{f}(x_{i-1},y)$,
$w=(w_{1},\ldots,w_{j})\in M_{f}(y,y_{1})\times\ldots\times
M_{f}(y_{j-1},B)$ we let $v\#_{\tau}w$ be the
element $(v_{1},\ldots, v_{i}\#_{\tau}w_{1},
\ldots, w_{j})\in M_{f}(x,x_{1})\times\ldots\times M_{f}(x_{i},y_{1})\times \ldots\times M_{f}(y_{j},B)$.
As a consequence of the parametrization of the
corners of $\overline{M}_{f}(x,B)$ as described in the Appendix we obtain that $f'$ so defined is a homeomorphism.
We also notice that if $((v,w),t)\sim (v',w'),t)$
with $t<f(y)$, then $(v\#_{\tau}w,t)\sim (v'\#_{\tau}w',t)$ and so
we also have $(f'((v,w),\tau),t)\sim (f'((v',w'),\tau),t)$. We now
let $\epsilon'<\epsilon$ and consider a smooth one parameter family of functions
$q^{s}_{y}:[0,f(x)]\to [0,\epsilon']$ such that for each $s\in [0,1]$, $q^{s}_{y}$ is decreasing and smooth, $q^{s}_{y}|_{[f(y),f(x)]}=0$, $q^{s}_{y}$ is strictly decreasing on $[0,f(y)]$,
$q^{s}_{y}(0)\leq\epsilon'$ and, moreover, for any fixed $t$, $q^{-}_{y}(t)$ is a function increasing in $s$ and $q^{0}\equiv 0$. We now define $f^{s}_{y}(v,t)=(f'(v,q^{s}_{y}(t)),t)$ and we let $f_{y}=f^{1}_{y}$.

We pursue with the construction of  $\widehat{M}'(x)$.
We let $D'_{s}(y)=Im(f^{s}_{y})$ (so that $D'_{1}(y)=D'(y)$, $D'_{0}(y)=D(y)$), $V^{s}(y)=\bigcup_{0\leq s'\leq s}D'_{s'}(y)$ (so that the slice of $V^{s}(y)$ of height
$t$ is a tubular neighbourhood of $\partial(y)$ in $\overline{M}_{f}(x,B)\times\{t\}$), $V^{s}=\bigcup_{y\not=B}V^{s}(y)$ and
$M_{s}'=\overline{M}_{f}(x,B)\times [0,f(x)]\backslash V^{s}$ (so that $M_{0}'=\overline{M}_{f}(x,B)\times [0,f(x)]$).
Notice that, for each $s>0$, and for each  $v\in M_{f}(x,B)$ the
path $s''_{v}$ is transversal to $\partial M_{s}'$.
We now define $\widehat{M}'_{s}=k(M_{s}')$, $E_{s}(y)=k(D'_{s}(y)\cap \overline{M}(x,B)\times [0,f(y)])$, $W_{s}(y)=k(V^{s}(y)\cap \overline{M}(x,B)\times [0,f(y)])$. By definition recall that $V^{s}(y)\cap \overline{M}(x,B)\times
[f(y),f(x)]=\partial(y)\times [f(y),f(x)]=D'_{s}(y)\cap\overline{M}(x,B)\times [f(y),f(x)]$. Because of this, as $f_{y}^{s}$ respects the relation $\sim$ and
 as the identifications producing $\overline{M}_{f}(x,y)\times \widehat{M}(y)\subset \partial\widehat{M}(x)$ occur only on the boundary of $\partial(y)\times [0,f(y)]$, it follows that $E_{s}(y)$ is a copy of $\overline{M}_{f}(x,y)\times\widehat{M}(y)$
which verifies $E_{s}(y)\cap \overline{M}_{f}(x,y)\times\widehat{M}(y)=\overline{M}_{f}(x,y)\times\ast$.
Clearly, $W_{l}(y)=\bigcup_{0\leq s\leq l}E_{s}(y)$ and $\widehat{M}'_{s}=
\widehat{M}(x)\backslash (\bigcup_{y\not=B} W_{s}(y))$. The description of
$W_{s}(y)$ shows that, by starting with the $y$'s
of lowest index, and proceeding by induction we may isotope $\widehat{M}'_{s}$
to $\widehat{M}(x)$ and thus these two spaces are homeomorphic.
To conclude the proof we only need to show that the paths $s'_{v}$ are
transverse to the boundary of $\widehat{M}'(s)$ and we may then take $\widehat{M}'(x)=\widehat{M}'(1)$. This transversality
is already clear for the points on the ``bottom" - the points
that belong to $M_{f}(x,B)\times \widehat{M}(B)$.
The transversality of the paths $s''_{v}$ to $\partial M_{s}'$ for
$v\in M_{f}(x,B)$  implies that for each such $v$, the path
$s'_{v}$ is transversal to $\partial\widehat{M}'_{s}$. As $\bigcup_{y}\partial(y)=\partial\overline{M}_{f}(x,B)$,
the only case that remains to be
discussed is that of transversality at the
points $k(v,t)\in \partial\widehat{M}'_{s}$ with $v$ in some $\partial(y)$.
By the description of $E_{s}(y)$, such a point $k(v,t)$ belongs to
$\overline{M}_{f}(x,y)\times \ast_{y}$ (where $\ast_{y}$ is the distinguished
point in $\widehat{M}(y)$), in particular $t=f(y)$.
We want to notice that, moreover, such a $k(v,t)$
actually belongs to $M_{f}(x,y)\times \ast_{y}$. Indeed,
a point $k(v,t)\in \partial\overline{M}_{f}(x,y)\times \ast_{y}$ has the property that there exists $x_{1}$, $f(x_{1})> f(y)$ so that $k(v,t)\in M_{f}(x,x_{1})\times \widehat{M}(x_{1})$.
This means that $t<f(x_{1})$ which implies $k(v,t)\in W_{s}(x_{1})$ and thus $k(v,t)\not\in \widehat{M}'_{s}$.
Therefore, we now consider $k(v,t)\in M_{f}(x,y)\times\ast_{y}$.
But from the transversality of $s''_{v}$ to $M_{f}(x,y)\times \{w\}\times \{f(y)\}$ for any $v=(v',w)\in M_{f}(x,y)\times \overline{M}_{f}(y,B)$ we immediately deduce
the transversality of $s'_{v}$ in this case and this concludes the proof of the lemma.

\subsubsection{Construction of $\xi$.}\label{subsubsection:constr_xi}
We first fix a representing chain system (recall
Definition \ref{def:repres})
$s_{xy}\in S_{ind(x)-ind(y)-1}(\overline{M}_{f}(x,y))$ for the moduli spaces
$M_{f}(x,y), \ x,y\in Crit(f)$. By using the description of $\partial\widehat{M}(x)$ and proceeding as in Lemma \ref{lem:repres_chain_syst},
we define by induction on
$ind_{f}(x)$ cubical chains
$\lambda_{x}\in S_{i}(\widehat{M}(x)),  x\in Crit_{i}(f)$
representing the fundamental class
of $(\widehat{M}(x),\partial\widehat{M}(x))$ and so that
\begin{equation}\label{eq:lambda_coherent}
\partial \lambda_{x}=\sum_{y} s_{xy}\times \lambda_{y}~.~
\end{equation}
Recall the map $\beta : \widehat{M}(x)\to P'(\widehat{M}(x)),\ y\to \overline{y}$ and
let $o':P'(\widehat{M}(x))\to P'(L)\to P'(\tilde{L})$
be induced by $\widehat{M}(x)\stackrel{o}{\longrightarrow} L\stackrel{q}{\longrightarrow} \tilde{L}$. Let also $\beta'=o'\circ\beta$.
We now define $\xi: \mathcal{R}_{\ast}\otimes \Z/2<Crit(f)>\to S_{\ast}(P'\tilde{L})$
by $\xi(x)=\beta'(\lambda_{x})$ for each $x\in Crit(f)$. It is clear
that this map respects the relevant filtrations. Due
to (\ref{eq:lambda_coherent}) it is also obvious
that  $\xi$ so defined is a chain map.

\subsubsection{$E^{2}(\xi)$ is an isomorphism.}\label{subsubsection:isomo_xi}
By construction $\xi$ is a morphism of $\mathcal{R}_{\ast}$-modules so it is
sufficient to show that $\xi'=E^{2}_{\ast,0}(\xi)$ is an isomorphism. For this purpose we notice that there is a natural evaluation map $\Upsilon:P'(\tilde{L})\to \tilde{L}$. By considering the map $id_{L}:L\to L$ as a trivial fibration
we see that $\Upsilon$ induces an isomorphism
$\Upsilon':E^{2}_{\ast,0}\to H_{\ast}(\tilde{L})$ and
that the composition $\Upsilon \circ \xi$ may be factored
as $$\mathcal{R}_{\ast}\otimes \Z/2<Crit(f)> \stackrel{r\otimes id}{\longrightarrow}
C'(f) \stackrel{u}{\longrightarrow} S_{\ast}(\tilde{L})$$
where $C'(f)$ is the chain complex defined as $C'(f)=S_{\ast}(\ast)\otimes \Z/2<Crit(f)>$ with differential $\partial'x=\sum_{y}r(s_{xy})y$ and with $r:S_{\ast}(\Omega\tilde{L})\to S_{\ast}(\ast)$ induced by the projection $\Omega'L\to \ast$ (as our cubical chains are normalized we have $S_{\ast}(\ast)=\Z/2$). Given that $r\otimes id$ induces an isomorphism $E^{2}_{\ast,0}(f)\to H_{\ast}(C'(f))$,
 our proof ends if we show that $u$ induces an isomorphism
in homology. Clearly, $u$ is defined by $u(x)=\Upsilon (\beta'(\lambda_{x}))$.
To prove that $u$ induces an isomorphism we proceed by induction. We let
$C'_{k}$ be the subcomplex of $C'(f)$ consisting of elements of degree at most $k$
and we assume that $u_{k}=u|_{C'_{k}}:C'_{k}\to S_{\ast}(\tilde{L}_{k})$ induces
an isomorphism in homology. For each $x\in Crit(f)$ the
chain $\lambda_{x}$ represents the fundamental class of
$(\widehat{M}(x),\partial\widehat{M}(x))$ and, moreover,
we have the homeomorpism indicated in (\ref{eq:local_hom}). This implies that
the couple of maps  $(u_{k+1},u_{k})$ induces an isomorphism
$H_{k+1}(C'_{k+1},C'_{k})\to H_{\ast}(\tilde{L}_{k+1},\tilde{L}_{k})$. By the
5-lemma this shows that $u_{k+1}$ induces an isomorphism and concludes the
proof of the theorem.
\end{proof}

\section{Applications.}\label{sec:appli}
As mentioned in the introduction the Serre spectral sequence has many non-trivial
differentials. Obviously, in view of Theorem \ref{theo:main} d. this shows that
there is an abundace of pseudo-holomorphic strips. In this section we make explicit
this statement and deduce a number of applications.

We consider here the same setting as before: $(M,\omega)$ is fixed as well as the Lagrangian submanifolds $L$ and $L'$
which are in general position  and satisfy
(\ref{eq:connectivity}) if not otherwise indicated. This condition is
dropped only in \S \ref{subsec:extension} where it will be replaced by requiring
that $L$ and $L'$ be hamiltonian isotopic and $\omega|_{\pi_{2}(M,L)}=0$.

We review shortly the other relevant notation to be used in this chapter.
In the presence of an almost complex structure which tames $\omega$, $J\in \mathcal{J}_{\omega}$,
we have the moduli spaces $\mathcal{M}_{J}(x,y)$, and $\mathcal{M}'_{J}(x,y)$ of,
respectively, un-parametrized and parametrized pseudo holomorphic strips joining
the intersection points $x, y\in L\cap L'$ as in
\S\ref{subsubsec:action_funct}. Moreover,  $\mathcal{M}_{J}'=\bigcup_{x,y}\mathcal{M}_{J}'(x,y)$, $\mathcal{M}_{J}=\bigcup_{x,y}\mathcal{M}_{J}(x,y)$.
In case $J$ is regular (which in the terminology
used before in the paper means that the pair $(0,J)$ is regular), then the Gromov compactifications
$\overline{\mathcal{M}}_{J}(x,y)$ satisfy (\ref{eq:Grom_comp}) and
the spectral sequence $EF(L,L')=EF(L,L';0)$
is defined as in \S\ref{subsubsec: ss_def}.
We denote by $\mathcal{J}_{reg}$ the set of those elements of $\mathcal{J}_{\omega}$ that are
regular. To simplify notation,
we drop the index $J$ if no confusion is possible. We recall that $CF(L,L')$ is the
usual Floer complex and $\mathcal{C}(L,L')$ is the extended complex
constructed in \S\ref{subsubsec: ss_def}.
Recall also from (\ref{eq:eval}) the
maps  $\Phi_{x,y}:\overline{\mathcal{M}}(x,y)\to \Omega'\tilde{L}$
where $\Omega'\tilde{L}$ is the space of Moore loops
and $\tilde{L}$ is obtained from $L$ by contracting to a point an
embedded path connecting the points in $L\cap L'$ as
in \S\ref{subsubsec:add_path_sp}. More explicitly, for
$u\in \overline{\mathcal{M}}(x,y)$, $\Phi_{x,y}(u)$ is the curve
traced by the strip $u$ on $L$ parametrized by the interval
$[0,\mathcal{A}(x)-\mathcal{A}(y)]$ (where $\mathcal{A}$ is the relevant action functional - see \ref{subsubsec:defor_pseudo}) viewed as a loop on $\tilde{L}$.  Finally, $\mathcal{R}_{\ast}=S_{\ast}(\Omega'\tilde{L})$.

\subsection{Global abundance of pseudo-holomorphic strips: loop space homology}

In all this subsection we work under the assumption that
(\ref{eq:connectivity}) is satisfied.
The point of view here is global. Roughly, we show that much of the algebraic topology of $\Omega L$ may be recovered from $\mathcal{M}$.
We fix the additional notation
$$\mathcal{K}=\bigcup_{x,y}Im(\Phi_{x,y})\subset \Omega'\tilde{L}$$
and we let $\widehat{\mathcal{K}}$ be the submonoid of $\Omega'\tilde{L}$ generated by $\mathcal{K}$.
Let $k:\widehat{\mathcal{K}}\to \Omega'\tilde{L}$ be the obvious inclusion.

\begin{cor} \label{cor:surj1} If $L$ and $L'$ are hamiltonian isotopic and $J\in \mathcal{J}_{reg}$, then
the algebra generated by the image of $$k_{\ast}: H_{\ast}(\mathcal{K};\Z/2)\to H_{\ast}(\Omega'\tilde{L}:\Z/2)$$  coincides with
$H_{\ast}(\Omega'\tilde{L}:\Z/2)$.
\end{cor}

\begin{proof}
In view of the definition
of the coefficients $a_{xy}\in S_{\ast}(\Omega'\tilde{L})$ as images through the
maps $\Phi_{x,y}$ of
a representing chain system for the moduli spaces $\mathcal{M}(x,y)$ (as described
in \S\ref{subsubsec:alg_constr}) we see that, in fact, $a_{xy}\in S_{\ast}(\widehat{\mathcal{K}})$. We denote this ring by $\mathcal{R}'$.
Therefore, the chain complex $\mathcal{C}(L,L')$ has coefficients in $\mathcal{R}'$ and, in general, in the construction of the spectral
sequence the ring $\mathcal{R}_{\ast}$ may be replaced by the smaller
ring $\mathcal{R}'$. This produces a spectral sequence $EF'(L,L')$ whose
$E^{2}$ term is $HF_{\ast}(L,L')\otimes H_{\ast}(\widehat{\mathcal{K}};\Z/2)$.
There is an obvious natural map$E(k)$ from this spectral sequence to the spectral sequence $EF(L,L')$. By Theorem \ref{theo:main} this last spectral sequence is isomorphic to the Serre spectral sequence $E^{r}_{p,q}$ of
$\Omega L\to PL\to L$. We shall prove that $k_{\ast}$ is surjective by induction.
Assume already shown that $k_{\ast}$ is surjective for $\ast< i$.
Let $a\in H_{i}(\Omega'\tilde{L};\Z/2)$. As the Serre spectral sequence converges
to the homology of an acyclic space this element viewed in $E^{2}_{0,i}$ has
to verify $[a]_{r}=d^{r}[b]_{r}\in E^{r}_{0,i}$ for some $r\geq 2$. We shall use again induction here over $r$. Therefore, assume that for all the elements $a'\in E^{2}_{0,i}$ which have the property $a'\in Im(d^{s})$ with $s<r$ we already know that $a'\in Im(k_{\ast})$. Let now
$b=\sum_{j}a_{j}\otimes x_{j}$ with
$a_{j}\otimes x_{j}\in \mathcal{R}'\otimes \Z/2<L\cap L'>$ (see \S
\ref{subsubsec: ss_def}). We know that $d^{0}(b)=0$. Therefore, the $a_{j}$'s
are cycles. As they are of degree stricly less than $i$ it follows that we may
assume that they are in the image of $k_{\ast}$. Together with the description
of the differential in $\mathcal{C}(L,L')$ this shows that
$[a]_{r}\in Im(E^{r}(k))$. But this shows that there exists $c\in E^{2}_{0,i}$
so that $[c]_{r}=0$ and $a+c\in Im(k_{\ast})$. However, our induction assumption
on $r$ implies that $c\in Im(k_{\ast})$ and so $a\in Im(k_{\ast})$ and this
concludes the proof.
\end{proof}

To state a closely related result, consider the (injective) map
\begin{equation}\label{eq:first_eval}
i:\mathcal{M}'\to \mathcal{P}(L,L'),\ i(u)(t)=u(0,t)
\end{equation}
and denote by $\tilde{\mathcal{M}}$ the (compact)  image of $i$. Consider the map
$e:\tilde{\mathcal{M}}\to \tilde{L}, e(u)=u(0,0)$.

\begin{cor}\label{cor:surj2} If $L$ and $L'$ are hamiltonian isotopic, $J\in \mathcal{J}_{reg}$, then
the map $e$ induces a surjective morphism
$$H_{\ast}(\Omega e):H_{\ast}(\Omega'\tilde{\mathcal{M}};\Z/2)\to
H_{\ast}(\Omega'\tilde{L};\Z/2)~.~$$
\end{cor}
\begin{proof}This is immediate from the previous corollary as $k:\widehat{\mathcal{K}}\to \Omega'\tilde{L}$ factors through  $\Omega e:\Omega'\tilde{\mathcal{M}}\to \Omega'\tilde{L}$.
\end{proof}

\begin{rem}\label{rem:fl_ho}
{\rm a. There exist many examples of maps $w:X\to Y$ so that one of $H_{\ast}(w)$, $H_{\ast}(\Omega w)$ is surjective but the other is not so the result at
\ref{cor:surj2} is non trivial.

b. As mentioned in the introduction it has been proven by Floer \cite{Fl} and Hofer \cite{Ho} that $H_{\ast}(e)$ is also surjective even in the
degenerate case. It is likely
that the surjectivity at  \ref{cor:surj2} remains true in the degenerate case.

c. Both corollaries may be viewed as stability results for moduli
spaces of pseudo-holomorphic strips: they are quite immediate for
negative gradient flow lines of Morse-Smale functions and therefore
they are unsurprising when the isotopy $\phi$ is small, however the fact that
the same properties are preserved even when making $\phi$ large are quite
non-trivial.}
\end{rem}

We end this subsection with a different, simple topological consequence.
Its content is that, generically, due to the presence of pseudo-holomorphic strips
the form $\omega$ ``sees'' much of $\pi_{1}(L\cup L')$.

There exists a generic class  $\mathcal{L}$ of lagrangians
$L'$ which are not only transversal and hamiltonian isotopic to $L$ as assumed
till now but also have the property that the abelian group generated by the obvious map $\mathcal{A}_{L,L'}|:L\cap L'\to \R$ is of maximal rank ($=\#(L\cap L'$)).
In other words, the action functional $\mathcal{A}_{L,L'}$ takes different values on
each of the intersection points $L\cap L'$ and, moreover, these values are linearly
independent.

Let $\Pi_{\omega, L, L'}:\pi_{2}(M,L\cup L')\to \R$ be defined by $I_{\omega}(u)=\int_{D^{2}}u^{\ast}\omega$.

\begin{cor}\label{cor:pi_2}
For $L'\in \mathcal{L}$, the image of $\Pi_{\omega, L, L'}$ is an abelian
group of rank at least $dim(\overline{H}_{\ast}(L;\Z/2))$ (where $\overline{H}(-)$
denotes reduced homology).
\end{cor}

\begin{proof} We fix $J\in \mathcal{J}_{reg}$. We will prove that there is a set  $\{x_{1}, x_{2},\ldots ,x_{m}\}\subset L\cap_{\eta} L'$ with $m=dim(\overline{H}_{\ast}(L;\Z/2))$ such that
for each $x_{i}$ there exists some $y_{i}$ and $u_{i}\in \mathcal{M}(x_{i},y_{i})\not=\emptyset$.
Given the definition of $\mathcal{L}$ this suffices to show the satement because
for $u_{i}\in \mathcal{M}'(x,y)$ we have $\mathcal{A}_{L,L'}(x_{i})-\mathcal{A}_{L,L'}(y_{i})=
\int_{\R\times [0,1]}u^{\ast}_{i}\omega$  which shows that the values $\Pi_{\omega,L,L'}(u_{i})$ are linearly independent (over $\Z$).

To simplify notation we shall say that
$x\in L\cap L'$ is a strip origin if there exists $y$ such that $\mathcal{M}(x,y)\not=\emptyset$.
We now let $a_{1},a_{2},\ldots a_{m}$ be a basis for $\overline{H}_{\ast}(L,\Z/2)$.
We pick chains $z_{i}\in CF_{\ast}(L,L')$
representing respectively the classes $a_{i}$. We write $z_{i}=\sum x^{i}_{j}$
where $x^{i}_{j} \in L\cap_{\eta} L'$. Notice that if a point $x\in L\cap L'$
of positive degree is not a strip origin, then its differential in the Floer complex
is null so its homology class $[x]\in HF_{\ast}(L,L')$ is well-defined. Moreover,
this homology class has to be null because $[x]$ viewed as an element
of $EF^{2}_{p,0}(L,L')$ survives to $E^{\infty}$. Because of this we may assume
that each $x^{i}_{j}$ appearing in the expression of the chains $z_{i}$ is a strip
origin. This implies the claim because if there are strictly less then $m$ distinct points among the $x^{i}_{j}$'s, then the family $\{z_{i}\}$ is linearly dependent which contradicts the fact that
the family $\{a_{i}\}$ is linearly independent.
\end{proof}

\subsection{Local pervasiveness of pseudo-holomorphic strips.}

With our machinery it is not hard to deduce that
through each point of $L\backslash L'$ passes some pesudo-holomorphic strip
(see for example Corollary \ref{cor:many_strips} below; this also follows from the results of Hofer and Floer mentioned in Remark \ref{rem:fl_ho}). The point of view here is however slightly different: what most interests us is to restrict the type of these strips that ``fill" $L$.   We again work under the assumptions
at (\ref{eq:connectivity}) .

\

We start with a useful, purely algebraic consequence of the construction of the spectral sequence $EF(L,L')$.
Assume that $g:L\to X$ is a continuous map.
Let $\Omega X\to E_{g}\to L$ be the pullback fibration
$g^{\ast}(\Omega X\to PX\to X)$. There is an obvious
map of fibrations which is induced by $g$:
$$\xymatrix@-5pt{
\Omega{L}\ar[r]^{\Omega g}\ar[d]& \Omega X\ar[d] \\
PL\ar[d] \ar[r]^{\overline{g}} & E_{g}\ar[d]\\
L\ar[r]_{id} & L}$$
This construction may be performed also by using Moore loops
instead of usual ones and we may as well replace $L$
by $\tilde{L}$. Therefore, if we denote the ring $S_{\ast}(\Omega'X)$ by $\mathcal{R}''$ we have a change of
coefficients map $g^{\#}:\mathcal{R}\to \mathcal{R}''$.
We may obviously use this map to define a complex $\mathcal{C}_{X}(L,L')$ as in
\S\ref{subsubsec: ss_def} which is obtained from $\mathcal{C}(L,L')$ by this change of coefficients. There is also
an associated spectral sequence $EF_{X}(L,L')$ into which $EF(L,L')$ maps by the map induced by $g^{\#}$. The properties of $EF_{X}(L,L')$ parallel those of $EF(L,L')$ and have the same proofs. In particular, property d. becomes:

\begin{cor}\label{cor:alg} If $L$ and $L'$ are hamiltonian isotopic and $J\in \mathcal{J}_{reg}$, then
the spectral sequence $EF_{X}(L,L')$ is defined and its terms of order greater
or equal than $2$ are isomorphic up to translation
to the corresponding terms of the Serre spectral
sequence of the fibration $\Omega X\to E_{g}\to L$.
\end{cor}

For the next corollary we shall assume that $L$ and $L'$ are hamiltonian isotopic and $J\in \mathcal{J}_{reg}$.
Recall that we have an isomorphism up to translation between
$HF_{\ast}(L,L')\approx H_{\ast}(L,\Z/2)$. To simplify notation we shall assume this isomorphism to be degree preserving.
We shall denote by $[1]\in HF_{0}(L,L')$ and $[L]\in HF_{n}(L,L')$
the generators of the respective homology groups.
There is an obvious evaluation map $\mathcal{E}:\mathcal{M}'\to L$
which is defined by $\mathcal{E}(u)=u(0,0)$ (it verifies $\mathcal{E}=e\circ i$
with $i,e$ as in (\ref{eq:first_eval})).
For two elements $x=\sum_{i}c_{i}x_{i}\in CF(L,L')$, $y=\sum_{j}d_{j}x_{j}\in CF(L,L')$ we let $$\mathcal{R}(x,y)=
\bigcup_{c_{i}\not=0,\ d_{j}\not= 0}\mathcal{E}(\mathcal{M}'(x_{i},x_{j}))~.~$$

Let $a$ be a
representative of the fundamental class $[L]$. Consider also
an element $b\in CF_{0}(L,L')$ which is the sum of all the intersection points
which appear with non-vanishing coefficients in some  representatives of
the homology class $[1]$.

\begin{cor}\label{cor:many_strips}
In the setting above, the set $\mathcal{R}(a,b)$ is dense in $L$. In particular, each $x\in L\backslash L'$ belongs to some pseudo-holomorphic strip of Maslov index at most $n$.
\end{cor}

\begin{proof}
Assume that the image of $\mathcal{R}(a,b)$ avoids a small  open disk $D\subset L$.
This implies that $\overline{\mathcal{R}(a,b)}\subset L\backslash D$.
We may assume that the path $w$ used to construct $\tilde{L}$ out of $L$
intersects $D$ in just a point.
We intend to apply Corollary \ref{cor:alg} to the map $g:L\to D/\partial D=S^{n}$
defined by contracting the closure of the complement of $D$ to a point. The basic
algebraic fact that we will be using is that in the Serre spectral sequence
of the induced fibration $\Omega S^{n}\to E_{g}\to L$ the differential
$d^{n}$ verfies $d^{n}[L]=[1]\otimes l$ where $l$ it the homology class ot the bottom
 sphere $S^{n-1}\hookrightarrow \Omega S^{n}$ (this inclusion is the adjoint of the identity). By Corollary \ref{cor:alg} the same statement is true for the spectral sequence $EF_{S^{n}}(L,L')$. Let $\mathcal{C}'_{S^{n}}(L,L')$ be the subcomplex of $\mathcal{C}_{S^{n}}(L,L')$ which is generated
by $\{x\in L\cap L' : \mu(a,x)>0\}\cup \{a\}$. The index filtration defines
a spectral sequence $EF'_{S^{n}}(L,L')$ which obviously maps into
$EF_{S^{n}}(L,L')$ by the map $E(i)$ induced by the inclusion
$i:\mathcal{C}'_{S^{n}}(L,L')\hookrightarrow \mathcal{C}_{S^{n}}(L,L')$.
As $a$ represents the fundamental class in $CF_{n}(L,L')$ it follows that
we also have $d^{n}([a])=[1]\otimes[l]$ in $EF_{S^{n}}'(L,L')$. Consider also the
subcomplex $\mathcal{C}''_{S^{n}}(L,L')$ which is generated by
$\{x\in L\cap L' : \mu(b,x)> 0\}$. The
quotient complex $\overline{\mathcal{C}}_{S^{n}}(L,L')=\mathcal{C}'_{S^{n}}(L,L')/
\mathcal{C}''_{S^{n}}(L,L')$ is well-defined and it admits an obvious filtration
such that the quotient map $p:\mathcal{C}'_{S^{n}}(L,L')\to \overline{\mathcal{C}}_{S^{n}}(L,L')$ preserves filtrations. Therefore, it
induces a morphism of spectral sequences $E(p):EF'_{S^{n}}(L,L')\to \overline{EF}_{S^{n}}(L,L')$. We notice that $(E(p))^{2}_{0,0}$ is
injective and so, in $\overline{EF}_{S^{n}}(L,L')$, we have again $d_{S^{n}}([a])=[1]\otimes l$. But, if in  $\mathcal{C}(L,L')$ we have $da=\sum_{i} k_{ay}\otimes y$, then the differential of $a$
 in $\mathcal{C}_{S^{n}}(L,L')$ is given by
$d_{S^{n}}a=\sum_{i} g^{\#}(k_{ay})\otimes y$ and as by assumption $\overline{\mathcal{R}(a,b)}$ avoids $D$
it follows that all the critical points $y$ which appear with non-zero coefficients in the expression of some representative of  $1$ have
$g^{\#}(k_{ay})=0$ which contradicts $d_{S^{n}}([a])=[1]\otimes l\not=0$.
  \end{proof}

For the next result recall from the introduction that, for any two Lagrangians $L,L'\subset M$ (not necessarily transversal), we  define as in \cite{Chek}, \cite{Oh} the isotopy energy of $L$ and $L'$  by $$\nabla (L,L')=\inf_{\phi\in \mathcal{H}, \phi(L)=L'} ||\phi||_{H}$$
where $\mathcal{H}$ is the group of Hamiltonian diffeomorphisms with
compact support and, as before, $|| - ||_{H}$ is Hofer's energy (see (\ref{eq:hofer_norm})). In case $L$ and $L'$ are not isotopic we take $\nabla(L,L')=\infty$.
It is easy to see that $\nabla(-,-)$ is symmetric and verifies
the triangle inequality.  Moreover,
it has been shown by Chekanov \cite{Chek} following earlier work by Oh \cite{Oh}
that $\nabla(-,-)$ is non-degenerate for arbitrary compact lagrangians in tame
symplectic manifolds thus providing a metric on any (Hamiltonian) isotopy equivalence class of Lagrangians.

\begin{cor}\label{cor:strips_energy}
Suppose that $L,L'\subset M$ are transversal, simply-connected lagrangians
embedded in the symplectic manifold $(M,\omega)$ with
$\omega|_{\pi_{2}(M)}=c_{1}|_{\pi_{2}(M)}=0$. If $L$ and $L'$ are hamiltonian
isotopic, then for any almost complex structure $J\in\mathcal{J}_{\omega}$ and for any
point $x\in L\backslash L'$ there exists a $J$-pseudo-holomorphic strip
$u: \R\times [0,1]\to M,\ u(\R,0)\subset L, \ u(\R,1)\subset L'$ so that
 $x\in Im(u)$, $$\int_{\R\times [0,1]}u^{\ast}\omega\leq \nabla(L,L')~.~$$ Moreover,
 when $J\in\mathcal{J}_{reg}$, there is a strip  $u$ as above which also verifies
 $\mu(u)\leq n$.
\end{cor}

\begin{proof}
The Gromov compactness
theorem applies to sequences $u_{n}$ of $J_{n}$-holomorphic
curves where $J_{n}\in \mathcal{J}_{\omega}$ is a sequence of almost complex structures which converges towards another almost complex structure $J\in\mathcal{J}_{\omega}$ \cite{McSal}. As any almost complex structure belonging
to the set $\mathcal{J}_{\omega}$ may be viewed as the limit of a sequence of regular almost complex structures, this implies
that it is sufficient to prove the statement when $J\in \mathcal{J}_{reg}$ and
so we assume this for the rest of the proof.

Recall the definition of the energy $E_{L,L',H}(u)$
of the elements of $u\in\mathcal{M}'_{L,L',J,H}$
from formula (\ref{eq:energy}). If $u\in \mathcal{M}'_{L,L',J,0}$ we have
$E_{L,L'}(u)=\int_{\R\times [0,1]}u^{\ast}\omega$. Let
$$\mathcal{M}^{a}_{L,L',J,H}=\{u\in \mathcal{M}'_{L,L',J,H} : E_{L,L',H}(u)\leq a, \mu(u)\leq n\}~.~$$ Our main interest is in $\mathcal{M}^{a}=\mathcal{M}^{a}_{L,L',J,0}$.
We now assume that $\phi$ is a hamiltonian isotopy such that $L'=\phi (L)$.
Given that the set of energies of the elements in $\mathcal{M}'$ is discrete
it is sufficient for our statement to show that the set $\mathcal{E}(\mathcal{M}^{||\phi||_{H}})$ is dense in $L$ and this is what we
shall show next.

Let $f:L\to \R$ be a Morse function with a single
maximum $P$ and a single minimum $Q$ and we let
$f(P)=\epsilon >0$, $f(Q)=0$. We pick a riemannian metric $\alpha$ so that
the pair $(f,\alpha)$ is Morse-Smale. Let $\overline{f},L_{f}, J^{f}$ be as in
\S\ref{subsubsec:Morse_fl_strips}. As in \S\ref{subsubsec:red_Morse}
we have a homeomorphisms of moduli spaces inducing an identification
of chain complexes:
$$\overline{l}_{f}: C^{f,\alpha}\to \mathcal{C}^{J^{f}}(L,L_{f};0)~.~$$
Notice also that the naturality results in \S\ref{subsubsec:naturality}
show that $\mathcal{A}_{L,L_{f}}(x)=f(x)$ for all $x\in Crit(f)=L\cap L_{f}$.
Fix some small $\delta >0$. By taking $\epsilon$ sufficiently small we see
that there exists a Hamiltonian $G:[0,1]\times M \to \R$ which is constant outside a compact set and so that $L_{f}=\phi_{1}^{G}(L')$ and $$Var(G)=\sup_{x,t}G(t,x)-\inf_{x,t}G(t,x)\leq ||\phi||_{H}+\delta~.~$$
By adding an appropriate constant to $G$ we may assume that $\inf_{x,t}G(t,x)=\epsilon$.
>From \S\ref{eq:natural_id} we see that we have homeomorphisms of moduli spaces
$b_{G}:\mathcal{M}_{L,L',J,0}\to \mathcal{M}_{L,L_{f},J'}$ and an identification of chain complexes (which
is action preserving - see \S\ref{subsubsec:naturality}):
$$\overline{b}_{G}:\mathcal{C}^{J}(L,L';0)\to \mathcal{C}^{J'}(L,L_{f};G)$$
where $(\phi_{1}^{G})^{\ast}(J')=J$.

Thus, in view of the definition of $b_{G}$, the image of $\mathcal{M}^{a}$ by the evaluation map coincides with the image
of $\mathcal{M}^{a}_{L,L_{f},J',G}$ and, therefore, to show the claim
it is enough to
prove that $\mathcal{E}(\mathcal{M}^{a}_{L,L_{f},J',G})$ is dense
(recall also that the energy and the action are related by formula
(\ref{eq:energy_action})).

Let $\Delta=Var(G)$. We consider also the Hamiltonian
$G'=G-\Delta-2\epsilon$. We can then define monotone homotopies
from $G$ to $0$ and from $0$ to $G'$.
Let $CF=CF(L,L_{f};G)$ be the Floer complex and, similarly,
let $CF'=CF(L,L_{f};G')$ and $C=CF(L,L_{f};0)=C^{f,\alpha}$.

It follows from
equation (\ref{eq:monotone}) and as in the proof of Theorem \ref{theo:comparison}  that the monotone homotopies mentioned above induce morphisms of chain complexes:
$$\mathcal{V}:\mathcal{C}^{J'}(L,L_{f};G)\to \mathcal{C}^{J^{f}}(L,L_{f};0)\ , \ V:CF\to C$$
and
$$\mathcal{W}:\mathcal{C}^{J^{f}}(L,L_{f};0)\to \mathcal{C}^{J'}(L,L_{f}; G')\ , \ W:C\to CF'$$
which are not action increasing.

For an element $x\in CF$,
$x=\sum_{i}c_{i}x_{i}$ with $c_{i}\in \Z/2$ and $x_{i}\in I(L, L_{f};G)$ we let $\overline{\mathcal{A}}_{G}(x)=\max_{c_{i}\not=0}
\mathcal{A}_{L,L_{f},G}(x_{i})$ and $\underline{\mathcal{A}}_{G}(x)=
\min_{c_{i}\not=0}\mathcal{A}_{L,L_{f},G}(x_{i})$.
We define similarly $\overline{\mathcal{A}}_{G'}(x')$
for $x'\in CF'$. Let now $a=\sum_{i}r_{i}a_{i}\in CF$ be a chain
representative of the fundamental class $[L]$ such that $\overline{\mathcal{A}}_{G}(a)$ is minimal among all such representatives.
Notice at this point that the complexes $\mathcal{C}(L,L_{f};G')$
and $\mathcal{C}(L,L_{f};G)$ are canonically identified (and similarly for $CF$ and $CF'$). We distinguish elements
of the two complexes by indicating this identification
as $ \mathcal{C}(L,L_{f};G)\ni x \to x'\in \mathcal{C}(L,L_{f};G')$ and
we clearly have $\mathcal{A}_{G}(x)=\mathcal{A}_{G'}(x')+\Delta+2\epsilon$. Let $c'\in CF'$ be defined by $c'=W(P)$ (we recall that $P$ is the unique maximum point of $f$). Then $c$ is also a representative of $[L]$ and, therefore, $\overline{\mathcal{A}}_{G}(c)\geq \overline{\mathcal{A}}_{G}(a)$. At the same
time, as $W$ is not action increasing, we have $\overline{\mathcal{A}}_{G'}(c')\leq \mathcal{A}_{L,L^{f}}(P)=f(P)=\epsilon$ which means that $\overline{\mathcal{A}}_{G}(a)\leq \Delta +3\epsilon\leq ||\phi||_{H}+\delta+3\epsilon$.

For two elements $x,y\in CF$ with
$x=\sum c_{i}x_{i}$, $y=\sum d_{i}y_{i}$ we put
$\mathcal{R}_{G}(x,y)=\bigcup_{i,j} \mathcal{E}(\mathcal{M}'_{L,L_{f},G}(x_{i},y_{j}))$. For the
proof of the corollary it suffices to show that the element
$b=\sum_{i}b_{i}\in CF$ defined as the sum of all the generators $b_{i}$ of $CF$ with $\mu(b_{i})=0$, $\mathcal{A}_{G}(b_{i})\geq 0$ has the property that the
set $\mathcal{R}_{G}(a,b)$ is dense in $L$. Indeed, each element
in $\mathcal{M}'_{L,L_{f},G}(a_{i},b_{j})$ has energy bounded by
$\overline{\mathcal{A}}_{G}(a)-\underline{\mathcal{A}}_{G}(b)\leq ||\phi||_{H}+\delta+3\epsilon$ and as we may take $\delta$ and $\epsilon$ arbitrarily
small and because the possible energy values form a discrete set this implies the claim.

To show the density of $\mathcal{R}_{G}(a,b)$ we proceed in a way similar to that
of Corollary \ref{cor:many_strips}. Therefore, we assume that this set is disjoint from a disk $D\subset L$ and we consider the associated map
$g:L\to D/\partial D=S^{n}$ and the associated change of coefficients
$g^{\#}:\mathcal{R}\to \mathcal{R}''$. We use the same conventions as in Corollary
\ref{cor:alg} and to shorten notation we let $\mathcal{C}_{1}=\mathcal{C}^{J'}_{S^{n}}(L,L_{f};G)$,
$\mathcal{C}(f)=\mathcal{C}^{J^{f}}_{S^{n}}(L,L_{f};0)$.
Consider the sub-complex $\mathcal{C}_{0}\hookrightarrow \mathcal{C}_{1}$
which is generated by the elements $x\in I(L,L_{f};G)$ such that $\mathcal{A}_{G}(x)<0$ and consider also
 the quotient $\mathcal{C}_{2}=\mathcal{C}_{1}/\mathcal{C}_{0}$. We notice that due to the monotonicity of the homotopy inducing $\mathcal{V}$ the map $\mathcal{V}_{S^{n}}$ factors as $\mathcal{C}_{1}\stackrel{\mathcal{V}'}{\longrightarrow}\mathcal{C}_{2}
\stackrel{\mathcal{V}''}{\longrightarrow}\mathcal{C}(f)$ where the first map
is the passage to quotient. Both $\mathcal{V}'$ and $\mathcal{V}''$ respect
filtrations and thus they induce spectral sequences
morphisms $EF(\mathcal{V}'): EF_{S^{n}}(L,L_{f};G)\to \overline{EF}$
(where $\overline{EF}$ is the spectral sequence induced by the degree filtration
on $\mathcal{C}_{2}$) and $EF(\mathcal{V}''):\overline{EF}\to E_{S^{n}}(f)$.
The composition of these two morphism being an isomorphism for $r\geq 2$
and as in $E_{S^{n}}(f)$ we have $d^{n}[L]=[1]\otimes [l]$ (with $l$ the class of the
bottom sphere in $H_{n-1}(\Omega S^{n})$ it follows that $d^{n}[a]=k\otimes [l]$
with $k\not=0$ in $\overline{EF}$. But the fact that
$\mathcal{E}(\mathcal{R}_{G}(a,b))$ avoids $D$ implies that all the coefficients $g_{i}\in \mathcal{R}$ of the $b_{i}$'s in the expression of the differential of $da\in \mathcal{C}^{J'}(L,L_{f};G)$ have the property that $g^{\#}(g_{i})=0$
and thus we arrive at a contradiction.
\end{proof}

\begin{rem}\label{rem:Maslov_ind}{\rm
a. There is a certain overlap between Corollary \ref{cor:strips_energy} and Corollary \ref{cor:many_strips}. However, the choice of the element $a$ in this last corollary is less restrictive than in the proof of  \ref{cor:strips_energy}.

b. Given a manifold $N^{n}$, the degree one map $N\to S^{n}$, produced by collapsing onto $D^{n}/\partial D^{n}$ where $D^{n}$ is a closed disk in $N$, is the simplest
possible example of a Thom-Pontryaguin map. From this point of view, Corollaries
\ref{cor:many_strips} and \ref{cor:strips_energy} are truly immediate consequences of the main theorem (compared to this theorem, the only new idea appears in
the proof of \ref{cor:strips_energy}).

c. There exist some other methods to produce Floer orbits joining the ``top and bottom classes'' in the Floer complex. An interesting such approach is provided
by Schwarz \cite{Sch3},\cite{Sch2} and is based on the pair of pants product. This suggests a
relation between our invariant and this product. Obviously, such a relation
is also to be expected for purely topological reasons.}
\end{rem}

\subsection{Non-squeezing.}\label{subsec:non_squeeze}

In this subsection we shall prove a number of geometric
consequences of the previous results.

We consider two closed Lagrangians $L, L'\subset M$. We assume that $L$ and $L'$ intersect transversely
and let $J\in\mathcal{J}_{\omega}$.  We do not assume for now that $L$ and
$L'$ are hamiltonian isotopic.
Recall from, \S\ref{subsubsec:non_squeeze} and from equation (\ref{eq:ham_distance}) the following notation
in which the areas are computed with respect
to the riemannian metric induced by $J$.
\begin{itemize}
\item[-] $\mathcal{S}(x,y)$ is the set of $C^{\infty}$ strips joining $x,y\in L\cap L'$.
\item[-] $a_{L,L'}(x,y)$ is the infimum of the areas of the
strips in $\mathcal{S}(x,y)$.
\item[-]$a_{k}(L,L';J)$ is the minimal area of a pseudo-holomorphic strip of index $k$.
\item[-]$A_{k}(L,L';J)$ is the maximal area of such a strip (these numbers are taken to be infinite if no such pseudo-holomorphic strips exist).
\item[-] $\epsilon(L,L';J)=\epsilon(L,L';J,0)$ is the minimal energy of some element in $\mathcal{M}'_{L,L',J}$ (this number is taken to be infinite if $\mathcal{M}'_{L,L',J}$ is void). We define another pair of associated numbers
 by $\overline{\epsilon}(L,L')=\sup_{J\in\mathcal{J}_{reg}}\epsilon(L,L';J)$
and $\underline{\epsilon}(L,L')=\inf_{J\in\mathcal{J}_{reg}}\epsilon(L,L';J)$.
\item[-] $\delta(L,L')$ is the maximal radius of a standard symplectic ball $B(r)$ which is symplectically embedded in $M$ with an image disjoint from $L'$ and
so that the image of $\R^{n}\cap B(r)$ is included in $L$.
\end{itemize}

These notions are well-defined independently of the connectivity
conditions in (\ref{eq:connectivity}) but for the remainder of this
subsection we assume that this condition is satisfied.

\begin{rem}\label{rem:var_area}{\rm a. The area of an element of $\mathcal{M}'_{L,L',J}$ coincides with its energy and also coincides with its symplectic area and, moreover, if $\mathcal{M}_{L,L',J}(x,y)\not=\emptyset$, then $a_{L,L'}(x,y)= E_{L,L',J}(u), \ \forall u\in \mathcal{M}_{L,L',J}'(x,y)$.

b. As has been observed by Fran\c{c}ois Lalonde, the value of $\delta (L,L')$ is not changed if to its definition we add the condition that the image of the center of $B(r)$ be equal to some fixed point $x\in L$. This is because if  $e:B(r)\to M$ is an embedding as required and such that
$e(0)=y$, then we may find a hamiltonian isotopy that carries $y$ to $x$, which is supported in the neighbourhood of a path joining $y$ to $x$ inside $L$ and which send $L$ to itself. Indeed, as $L$ is compact, by dividing the path joining $y$ to $x$ into small enough pieces, we may assume that both $x$ and $y$ belong to a
standard coordinate chart ressembling $(\C^{n},\R^{n})$, $x,y\in \R^{n}$
and in this situation the problem is trivial.

c. We have the obvious inequalities: $\underline{\epsilon}(L,L')\geq \min\{\mathcal{A}_{L,L'}(x)-\mathcal{A}_{L,L'}(y) : x,y\in L\cap L', \mathcal{A}_{L,L'}(x)>\mathcal{A}_{L,L'}(y)\}>0~.~$}
\end{rem}

>From Corollary \ref{cor:many_strips} and Remark \ref{rem:var_area} c. we deduce that, if, in addition, $L$ and $L'$ are Hamiltonian isotopic and $J\in\mathcal{J}_{reg}$ then:
\begin{equation}\label{eq:inequality1}\infty > A_{n}(L,L';J)\geq a_{n}(L,L';J)\geq \epsilon(L,L';J)\geq \underline{\epsilon}(L,L')>0~.~
\end{equation}

Under the same assumptions, we obtain from Corollary \ref{cor:strips_energy}
that \begin{equation}\label{eq:non_vanishing_trans}
\nabla(L,L')\geq a_{n}(L,L';J)~.~
\end{equation}

Stronger inequalities follow.

\begin{cor}\label{cor:non-squeeze}
Assume $L$, $L'$ are two simply-connected Lagrangian submanifolds of $(M,\omega)$ and suppose $\omega|_{\pi_{2}(M)}=0=c_{1}|_{\pi_{2}(M)}$.
\begin{itemize}
\item[i.] If $L$ and $L'$ intersect transversely and $J\in \mathcal{J}_{reg}$, then we have $$\nabla(L,L')\geq \frac{\pi}{2}\ \delta(L,L')^{2}~,~$$ $$A_{n}(L,L';J)\geq \frac{\pi}{2}\ \delta(L,L')^{2}~.~$$
\item[ii.] Additionally, suppose that $L''$ is another Lagrangian transversal to $L$ and assume $J'$ is another almost complex structure
such that $\mathcal{C}^{J'}(L,L'')$ is defined.
If $\mathcal{C}^{J'}(L,L'')$ does not admit $\mathcal{C}^{J,\eta}(L,L')$ as
a retract(for some choice of $\eta$), then $$\nabla(L',L'')\geq \epsilon(L,L';J)/4~.~$$
In particular, the energy needed to diminish the number of intersection points between $L$ and $L'$ by a hamiltonian isotopy is at least $\underline{\epsilon}(L,L')/4$.
\end{itemize}
\end{cor}

\begin{proof}
The proof of i. is a rapid consequence of Corollary \ref{cor:strips_energy}
combined with an argument classical in symplectic topology since the work of Gromov. We assume $L$ and $L'$ hamiltonian isotopic by an isotopy $\phi$.
Fix a (standard) ball $B(r)$ and an embedding $e:B(r)\to M$ as in
the definition of  $\delta(L,L')$. Thus $e(B(r))\cap L'=\emptyset$, $e^{-1}(L)=\R^{n}\cap B(r)$, $e(0)=x\in L$.
 Fix a small $\delta>0$. There exists an almost complex structure
$J\in \mathcal{J}_{\omega}$ on $M$ such that $e^{\ast}J|_{B(r-\delta)}=J_{0}$ where $J_{0}$ is the canonical almost complex structure on $B(r)$ (in fact, $J$ is constructed by extending the push forward of $J_{0}$).
It results from Corollary \ref{cor:strips_energy}
that there exists a  $J$-pseudo-holomorphic strip $u\in\mathcal{M}'_{L,L',J}$
that passes through $x$ and verifies $\int_{\R\times [0,1]}u^{\ast}\omega
\leq ||\phi||_{H}$, $\mu(u)\leq n$. We now consider $v=e^{-1}(u\cap B(r-\delta))$. This is a $J_{0}$-pseudo-holomorphic curve in $B(r-\delta)$ whose boundary lies on  $\partial B(r-\delta)\cup \R^{n}$ and whose area is bounded from above by $||\phi||_{H}$. By analytic continuation, this
curve extends to a $J_{0}$-pseudo-holomorphic curve $\overline{v}$ whose boundary is contained in $\partial B(r-\delta)$, which contains $0$ and whose
area is the double of that of $v$.
By the classical isoperimetric inequality we get that the area
of $\overline{v}$ is at least $\pi(r-\delta)^{2}$. Thus the area
of $u$ is at least $\pi(r-\delta)^{2}/2$ and this shows
$||\phi||_{H}\geq \pi r^{2}/2$ and implies the inequalities at the first point.

Point ii. is a reformulation of Corollary \ref{cor:comp_lagran}.
\end{proof}

\begin{rem}\label{rem:nondeg_nabla} {\rm  Point i. implies that $\nabla(-,-)$ is nondegenerate
(in our setting): indeed, if $L\not=L''$ then we may find a small symplectic ball $B(r)$ embedded in $M$ as in the definition of $\delta(-,-)$. If $L$ and $L''$ are not transversal we may perturb $L''$ to a lagrangian $L'''$ which is transversal to $L$ without touching $B(r)$. Therefore, using the triangle inequality, we have $\nabla(L,L'')\geq \pi r^{2}/2 -\delta(L'',L''')$. As the perturbation of $L''$ can be made as small as needed we get that $\nabla(L,L'')>0$.
 }

\end{rem}

We conjecture that the inequality on the left in Corollary \ref{cor:non-squeeze} i.
is true for any pair of closed lagrangian submanifolds of a closed symplectic manifold. We shall prove it in the next section under weaker assumptions than (\ref{eq:connectivity}).

\subsection{Relaxing the connectivity conditions.} \label{subsec:extension}
There are some obvious extensions of our construction - for example by using the
orientations of the various moduli spaces involved we could use $\Z$ coefficients
instead of $\Z/2$ coefficients.

The purpose
of this subsection is to discuss a different extension which is quite useful.
This concerns replacing the rather stringent connectivity requirements in (\ref{eq:connectivity}) by the assumption that $L$ and $L'$ are hamiltonian isotopic and
\begin{equation}\label{eq:connect2}
\omega|_{\pi_{2}(M,L)}=0~.~
\end{equation}

As we shall see, adapting the construction of our invariant to this setting turns out to be reasonably straightforward even if the result of the construction is less elegant than before (reason that have made us postpone this variant of the construction till this moment). As a consequence, the proofs of Corollaries \ref{cor:strips_energy} and \ref{cor:non-squeeze} i. remain valid in this case and therefore we obtain the
following strengthening of these two results.

\begin{cor}\label{cor:general_strips}
Assume that $L^{n}\subset (M^{2n},\omega)$ is a closed Lagrangian submanifold such that $\omega|_{\pi_{2}(M,L)}=0$. If $L'\subset M$ is a second Lagrangian in the
same Hamiltonian isotopy class as $L$, then for each point $x\in L$
and each almost complex structure structure $J\in \mathcal{J}_{\omega}$ there
exists a pseudo-holomorphic strip $u\in \mathcal{M}'_{L,L',J}$ such that
$x\in Im(u)$, $\int u^{\ast}\omega \leq \nabla(L,L')$ and, additionally, when $J\in \mathcal{J}_{reg}$,
$\mu(u)\leq n$.
In particular, $\nabla(L,L')\geq \frac{\pi}{2}\delta(L,L')^{2}$ and, moreover, if $J\in\mathcal{J}_{reg}$
we also have $ \frac{\pi}{2}\delta(L,L')^{2}\leq A_{n}(L,L')$.
\end{cor}

\begin{rem}\label{rem:energy_capacity_disj}{\rm
 As in Remark \ref{rem:nondeg_nabla}, the inequality in the last corollary
recovers (under the assumption $\omega|_{\pi_{2}(M,L)}=0$)
 Chekanov's result claiming  that $\nabla (L,L')$ is a distance on the hamiltonian isotopy class of $L$. Moreover, the same inequality
also implies that for any symplectic
manifold $M$ with $\omega|_{\pi_{2}(M)}=0$, the disjunction energy of a subset
$A\subset M$ is greater than half the Gromov radius of $A$, a result proven
for all symplectic manifolds by Lalonde and McDuff \cite{Lalonde_McDuff}.
Indeed, for this last result, assume that $\phi$ is a hamiltonian
isotopy of $M$ that disjoins $A$ from itself. Suppose also that $A$ contains
a standard ball $B(r)$ or radius $r$. We may assume that $M'=graph(\phi)\subset (M\times M,\omega\oplus -\omega)$ intersects the diagonal $\Delta$ transversely.
The fact that $\phi$ disjoins $A$ from itself implies
that $B(r)\times B(r)$ is disjoint from $M'$. Given that $B(r)\times B(r)$
contains a standard $4n$-dimensional ball $B'(r)$ of radius $r$ centered on $\Delta$
and such that $B'(r)\cap \Delta$ is included in the image of $\R^{2n}$,
the Corollary \ref{cor:general_strips} implies that $||\phi||_{H}\geq \frac{\pi}{2}r^{2}$.
}\end{rem}

\begin{proof}
We start by noting that no bubbling is possible under the assumption (\ref{eq:connect2}).
 As in the proof of Corollary \ref{cor:strips_energy}, to prove our claim it is enough  to show that for any \emph{regular} pair  $(H,J)$ where $J\in \mathcal{J}_{\omega}$, $H:[0,1]\times M\to \R$, $\phi^{H}_{1}(L)=L'$, and for
any $x\in L$ there exists $u\in \mathcal{M}_{L,L',J,0}'$ such that $\int u^{\ast}\omega \leq Var (H)$, $\mu(u)\leq n$ and $x\in Im(u)$ (the relevant
moduli spaces as well as the notion of regular pairs are defined in this setting
in the same way as in \S\ref{subsec:recall}).

We fix such a regular pair $(H,J)$ and
we start to adapt the construction of $EF(L,L')$ to our new setting.
\subsubsection{The action functional.}
Consider an additional Hamiltonian $G:[0,1]\times M\to \R$. We first verify
that the action functional $\mathcal{A}_{L,L',G}$ from (\ref{eq:action}) continues to be well defined in our new setting.
Indeed, let $\eta_{0}=z_{0}\in L\cap L'$ such that the path $(\phi^{H}_{t})^{-1}(z_{0})$ is null in $\pi_{1}(M,L)$. We need to show that if $\overline{x}:[0,1]\to \mathcal{P}_{z_{0}}(L,L')$ is a path joining $z_{0}$ to $x\in \mathcal{P}_{z_{0}}(L,L')$, then the expression in (\ref{eq:action}) only depends
on $x$ and not on $\overline{x}$. This obviously comes down to showing that
$\int \overline{x}^{\ast}\omega$ only depends of $x$. For this we use a hamiltonian isotopy $\phi^{H'}$ inverse of $\phi^{H}$ such that $\phi^{H'}_{1}(L')=L$ and the map
$b_{H'}:\mathcal{P}_{z_{0}}(L,L')\to \mathcal{P}_{\ast}(L,L)$ defined by $(b_{H'}(x))(t)=\phi^{H'}_{t}(x(t))$ as in \S\ref{subsubsec:naturality}.
By applying the action functional computation in \S\ref{subsubsec:naturality}
we see that proving that $\int \overline{x}^{\ast}\omega$ only depends
of $x$ is equivalent to showing that for any $\overline{x}':[0,1]\to \mathcal{P}_{\ast}(L,L)$
with $\overline{x}'(0)=z_{0}'=b_{H'}(z_{0})$
the integral $\int (\overline{x}')^{\ast}\omega$ only depends
of $x'=\overline{x}'(1)$. But now as $[z_{0}']=0\in\pi_{1}(M,L)$ we consider
$d:[0,1]\to \mathcal{P}_{\ast}(L,L)$ which contracts $z_{0}'$ to the constant path
and if $\overline{x}''$ is a second path with the same ends as $\overline{x}'$, then
we may consider the concatenation $d\#\overline{x}'\#(\overline{x}'')^{-1}\#d^{-1}$
and we notice that this represents geometrically a map defined on a disk
whose boundary rests on $L$. Therefore, the integral of $\omega$ on this disk vanishes
and this implies that $\int (\overline{x}')^{\ast}\omega = \int (\overline{x}'')^{\ast}\omega$.

\subsubsection{The Maslov index.}
The main other difficulty that remains to be solved is that, when (\ref{eq:connect2}) replaces (\ref{eq:connectivity}), the Maslov index of a strip $u\in \mathcal{M}(x,y)$ as defined in \S\ref{subsubsec:Maslov_ind} does not only depend on the ends of the strip, $x,y$.
The space $\mathcal{P}_{z_{0}}(L,L')$ is not in general simply connected.
We let  $\Pi=\pi_{1}(\mathcal{P}_{z_{0}}(L,L'))$. There is a natural morphism
$\mu:\Pi\to \Z$ which is defined as follows: we fix $x\in I(L,L';G)$ and we consider a $C^{\infty}$
path $\gamma:[0,1]\to \mathcal{P}_{z_{0}}(L,L')$ such that $\gamma(0)=x=\gamma(1)$ and $[\gamma]=g\in \Pi$.
We then put $\mu(g)=\mu(\gamma)$ (where $\mu(\gamma)$ is computed by viewing $\gamma$ as
a ``strip" joining $x$ to $x$ and by using the method described in
\S\ref{subsubsec:Maslov_ind}). It is easy to see that this is well-defined and that it defines a
homomorphism. Let $\pi$ be the image of this homomorphism (obviously $\pi$ is isomorphic to $\Z$ or trivial) and let $Ker(\mu)$ be its kernel.
There exists a regular covering $\tilde{\mathcal{P}}$ of $\mathcal{P}_{z_{0}}(L,L')$ with covering projection
$p:\tilde{\mathcal{P}}\to \mathcal{P}_{z_{0}}(L,L')$, covering group equal to $\pi$ and
such that $\pi_{1}(\mathcal{\tilde{P}})=Ker(\mu)$.

Denote by $\tilde{I}(L,L',G)=p^{-1}(I(L, L';\eta_{0}, G))$ and
let $x,y\in \tilde{I}(L,L',G)$. We define $\mu(x,y)=\mu(u)$ where
$u\in \mathcal{S}(p(x),p(y))$ verifies $u=p(u')$ with
$u':[0,1]\to \tilde{\mathcal{P}}$ a path joining $x$ to $y$.  Notice that with this definition we have $\mu(x,y)=\mu(gx,gy)$ and $\mu(gx,y)=\mu(g)+\mu(x,y)$
for any $g\in\pi$ (we consider that $\pi$ acts on $\tilde{\mathcal{P}}$ on the left).  Fix $x_{0}\in \tilde{I}(L,L',G)$. We also
define an absolute Maslov index for the points  $y\in\tilde{I}(L,L',G)$ by letting
$\mu(y)=\mu(y,x_{0})$ - clearly this depends on the choice of $x_{0}$.

We end this sub-subsection with a remark that will be useful later on.
There exists a natural map $j_{L}: L\to \mathcal{P}_{z_{0}}(L,L')$ which is
defined by $j_{L}(x)=\phi^{H}_{t}(x)$. This map has the property that $l\circ j_{L}=id_{L}$ where $l(\gamma)=\gamma(0)$. Therefore, we may view $\pi_{1}(L)$
as a subgroup of $\Pi$. The remark in question is that $j_{L}(\pi_{1}(L))\subset Ker(\mu)$. To verify this, notice that this property is homotopy invariant and so
it is sufficient to check it in case $L'$ is close to $L$ and is
the image of some $df$ where $f:L\to R$ is some Morse function which
is $C^{2}$-small. In this case the relative Maslov index agrees with the relative Morse index and so does only depend on the ends of strips and not on the strips
themselves.

\subsubsection{The moduli spaces.}
We pursue the construction by defining for $x,y\in\tilde{I}(L,L',G)$ the moduli spaces  $\mathcal{N}_{L,L',J,H}(x,y)=\{u\in C^{\infty}(\R, \widetilde{\mathcal{P}}) : (p\circ u)\in\mathcal{M}_{L,L',J,G}(p(x),p(y)),\\ u(-\infty)=x,\  u(+\infty)=y \}$. These moduli spaces behave in a way perfectly similar
to the behaviour of $\mathcal{M}(x,y)$ (when condition (\ref{eq:connectivity})
is satisfied) as described in   \S\ref{subsubsec:Maslov_ind}, \S\ref{subsubsec:naturality} and \S\ref{subsubsec:gromov_comp}. In particular,
the manifolds with corners structure of their compactifications remains true
(the proof discussed in \S\ref{sec:appendix} applies in this case without modification).
We also have, $0< E_{L,L',G}(p\circ u)=\mathcal{A}(p(x))-\mathcal{A}(p(y))$
when $u\in \mathcal{N}(x,y)$, $\mu(x,y)-1=dim(\mathcal{N}(x,y))$,
and the formula
$$\partial\overline{\mathcal{N}}(x,y)=\bigcup_{z} \overline{\mathcal{N}}(x,z)\times \overline{\mathcal{N}}(z,y)$$
remains valid. In particular, we claim that only a finite number
of non-trivial terms appear in this union.
This is because for any $B>0$ and any $x,y\in I(L,L';z_{0},G)$
there is at most a finite number of homotopy classes $[u]$ of paths in $\mathcal{P}_{z_{0}}(L,L')$ that
join $x$ to $y$ and are represented by
strips $u\in \mathcal{M}(x,y)$ with $E_{L,L',G}(u)\leq B$
(otherwise, by passing to a convergent sub-sequence of strips, each of a different
homotopy class, Gromov compactness would be contradicted). This means that, for any
$x\in \tilde{I}(L,L',G)$, $y\in I(L,L';\eta_{0},G)$ there are at most
a finite number of points
$z\in p^{-1}(y)$ such that $\mathcal{N}(x,z)\not=\emptyset$. As the number
of points in $I(L,L';\eta_{0},G)$ is finite
this means that, for fixed $x\in \tilde{I}(L,L',G)$, there are only finitely
many non-vanishing spaces $\mathcal{N}(x,z)$
and implies the claim.

The spaces $\mathcal{N}(x,y)$ have the additional property that they
are equivariant in the sense that the left action of $\pi$ induces
a homeomorphism $\mathcal{N}(x,y)\stackrel{g}{\longrightarrow} \mathcal{N}(gx,gy)$ for any $g\in \pi$.

\subsubsection{The extended Floer complex and the spectral sequence.}\label{subsubsec:extend_ss_seq}
 The next step is to construct an extended Morse complex by following the method in \S\ref{subsec:constr_SS}. We choose elements
$\tilde{x}\in p^{-1}(x)$ for each $x\in I(L,L';\eta_{0},G)$ so that all the other
elements in $p^{-1}(x)$ can then be uniquely written as $g\tilde{x}$, $g\in \pi$.
We want to construct a representing chain system $s_{xy}\in S_{\mu(x,y)-1}\overline{\mathcal{N}}(x,y)$, $x,y\in \tilde{I}(L,L',G)$ for the
moduli spaces $\mathcal{N}(x,y)$.
We consider a finite, increasing sequence
of strictly positive numbers
$\Delta_{1},\Delta_{2},\ldots,\Delta_{q}$ such that
for any $a,b\in I(L,L;G)$ there exists an $i$ verifying $|\mathcal{A}(a)-\mathcal{A}(b)|=\Delta_{i}$.
Assume by induction that the $s_{xy}$ have been constructed for all $x,y$
such that $\mathcal{A}(p(x))-\mathcal{A}(p(y))\leq \Delta_{k}$.
We then consider a couple $\tilde{x},y$ such that $\mathcal{A}(p(\tilde{x}))-
\mathcal{A}(p(y))=\Delta_{k+1}$ and $\mathcal{N}(\tilde{x},y)\not=\emptyset$ (there are finitely many such couples as mentioned
before). We proceed as in Lemma \ref{lem:repres_chain_syst} to construct
$s_{\tilde{x}y}$ and we then define $s_{(g\tilde{x})(gy)}$ to be
$g(s_{xy})$.

To continue with the construction we
fix a map $s:L\to X$ such that $X$ is simply connected and
$s$ carries the $0$-ends of the paths in $I(L,L';G)$ to a distinguished base-point $\ast$ in $X$ (in the original
construction the role of $X$ was played by $\tilde{L}$). In our applications
$s$ will be a degree one map (and so $X=S^{n}$).
We have an obvious map $\tilde{s}:\tilde{\mathcal{P}}\stackrel{p}{\longrightarrow} \mathcal{P}_{z_{0}}(L,L')\stackrel{l}{\longrightarrow} L\to X$
(where $l(\gamma)=\gamma(0)$).
This map takes the representing chain system $s_{xy}$ to cubical chains $u_{xy}\in S_{\ast}(\Omega'X)$.
The advantage of using $S_{\ast}(\Omega'X)$ with $X$ simply connected is that $\Z[\pi]$ acts
trivially on $\Omega' X$. Moreover, this loop space is connected in contrast to
$\Omega L$ if $\pi_{1}(L)\not=0$. We also notice $u_{xy}=u_{(gx)(gy)}$ for all $g\in \pi$.

Denote by $\mathcal{R}^{X}_{\ast}$ the ring $S_{\ast}(\Omega'X)$ and let $\tilde{\mathcal{R}}= \mathcal{R}^{X}_{\ast}\otimes \Z[\pi]$ be graded as a tensor  product and endowed with the differential coming from the first factor.
For $x,y\in I(L,L',G)$ consider $v_{xy}\in \tilde{\mathcal{R}}$
be defined by $v_{xy}=\sum  u_{\tilde{x}g(\tilde{y})}\otimes g$.
The extended Floer complex in this situation will be denoted
by $\tilde{\mathcal{C}}^{J,\zeta}(L,L';G)$. It  is a free $\tilde{\mathcal{R}}$-
chain complex with generators the elements of $I(L,L',G)$.
Its grading is defined as follows. Recall that we have already fixed an absolute index for
the points in $\tilde{I}(L,L',G)$. As the generators of $\Z[\pi]\otimes \Z/2<I(L,L';G)>$
are in bijection with the elements in $\tilde{I}(L,L',G)$
this absolute index gives a grading to $\Z[\pi]\otimes \Z/2<I(L,L',G)>$ and the grading on $\tilde{\mathcal{C}}$
is the canonical tensor product grading (this is compatible with the action of $\tilde{\mathcal{R}}$ because
$\mu(gx,y)=\mu(g)+\mu(x,y)$).
Its differential is defined by $dx=\sum_{y}v_{xy}\otimes y$ (using $u_{xy}=u_{(gx)(gy)}$ it is easily seen
that $d^{2}=0$). We denote by $|-|$ the grading defined before and
we let $F^{k}=\{a\in \Z[\pi]\otimes \Z/2<I(L,L',G)> : |a|\leq k\}$.

The spectral sequence in this more general setting, $\tilde{E}F(L,L';G)$,
is induced by the filtration $\tilde{F}^{k}\tilde{\mathcal{C}}=
\mathcal{R}^{X}_{\ast}\otimes F^{k}$. This spectral sequence is not anymore a first
quadrant spectral sequence in general but rather an upper semi-plane sequence. It is however
a rather well-behaved spectral sequence because $\pi$ acts by translation parallel to the
$x$-axis on this sequence and, as consequence of the fact that $\tilde{C}$ is a free, finitely generated  $\mathcal{R}^{X}_{\ast}\otimes \Z[\pi]$-module, the sequence is equivariant with respect to this action
(in the sense that $d^{r}(ga)=gd^{r}a$).

\subsubsection{The Morse case.}
It is easy to see that all the other properties of $EF(L,L';G)$ - with the proofs
provided in \S\ref{subsubsec:inv} - extend to the case of $\tilde{E}F(L,L';G)$
without difficulty (the key point of course being that the action functional
continues to be well defined and it is equivariant with respect to the action of $\pi$). In particular,
the pages of order greater or equal than $2$ of the spectral sequence are invariant up to translation and
as the isomorphisms in question are also naturally $\pi$-equivariant we obtain that
$\tilde{E}F^{r}(L,L'), r\geq 2$ is an invariant up to translation which is $\pi$-equivariant.  In particular $\tilde{E}F^{r}(L,L')$ is isomorphic up to
translation for $r\geq 2$ to the analogue spectral sequence arising from a Morse function on $L$. Therefore, the last stage consists in detecting what it the output of this construction in the  Morse function context (similar to \S\ref{subsubsec:Morse_Serre_ss}).  Consider the pull-back covering of $L$ which is induced from $\tilde{\mathcal{P}}
\to \mathcal{P}_{z_{0}}(L,L')$ by the map $j_{L}:L\to \mathcal{P}_{z_{0}}(L,L')$.
 We denote this covering by $p:\hat{L}\to L$.
 As $j_{L}(\pi_{1}(L))\subset Ker(\mu)$
this covering is trivial. Therefore, $\hat{L}$ is homeomorphic to $\pi\times L$.

Fix a Morse function $f:L\to \R$ and let $\hat{f}=f\circ p$. Moreover, consider the fibration $\Omega' X\to \hat{E}\to \hat{L}$ which is the pullback of the fibration  $\zeta:\Omega'X\to P'X\to X$ over the map $s\circ p$ (recall that we have fixed $s:L\to X$ in \S\ref{subsubsec:extend_ss_seq}). This fibration consists simply of $\pi$-copies of the fibration $\Omega'X\to E\to L$ which is the pull-back of $\zeta$ over $s$.  Let $EF_{X}$ be the Serre spectral sequence of this last
fibration. The argument in \S\ref{subsubsec:Morse_Serre_ss} immediately implies that the spectral sequence constructed by following the method above in this Morse case  - as described in \S\ref{subsubsec:Morse_ss} when (\ref{eq:connectivity}) is satisfied -
is isomorphic up to translation to the spectral sequence $\Z[\pi]\otimes EF_{X}$.

Therefore, when $r\geq 2$, we have a $\pi$-equivariant isomorphism up to translation between $\tilde{E}F^{r}(L,L')$ and $\Z[\pi]\otimes EF_{X}^{r}$.

\subsubsection{End of the proof.}
Once the machinery above is constructed we take $X=S^{n}$ and $s:L\to S^{n}$ to be
a degree one map and the proof in Corollaries \ref{cor:strips_energy}
and \ref{cor:non-squeeze} i. proceeds without change.

\end{proof}

\appendix

\section[Structure of the Floer Moduli spaces]
{Structure of manifolds with corners on Floer moduli spaces.}
\label{sec:appendix}

\subsection{Introduction.}

Let $(M,\omega)$ be a symplectic manifold which is convex at infinity and
let  $L$ and $L'$ be two simply connected compact Lagrangian sub-manifolds
in $M$. The symplectic form $\omega$ will be supposed to vanish on
$\pi_{2}(M)$, so that there are no symplectic spheres in $M$,
nor symplectic disks attached to $L$ or to $L'$. Similarly, we assume that
the first Chern class $c_{1}(M)$ vanishes on $\pi_{2}(M)$
so that the Maslov index of Floer trajectories only depends on the ends of the trajectories.
Suppose moreover that $L$ and $L'$ intersect transversally.

The purpose of this section is to endow the moduli spaces, $\overline{\M}(x,y)$,
of Floer trajectories - pseudo-holomorphic
strips, in our case -  with the topological structure of ``manifolds with corners''
(see Definition \ref{def:corner manifold}).
In \cite{Fl2}, A.~Floer introduced the gluing construction to
treat the case of relative index $1$. His work extends almost verbatim to
the case of higher relative indexes, but some particular care is needed when the number
of breaking points is bigger than one. In this case, Floer's argument - as he described it
in \cite{Fl4} (Proposition 2d.1.) - only provides stratum by stratum
homeomorphisms, i.e. local maps of the form
$$
\M(x_{0},x_{i_{1}})\times\dots\times\M(x_{i_{r}},x_{k})\times(0,1)^{r}
\xrightarrow{\varphi_{i_{1},\dots,i_{r}}}\M(x_{0},x_{k}).
$$
instead of a map defined up to the boundary, i.e. a local map of the form
$$
\M(x_{0},x_{1})\times\dots\times\M(x_{k-1},x_{k})\times[0,1)^{k-1}
\xrightarrow{\varphi}\overline\M(x_{0},x_{k})
$$
where the ``$k$-fold broken" trajectories are identified with elements of the form
$(u_{1},\ldots,u_{k})\times \{0\}$ and the map $\varphi$ provides
a ``cornered neighborhood" of these trajectories in the sense
that $\varphi$ preserves the natural stratifications on the two sides.
To build this last map out of the former ones, some gluing compatibility
conditions have to be fulfilled. Verifying this conditions
is not obvious, in essence, because the gluing construction relies on an application
of an  implicit function theorem. The question of defining some structure on
moduli spaces of pseudo-holomorphic curves, at least such as to produce
a (virtual) fundamental class, appears as a key point in most applications
of pseudo-holomorphic curves to symplectic geometry, and, in particular, in the
definition of Gromov Witten invariants.

This point has been treated
(both in the context of Floer homology and Gromov Witten
invariants) by different authors (\cite{Liu_Tian},\cite{Ruan},\cite{Fuk},\cite{Li_Tian},\cite{Siebert})
in a very general setting (allowing bubbles). However, this goes far beyond what is required
for the present paper, and we have not been able to find in the
literature a simple and explicit proof of the ``manifold with corners" structure for the moduli
spaces of pseudo-holomorphic strips. For this reason as well as to
make the paper self-contained, we include one here.

We make use of now classical ideas and techniques introduced by different authors
 (\cite{Fl2},\cite{McSal},\cite{Sch}, \cite{Siko1}, \cite{Liu_Tian})
 and we shall follow rather
closely the work of J.-C.~Sikorav \cite{Siko1} about the gluing
construction for compact Riemann surfaces, and adapt it to our strips.
The more recent and much more general technique introduced in \cite{Ho-Ze-Wi}
offers a more conceptual approach to gluing problems of this type.

\medskip

Recall that the action functional $\A$ is defined as follows~: choose a
path $\gamma_{0}$ from $L$ to $L'$, let $\PathSpace_{\gamma_{0}}(L,L')$
be the component of $\gamma_{0}$ in $\PathSpace(L,L')$ and for
$\gamma\in\PathSpace_{\gamma_{0}}(L,L')$, set
$$
\A(\gamma)=-\int_{[0,1]^{2}}\bar\gamma^{*}\omega
$$
where $\bar\gamma$ is a path from $\gamma_{0}$ to $\gamma$ in
$\PathSpace(L,L')$.

Let $J$ be an almost complex structure that tames $\omega$. For
$x_{+},x_{-}\in L\cap L'$, a parametrized Floer trajectory from $x_{-}$
to $x_{+}$ is a map $u:\Sigma\to M$, where
$$
\Sigma=\{z=s+it\in\C, 0\leq\im(z)\leq 1\}
$$
such that
\begin{gather}
  \dbar_{J}u:=du+J(u)du{\i}=0             \label{eq:J-holo}\\
  u(\R\times\{0\})\subset L
\quad\text{ and }\quad
  u(\R\times\{1\})\subset L'\label{eq:BC}\\
  \lim_{s\to\pm\infty}u(s,t)=x_{\pm}(t) \label{eq:ends of trajectories}.
\end{gather}

To shorten notation, the product of a real interval with $[0,1]$ will
be denoted by double brackets~:
$\Sigma=\strip]-\infty,+\infty\strip[=(-\infty,+\infty)\times[0,1]\subset\C$.

\

A non parametrized  Floer trajectory is the orbit of a parametrized one
under the action of $\R$ by translation in the $s$ direction.

Let $\Mp_{J}(x_{-},x_{+})$ be the space of all parametrized Floer
trajectories from $x_{-}$ to $x_{+}$~:
\begin{align}
\Mp(x_{-},x_{+})&=\left\{u:\Sigma\to M,
  \eqref{eq:J-holo}
  \eqref{eq:BC}\text{ and }
  \eqref{eq:ends of trajectories}
\right\}\\
\M(x_{-},x_{+})&=\Mp(x_{-},x_{+})/\R\\
\overline\M_{J}(x_{-},x_{+})&=\text{ Gromov compactification of }\M(x_{-},x_{+})
\end{align}

\begin{rem} {\rm
  For technical reasons, it is sometimes useful to work with ``time dependent''
  almost complex structures $(J_{s})_{s\in\R}$. We will use families that are
  constant near $\pm\infty$. The equation \eqref{eq:J-holo} is then
  naturally replaced by
  \begin{equation}
    \label{eq:J(s)-holo}
    du(s,t)+J(s,u(s,t))\ du(s,t)\ {\i}=0
  \end{equation}

  Finally, hamiltonian perturbations of $L$ and $L'$ are also needed, and
  to take advantage of particular Hamiltonians, we need to
  keep track of them. The action functional is then replaced by
  \eqref{eq:action}, and the $J$-holomorphy equation is replaced by the
  non-homogeneous one~:
  $$
    \frac{\partial u}{\partial s}+J(u)\frac{\partial u}{\partial t}
    =-\nabla H(t,u(s,t))
  $$
  We will first deal with the homogeneous case, and add comments about
  the non homogeneous one in the last section.
}\end{rem}

Let us fix one more notation~: the linearization of the operator
$\dbar_{J}:u\mapsto du+Jdu\i$ at a map $u$ such that $\dbar_{J}u=0$ does
not depend on the connection used to compute it~: this defines an
operator $D_{u}:\L^{1,p}(u^{*}TM)\to \L^{p}(\Lambda^{01}u^{*}TV)$. The
$J$-holomorphic map $u$ is said to be \emph{regular}, if $D_{u}$ is
onto.

\begin{rem}{\rm   The topological assumptions $\omega(\pi_{2}(M))=0$ and
  $c_{1}(\pi_{2}(M))$ are only used to obtain a global structure on the
  moduli spaces, but the underlying ``gluing'' construction is purely
  local and only relies on the regularity of the trajectories under
  consideration.}
\end{rem}

\bigskip

We start with a simple definition~:
\begin{defi}
  \label{def:corner manifold}
  A topological space $X$ is said to have a structure of manifold with
  corners if there is a partition $X=\sqcup_{i\in I} X_{i}$ of $X$ into
  manifolds $X_{i}$ of dimension $i$ such that for each point $x\in
  X_{i}$ there is a neighbourhood $U_{x}$ of $x$ in $X$, and a local
  homeomorphism $\phi:U_{x}\to V_{0}\subset\R^{i}\times [0,1)^{n-i}$
  whose restriction to each $U_{x}\cap X_{j}$ is a
local
  diffeomorphism to the
  $j$-dimensional stratum of $\R^{i}\times [0,1)^{n-i}$
(which is defined
to be the disjoint union of the products of the form
$\R^{i}\times I_{1}\times I_{2}\times\ldots \times I_{n-i}$
where $I_{k}$ is either $\{0\}$ or $(0,1)$ and there are precisely $j-i$
non-zero terms).
\end{defi}

We are interested in such a structure to obtain a homological
representation of the stratification of $\overline\M(x,y)$, as described
in Lemma \ref{lem:repres_chain_syst}.

\begin{theo}
\label{thm:corner structure on Mbar}
For a generic choice of $J$, and two intersection points $x$ and $y$ with
$x\neq y$, the space $\overline\M_{J}(x,y)$ admits a
structure of manifold with corners of dimension $\mu(x,y)-1$, whose $k$
 co-dimensional stratum is the space of trajectories broken at $k$
intermediate points.
\end{theo}

This statement requires a topology on $\overline{\M}(x,y)$ which we
recall now.

For convenience, if $(x_{0},\dots,x_{k})$ is a sequence of intersection
points of decreasing indexes, let
\begin{equation}
 \M(x_{0},\dots,x_{k})=\M(x_{0},x_{1})\times\dots\times\M(x_{k-1},x_{k})
 \label{eq:M(x0,...,xk)=M(x0,x1) x...x M(xk-1,xk)}~.~
  \end{equation}


Consider a curve $C_{\infty}\in\M(x_{0},\dots,x_{k})$, and a
sequence $C_{n}$ of curves in $\M(x_{0},x_{k})$.

To express the convergence of $C_{n}$ to $C_{\infty}$ in the sense of the Gromov
topology, pick some parametrizations
$u_{\infty}:\Sigma_{\infty}=\sqcup_{i=1}^{k}\Sigma_{\infty,i}\to M$ of
$C_{\infty}$ and $u_{n}:\Sigma\to M$ of $C_{n}$. Then, for each $i\in\{1,\dots,k\}$, there is a unique
$s_{n,i}\in\Sigma$ such that
$$
 \A(u_{n}(s_{n,i},\cdot))=\A(u_{\infty,i}(0,\cdot)).
$$

The sequence $C_{n}$ is said to converge to $C_{\infty}$ in the Gromov topology
if, for all $i\in\{1,\dots,k\}$, $u_{n}(\cdot-s_{n,i},\cdot)$
$\RegClass^{0}$
converges to $u_{\infty,i}$ on all compact subsets of $\Sigma$.

This definition naturally extends to broken trajectories~: a sequence of broken trajectories converges,
if the topology of the domain stabilizes and each smooth component
converges in the previous sense.

\subsection{Sketch of the construction}

Let $C_{\infty}\in\M(x_{0}\dots x_{k})$, $C_{\infty}=(C_{\infty,1},\dots
,C_{\infty,k})$ be a regular curve, by which we mean that the
linearization of the Cauchy Riemann equation \eqref{eq:J-holo} on each
component is onto. Theorem \ref{thm:corner structure on Mbar}, will
be proved by constructing a local chart centered at $C_{\infty}$~:
\begin{equation}
  \label{eq:local charts on Mbar}
  \phi:
  \begin{array}{ccc}
    \overline\M(x_{0},x_{k})
    &\to&
    (1,+\infty]^{k-1}\times T_{C_{\infty}}\M(x_{0}\dots x_{k})\\
    C&\mapsto&(\rho(C),\pi_{\rho}(C))
  \end{array}
\end{equation}
satisfying the conditions of definition \ref{def:corner manifold}.

The first component of this map will be called the gluing
parameter and is a small perturbation of the parameter $\rho$ defined as follows~: we choose one regular level
$a_{i}$ ($a_{i}=\frac{\A(x_{i-1})+\A(x_{i})}{2}$ for instance) of the
action functional between each pair of critical values $\A(x_{i})$. Then we
measure the time a trajectory needs to run from one level to the next one
by setting
\begin{equation}
  \label{eq:gluing parameter}
  \rho:
  \begin{array}{ccc}
    \overline\M(x_{0},x_{k})&\to&(1,+\infty]^{k-1}\\
    C&\mapsto&(2(s_{2}-s_{1}),\dots,2(s_{k}-s_{k-1}))
  \end{array}
\end{equation}
(the only purpose of this $2$ is to simplify future notation) where $s_{i}$ is
defined by $\A(u(s_{i},\cdot))=a_{i}$ for some parametrisation
$u:\bigsqcup_{\alpha}\Sigma_{\alpha}\to M$ of $C$, with the convention
that $s_{i}-s_{i-1}=+\infty$ if $s_{i}$ and $s_{i-1}$ do not belong to
the same component. Of course, each $s_{i}$ depends on the choice of $u$,
but the differences $s_{i}-s_{i-1}$ do not.

\begin{rem}
{\rm  Because of their geometric meaning, our gluing parameters will tend to
  $+\infty$ when the curve splits, and hence $0$ is replaced by
  $+\infty$ and $[0,1)$ by $(1,+\infty]$ in definition \ref{def:corner
    manifold}. As a consequence, a gluing parameter
  $\rho=(\rho_{1},\dots,\rho_{k-1})$ will be considered as ``large'' if
  \emph{all} its component are large, or equivalently, if $e^{-\rho}$ is
  small. In view of this, we define $|\rho|=\inf{\rho_{i}}$.
}\end{rem}

The second component of the local chart is less explicit. Let $\M_{\rho}$
denote the fiber of $\rho$. Then, it is easy to see that, for each $\rho$ large enough,
$\M_{\rho}$ is locally diffeomorphic to $T_{C_{\infty}}\M(x_{0},\dots,x_{k})$.
The main difficulty is to control the dependence of these diffeomorphisms
with respect to $\rho$. In particular, they need to be constructed in such a way that
the map $\phi(C)=(\rho(C),\pi_{\rho}(C))$ is a homeomorphism on its image.

To achieve this we proceed as follows~:

\begin{enumerate}
\item
We intend to define an inverse to $\phi$. For this we start by defining
a pre-gluing map $\rho\mapsto w_{\rho}$~: use a gluing parameter
$\rho$ and cutoff functions to glue the different components of $C_{\infty}$ into
a map $w_{\rho}$, that is ``approximately'' a Floer trajectory matching
the transit times $\rho$. Check that any ``exact'' Floer trajectory $C$,
close enough to $C_{\infty}$ and matching the transit times $\rho$, is in
fact a $\L^{1,p}$-small perturbation of $w_{\rho}$.

\item
We then setup the $J$-holomorphy equation for perturbations $\xi\in\Gamma^{1,p}(w_{\rho}^{*}TM)$
of $w_{\rho}$ as a non linear PDE, $\dbartld_{w_{\rho}}\xi=0$. We check that
for $\xi=0$, $\dbartld_{w_{\rho}}(0)$ vanishes as $|\rho|$ tends to
$+\infty$, and that the linearization $D_{\rho}$ of this operator at
$\xi=0$ has a uniformly bounded right inverse. An implicit function
theorem ``near infinity'' ensures then that, for $\rho$ large enough, the equation
$\dbartld_{w_{\rho}}\xi=0$ is regular and provides local diffeomorphisms
$\zeta_{\rho}$ between $\M_{\rho}$ and $\ker D_{\rho}$.

\item
One difficulty at this point is that the topology of the base of $w_{\rho}^{*}TM$
strongly depends on $\rho$, making the spaces $\ker D_{\rho}$ difficult
to compare. However, these bundles are almost the same ``away from the
nodes''~: restricted to compact sets that avoid the nodes, and pushed by
parallel transport, their sections can be considered sections of one
constant bundle. An $\L^{2}$ projection on $\ker
D_{\infty}=T_{C_{\infty}}\M(x_{0},\dots,x_{k})$ induces then an isomorphism
$\pi:\ker D_{\rho}\to\ker D_{\infty}=T_{C_{\infty}}\M(x_{0},\dots,x_{k})$.

\item
Define $\pi_{\rho}=\pi\circ\zeta_{\rho}$. Check that the map
$(\rho,\pi_{\rho})$ and its inverse defined by gluing verify the
conditions in Definition \ref{def:corner manifold}.

\end{enumerate}

 \begin{rem}{\rm We identify below perturbations of $w_{\rho}$  to sections of
  $w_{\rho}^{*}(TM)$ via the exponential diffeomorphism associated to a fixed
  metric $g_{L,L'}$ for which $L$ and $L'$ are totally geodesic.}\end{rem}

\subsection{Pre-gluing}
We first introduce some notation for the gluing operation at the source
level, which are all summarized in figure \ref{fig:Gluing strips}. Then
the ``pre''-gluing operation will be described.

\subsubsection{Gluing strips} Let
$\Sigma_{\infty}=\Sigma_{\infty,1}\sqcup\dots\sqcup\Sigma_{\infty,k}$
denote the disjoint union of $k$ copies of the standard strip, $\Sigma=\strip]-\infty,+\infty\strip[$.

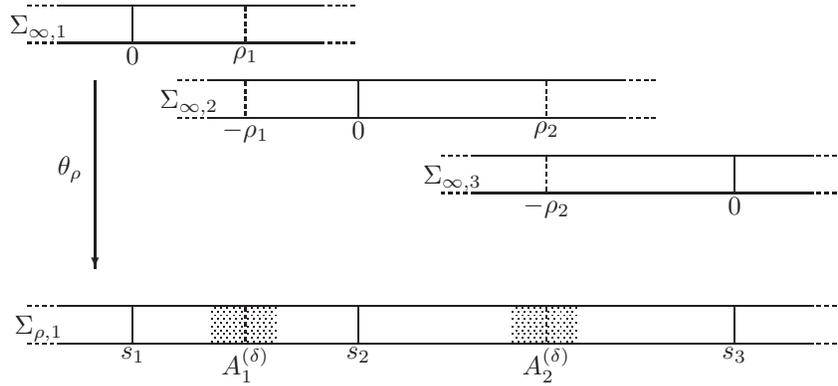
\begin{figure}[htbp]
  \begin{center}
\setlength{\unitlength}{1mm}%
\begin{picture}(110,52)(-6,-5)
  \def\drawstrip(#1,#2)#3#4{%
    \put(#1,#2){%
      \multiput(0,0)(0,5){2}{
        \multiput(0,0)(-1,0){5}{\line(1,0){.5}}
        \put(0,0){\line(1,0){#3}}
        \multiput(#3,0)(1,0){5}{\line(1,0){.5}}
        }
      \put(-2.5,2){\makebox(0,0)[c]{#4}}
      }
    }
  \def\drawstripbar(#1,#2)#3{%
    \put(#1,#2){\line(0,1){5}\put(0,-.5){\makebox(0,0)[t]{#3}}}
    }
  \def\drawstripbardot(#1,#2)#3{
    \put(#1,#2){\multiput(0,0.3)(0,1){5}{\line(0,1){.5}}
                \put(0,-.5){\makebox(0,0)[t]{#3}}
                }
    }
  \def\drawstripzone(#1,#2)#3{
    \linethickness{.1mm}
    \put(#1,#2){
      \multiput(-0.5,0.25)(0,1){5}{\multiput(0,0)(0.5,0.5){2}{\bezier{10}(-4,0)(0,0)(4,0)}}
      \put(0,-.5){\makebox(0,0)[t]{#3}}
      }
    }
\drawstrip( 0, 0){100}{$\Sigma_{\rho,1}$}
\drawstrip( 0,40){35}{$\Sigma_{\infty,1}$}
\drawstrip(20,30){55}{$\Sigma_{\infty,2}$}
\drawstrip(55,20){45}{$\Sigma_{\infty,3}$}
\drawstripbar(10,0){$s_{1}$}
\drawstripbar(40,0){$s_{2}$}
\drawstripbar(90,0){$s_{3}$}
\drawstripbar(10,40){$0$}
\drawstripbar(40,30){$0$}
\drawstripbar(90,20){$0$}
\drawstripbardot(25,40){$ \rho_{1}$}
\drawstripbardot(25,30){$-\rho_{1}$}
\drawstripbardot(25, 0){}
\drawstripbardot(65,30){$ \rho_{2}$}
\drawstripbardot(65,20){$-\rho_{2}$}
\drawstripbardot(65, 0){}
\put(5,35){\vector(0,-1){25}}
\put(0,22.5){$\theta_{\rho}$}
\drawstripzone(25,0){$A_{1}^{(\delta)}$}
\drawstripzone(65,0){$A_{2}^{(\delta)}$}
\end{picture}
    \caption{The strip $\Sigma_{\rho}$ constructed from pieces of
      $\Sigma_{\infty,i}$.}
    \label{fig:Gluing strips}
  \end{center}
\end{figure}

Consider now a gluing parameter
$\rho=(\rho_{1},\dots,\rho_{k-1})\in(1,+\infty]^{k-1}$. We define
$\Sigma_{\rho}=\Sigma_{\rho,1}\sqcup\dots\sqcup\Sigma_{\rho,\alpha}$ as
the disjoint union of $\alpha=1+\#\{i, 1\leq i\leq k-1,
\rho_{i}=+\infty\}$ copies of the standard strip which are obtained
by gluing together pieces of $\Sigma_{\infty}$ as described in \ref{fig:Gluing strips}.
Explicitly, if $\rho_{i}\not=\infty$ we paste $\Sigma_{\infty,i}\backslash (\strip]\rho_{i},+\infty\strip[)$ and $\Sigma_{\infty,i+1}\backslash (\strip]-\infty,-\rho_{i}\strip[)$
by identifying $\{\rho_{i}\}\times [0,1]\subset \Sigma_{\infty,i}$ with
$\{-\rho_{i}\}\times [0,1]\subset \Sigma_{\infty,i+1}$. The pasting function
is the obvious translation on the first coordinate in $\strip]-\infty,+\infty\strip[$
and the identity on the second. In case $\rho_{i}=\infty$ no gluing occurs between $\Sigma_{\infty,i}$ and $\Sigma_{\infty,i+1}$. In particular, if both $\rho_{i}$ and $\rho_{i+1}$ are infinite, then $\Sigma_{\infty,i+1}$ represents itself a component of $\Sigma_{\rho}$.
We let $\theta_{\rho,i}:\strip]-\rho_{i-1},\rho_{i}\strip[\subset \Sigma_{\infty,i}\to \Sigma_{\rho}$ be the obvious inclusion and we put $s_{i}=\theta_{\rho,i}(0)$.
Moreover,  for $\delta>0$,  the neighbourhood of size $\delta$ of the ``gluing''
region will be denoted by $A_{i}^{(\delta)}$~:
\begin{equation}
  \label{eq:A_i}
  A_{i}^{(\delta)} =
  \strip]s_{i}+\rho_{i}-\delta,s_{i+1}-\rho_{i}+\delta\strip[
\end{equation}
(in other words, if $\rho_{i}=+\infty$, then
$A_{i}^{(\delta)}=\emptyset$, otherwise $A_{i}^{\delta} =
\strip]s_{i}+\rho_{i}-\delta,s_{i}+\rho_{i}+\delta\strip[$).

\subsubsection{Pre-gluing maps}
\begin{figure}[htbp]
  \begin{center}
\setlength{\unitlength}{1mm}%
\begin{picture}(110,30)(-6,-5)
\def\drawstrip(#1,#2)#3#4{%
  \put(#1,#2){%
    \multiput(0,0)(0,5){2}{
      \multiput(0,0)(-1,0){5}{\line(1,0){.5}}
      \put(0,0){\line(1,0){#3}}
      \multiput(#3,0)(1,0){5}{\line(1,0){.5}}
      }
    \put(-2.5,2){\makebox(0,0)[c]{#4}}
    }
  }
\def\drawstripbar(#1,#2)#3{%
  \put(#1,#2){\line(0,1){5}\put(0,-.5){\makebox(0,0)[t]{#3}}}
  }
\def\drawstripbardot(#1,#2)#3{
  \put(#1,#2){\multiput(0,0.3)(0,1){5}{\line(0,1){.25}}
    \put(0,-.5){\makebox(0,0)[t]{#3}}
    }
  }
\def\cutoff(#1,#2){
  \put(#1,#2){%
    \bezier{0}(0,10)(1,10)(1,5)
    \bezier{0}(1,5)(1,0)(2,0)
    }
  }
\def\cuton(#1,#2){
  \put(#1,#2){%
    \bezier{0}(0,0)(1,0)(1,5)
    \bezier{0}(1,5)(1,10)(2,10)
    }
  }
\drawstrip( 0, 0){100}{$\Sigma_{\rho}$}
\drawstripbar(10,0){$s_{i-1}$}
\drawstripbar(40,0){$s_{i}$}
\drawstripbar(90,0){$s_{i+1}$}
\drawstripbardot(17.5, 0){}\drawstripbardot(17.5, 5){}
\drawstripbardot(25.0, 0){}
\drawstripbardot(32.5, 0){}\drawstripbardot(32.5, 5){}
\drawstripbardot(52.5, 0){}\drawstripbardot(52.5, 5){}
\drawstripbardot(65.0, 0){}
\drawstripbardot(77.5, 0){}\drawstripbardot(77.5, 5){}
\put(0,6){
\put(-2.5,0){\makebox(0,10)[c]{$w_{\rho}$}}
\put(0,10){\line(1,0){15.5}}
\put(7.75,11.5){\makebox(0,0)[c]{$u_{\infty,i-1}$}}
\cutoff(15.5,0)
\put(17.5,0){\line(1,0){15}}
\put(25,1.5){\makebox(0,0)[c]{$x_{i}$}}
\cuton(32.5,0)
\put(34.5,10){\line(1,0){16}}
\put(42.5,11.5){\makebox(0,0)[c]{$u_{\infty,i}$}}
\cutoff(50.5,0)
\put(52.5,0){\line(1,0){25}}
\put(65,1.5){\makebox(0,0)[c]{$x_{i+1}$}}
\cuton(77.5,0)
\put(79.5,10){\line(1,0){20}}
\put(89.5,11.5){\makebox(0,0)[c]{$u_{\infty,i+1}$}}
}
\multiput(46.25,0)(12.5,0){4}{
  \put(-1.25,2){\vector(-1,0){5}}
  \put( 1.25,2){\vector( 1,0){5}}
  \put(0,0){\makebox(0,5)[c]{$\scriptscriptstyle\frac{\rho_{i}}{2}$}}
  }
\multiput(50.5,10)(0,5){2}{\drawstripbardot(0,0){}}
\multiput(52.5, 0)(0,5){4}{\drawstripbardot(0,0){}}
\put(51.5,18){
  \put(-1,0){\vector( 1,0){0}}
  \put( 1,0){\vector(-1,0){0}}
  \put(-2,0){\line( 1,0){4}}
  \put(0,0.5){\makebox(0,0)[cb]{$\scriptscriptstyle1$}}
  }
\end{picture}
    \caption{The pre-gluing construction.}
    \label{fig:w_rho}
  \end{center}
\end{figure}
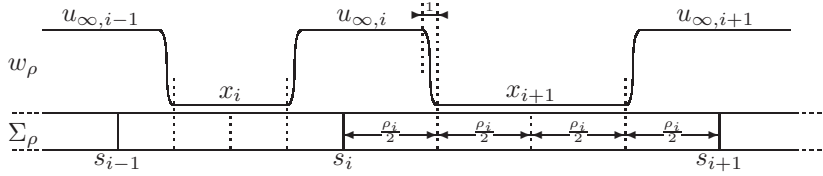

Let $u_{\infty}:\Sigma_{\infty}\to M$ be the parametrization of
our regular curve $C_{\infty}$ which reaches level $a_{i}$ at time $0$
on each component. We now
use the usual technique of cut-off functions to define an ``almost''
trajectory $w_{\rho}:\Sigma_{\rho}\to M$ which is ``close'' to
$C_{\infty}$, and agrees with the transit times $\rho$. The construction
is summarized in figure \ref{fig:w_rho}.

Let $\cutoff{}{}:\R\to\R$ be a smooth function such that
$\cutoff{}{}(t)=1$ if $t\leq0$ and $\cutoff{}{}(t)=0$ if $t\geq 1$~;
given two distinct real numbers $a,b$, we define the function
$\cutoff{a}{b}(t)=\cutoff{}{}(\frac{t-a}{b-a})$.

For $|s|$ large enough, $u_{\infty,i}(s,t)$ can be written in the form
$$
\begin{cases}
 u_{\infty,i}(s,t)=\exp_{x_{i-1}}(\xi_{\infty,i}^{-}(s,t))& (s<<0)\\
 u_{\infty,i}(s,t)=\exp_{x_{i  }}(\xi_{\infty,i}^{+}(s,t))& (s>>0)
\end{cases}
$$
where $\xi_{\infty,i}^{+}\in T_{x_{i}}M$, and $\xi_{\infty,i}^{-}\in
T_{x_{i-1}}M$.

Then, for $\rho$ large enough, we define $w_{\rho}:\Sigma_{\rho}\to M$ by
setting, if $z=\theta_{i}(s,t)$~:
\begin{align}
  w_{\rho}(z)=
    \begin{cases}
      \exp_{x_{i-1}}(\cutoff{-(\rho_{i-1})/{2}+1}{-(\rho_{i-1})/{2}}\ \xi^{-}_{\infty,i  }(s,t))
      &\text{ if } s\leq -\frac{\rho_{i-1}}{2}+1\\
      u_{\infty,i}(s,t)
      &\text{ if }-\frac{\rho_{i-1}}{2}+1\leq s\leq \frac{\rho_{i}}{2}-1\\
      \exp_{x_{i  }}(\cutoff{\rho_{i}/{2}-1}{\rho_{i}/{2}}\ \xi^{+}_{\infty,i  }(s,t))
      &\text{ if } s\geq \frac{\rho_{i}}{2}
    \end{cases}
\end{align}
Notice that, with this definition, $w_{\rho}$ is constant around
the gluing region, namely for $s_{i}+\rho_{i}/2\leq s\leq
s_{i+1}-\rho_{i}/2$.

This map $w_{\rho}$ is as smooth as $u_{\infty}$. It is $J$-holomorphic
in the exterior of $\cup A_{i}^{(\rho_{i}/2+1)}$ and, because of the exponential
decay \cite{RobSal} of $u_{\infty}$ near each $x_{i}$, there
are non negative constants $A$ and $\lambda$ such that
\begin{equation}\label{eq:dbar w_rho}
\Vert\dbar_{J} w_{\rho}(z)\Vert_{1,p}\leq Ae^{-\lambda|\rho|}
\end{equation}

\subsection{Holomorphic perturbations of $w_{\rho}$}\label{sec:analysis in Cn}

\subsubsection{Some analytic properties of the standard Cauchy Riemann equation}
\label{subsec:Asymptotic Op}
All the analytic results needed in the sequel
concern the standard Cauchy Riemann equation applied to $\Sigma$ in $\C^{n}$
with boundary conditions imposed by two transversal Lagrangian linear subspaces $\Lambda,
\Lambda'\subset \C^{n}$. We gather these results here.

We first recall the following lemma from \cite{Sch3}(Theorem 3.1.13). It
is relatively easy for $p=2$ but much more delicate for $p>2$~:
\begin{lem}
\label{lem:operator P}
  For all $p\geq 2$, there is a bounded operator $P$~:
  \begin{equation}
    \label{eq:P}
    \L^{p}(\Sigma,\C^{n})
    \xrightarrow{P}
    \L^{1,p}((\Sigma,\partial\Sigma),(\C^{n},\Lambda,\Lambda'))
  \end{equation}
  such that $\dbar\circ P=\Id$. In particular, there exists a constant $c(p)$
  such that for all $f\in
  \L^{1,p}((\Sigma,\partial\Sigma),(\C^{n},\Lambda,\Lambda'))$~:
  \begin{equation}
    \label{eq:||f||1p <c(p) ||dbar f||p}
    \Vert f\Vert_{1,p}\leq c(p)\Vert \dbar f\Vert_{p}
  \end{equation}
\end{lem}
The proof of this lemma in \cite{Sch3} is given for tubes $\R\times
S^{1}$ instead of strips, with appropriate deformation of the
Cauchy-Riemann equation on the ends, but the boundary case is a strict
analogue, where the invertibility assumption for the asymptotic operators
on the ends is replaced by the transversality of $\Lambda$ and $\Lambda'$.

\bigskip

We will also have to estimate the ``growth''  of holomorphic strips. The
following lemma is a corollary of \cite{RobSal} in the particular case of
an integrable almost complex structure. However we recall the proof to
fix any ambiguity about constant dependency.

\begin{lem}
\label{lem:std holo half strips}
  Let $\Lambda$ and $\Lambda'$ be two transversal Lagrangian linear
  subspaces
  in $\C^{n}$, and $f:[0,+\infty)\times[0,1]\to\C^{n}$ an $\L^{1,p}$
  holomorphic strip such that $f([0,+\infty)\subset\Lambda$,
  $f([0,+\infty)+i)\subset\Lambda'$. Then there are constants $C$ and
  $\delta>0$, depending only on the relative positions of $\Lambda$ and
  $\Lambda'$ such that~:
  $$
    \forall (s,t)\in [0,+\infty) \times[0,1],\
    \Vert f(s+it)\Vert\leq C \Vert f\Vert_{1,p}e^{-\delta s}
  $$
\end{lem}

\begin{proof}
 First, remark that $A={\i}\frac{\partial }{\partial
   t}:\L^{1,2}([0,1],\C^{n},\Lambda,\Lambda')\to \L^{2}([0,1],\C^{n})$
 is a self adjoint operator. Let $\alpha(s)=\int\Vert
 f(s,t)\Vert^{2}dt$.  Since both $f$ and $\frac{\partial f}{\partial
   s}$ belong to $\L^{1,2}([0,1],\C^{n},\Lambda,\Lambda')$, we have
 $\dot{\alpha}(s)=-2\langle f,Af\rangle$ and
 $\ddot{\alpha}(s)=4\int_{0}^{1}\Vert \frac{\partial f}{\partial
   t}\Vert^{2}dy$.  Moreover, since $\Lambda$ and $\Lambda'$ intersect
 transversally, $A$ is bijective with bounded inverse (namely
 $A^{-1}(g)= \int_{0}^{t}-\i g(x)dx-\pi_{\Lambda}(\int_{0}^{1}-\i
 g(x)dx)$ where $\pi_{\Lambda}$ is the projection on $\Lambda$ in the
 direction of $\Lambda'$). Hence there is a constant $\delta$ such
 that $\int_{0}^{1}\Vert \frac{\partial g}{\partial y}\Vert^{2}dy\geq
 \delta^{2}\int_{0}^{1}\Vert {g}\Vert^{2}dy$ for all functions
 $g\in\L^{1,2}([0,1],\C^{n},\Lambda,\Lambda')$.

 From $\ddot \alpha\geq4\delta^{2}\alpha$, we derive
 $\dot{\alpha}+2\delta \alpha\leq0$~: otherwise, if $\beta(s_{0})=\dot
 \alpha(s_{0})+2\delta \alpha(s_{0})>0$, then, as
 $\dot\beta\geq2\delta\beta$, we have $\beta(s)>0$ $\forall s\geq
 s_{0}$, and $\beta(s)\geq\beta(s_{0})e^{2\delta(s-s_{0})}$. Then
 $\alpha(s)\geq K e^{2\delta s}+B$ for some $K>0$ which is impossible.

 Therefore, $e^{2\delta s }\alpha$ is decreasing, and for $s\geq 1$, we
 have
 $\alpha(s)\leq \alpha(1)e^{-2\delta (s-1)}$, \textit{i.e.}
$ \Vert  f(s,\cdot)\Vert_{2}\leq\Vert f(1,\cdot)\Vert_{2}e^{-\delta(s-1)}$.
 The same argument applied to $\frac{\partial f}{\partial s}$ (which is also holomorphic
 and verifies the needed boundary conditions) leads to
 the estimate
 $$
 \Vert f(s,\cdot)\Vert_{1,2}\leq \Vert
 f(1,\cdot)\Vert_{1,2}\ e^{-\delta(s-1)}
 $$

 Now, using Sobolev embedding, we have $\Vert f(s,t)\Vert\leq K_{1} \Vert
 f(s,\cdot) \Vert_{1,2}$, and on the other hand, since $f$ is
 holomorphic in $\strip[0,2\strip]$, Schwarz' lemma implies $\Vert
 f(1,\cdot)\Vert_{1,2}\leq K_{2}\Vert f\Vert_{\infty,\strip[0,2\strip]}\leq
 K_{2}\Vert f\Vert_{1,p}$ with a uniform constant $K_{2}$. Finally, there is a
 uniform constant $C$ such that,
for $s\geq1$ :
 $$
 \Vert f(s,t)\Vert\leq C \Vert f\Vert_{1,p}\ e^{-\delta s}.
 $$
 The existence of such constants for $0\leq s\leq 1$ is obvious so this ends the proof of the lemma.
\end{proof}

We will also need a bounded version of this lemma, which is a sort of
maximum principle~:

\begin{lem}
\label{lem:std holo bnded strips}
  Let $\Lambda$ and $\Lambda'$ be two transversal Lagrangian linear
  subspaces in $\C^{n}$, and $f:\strip[a,b\strip]\to\C^{n}$ a
  holomorphic strip such that $f([a,b])\subset\Lambda$,
  $f([a,b]+i)\subset\Lambda'$, $b-a>2$. Then there are constants $C$ and
  $\delta>0$, depending only on the relative positions of $\Lambda$ and
  $\Lambda'$ such that~:
  $$
    \forall (s,t)\in [a,b]\times[0,1],\
    \Vert f(s+it)\Vert\leq
    C \Vert f\Vert_{\infty(\strip[a,a+2\strip]\cup\strip[b-2,b\strip])}e^{-\delta \min(s-a,b-s)}
      $$
\end{lem}

\begin{proof}
  As before, let $\alpha(s)=\int_{0}^{1}\Vert f(s,t)\Vert^{2}dt$, and
  $\beta_{+}=\dot\alpha+2\delta\alpha$. Let $a'=\inf\{s\in[a,b],
  \beta_{+}(s)>0\}\cup\{b\}$. Then $\beta_{+}\leq 0$ on $[a,a']$ and
  $\alpha$ has an exponential decay~: $\forall s\leq a',
  \alpha(s)\leq\alpha(a)e^{-2\delta(s-a)}$. In the same way, let
  $\beta_{-}=\dot\alpha-2\delta\alpha$. Then
  $\dot\beta_{-}\geq-2\delta\dot\beta_{-}$, so that if
  $\beta_{-}(s_{0})<0$ then $\beta_{-}(s)<0$ for all $s<s_{0}$, and
  setting $b'=\sup\{s, \beta_{-}(s)<0\}\cup\{a\}$, $\alpha$ has an
  exponential growth on $[b',b]$~: $\forall s\in[b',b], \alpha(s)\leq
  e^{-2\delta(b-s)}\alpha(b)$. On $[a,a']\cup[b',b]$, we have
  \begin{equation}
    \label{eq:L2(s) estimate for bnded strips}
    \alpha(s)\leq e^{-2\delta(s-a)}\alpha(a)+e^{-2\delta(b-s)}\alpha(b)
  \end{equation}
  If  $a'<b'$, we still have to deal with $[a',b']$~: recall that
  $e^{-2\delta s}\beta_{+}(s)$ is increasing, so $\beta_{+}(s)\leq
  e^{-2\delta (b'-s)}\beta_{+}(b')$, and integrating once more,
  $e^{2\delta s}\alpha(s)\leq e^{2\delta a'}\alpha(a')+
  \frac{e^{-2\delta(b'-2s)}}{4\delta}\,\beta_{+}(s)$. Finally,
  $\alpha(s)\leq e^{-2\delta(s-a')}\alpha(a') +
  \frac{e^{-2\delta(b'-s)}}{4\delta}\,\beta_{+}(b')$. Moreover, on
  $[a',b']$ we have $-2\delta\alpha\leq\dot\alpha\leq2\delta\alpha$, and
  hence $\beta_{+}(b')\leq 4\delta\alpha(b')$, and \eqref{eq:L2(s)
    estimate for bnded strips} still holds for $s\in[a',b']$.

  Finally, the same argument applied to $\frac{\partial f}{\partial s}$
  gives the estimate
  \begin{equation}
    \Vert f(s,\cdot)\Vert_{1,2}
    \leq \Vert f(a,\cdot)\Vert_{1,2}\,e^{-\delta(s-a)}
       + \Vert f(b,\cdot)\Vert_{1,2}\,e^{-\delta(b-s)}~.~
  \end{equation}
  The proof then ends like in the previous lemma, replacing $a$ and $b$ by
$a+1$ and $b-1$, and using Schwarz lemma.
\end{proof}

\subsubsection{$\L^{1,p}$-smallness}

Let $C$ be a Floer trajectory close enough to $C_{\infty}$ and let
$\rho=\rho(C)$. Let $u:\Sigma_{\rho}\to M$ be the parametrization of $C$
such that, on each component, the first level $a_{i}$ which is encoutered
is reached at time $s=0$. Then, $u$ can be written in the form~:
$$
u(s,t)=\exp_{w_{\rho}(s,t)}(\xi(s,t))
$$
where $\xi$ is a section of $w_{\rho}^{*}TM$, satisfying appropriate
boundary conditions,\\
$\A\big((\exp_{w_{\rho}}\xi)(s_{i},\cdot)\big)=a_{i}$, and small in
$\L^{\infty}$ norm. We want to prove that $\xi$ is also small in
$\L^{1,p}$-norm.
\begin{lem}\label{lem:L1p smallness}
  For all $\epsilon>0$, there exist constants $R,\eta,\eta'>0$ such that
  if $C\in V_{\eta,\eta',R}(C_{\infty})$ (\textit{i.e.} $|\rho|>R$,
  $\Vert\xi\Vert_{\infty}<\eta$ and
  $\Vert\xi\Vert_{\RegClass^{1}(\theta_{\rho,i}(\strip[-R,R\strip]))}<\eta'$), then
  $\Vert\xi\Vert_{1,p}<\epsilon$.
\end{lem}

\begin{proof}
Consider a small neighbourhood $U_{\eta}=\cup U_{i}$ of the points
$x_{i}$, and a large $R>0$ such that
$u_{\infty,i}(\strip]-\infty,-R\strip[)\subset U_{i-1}$ and
$u_{\infty,i}(\strip]R,+\infty\strip[)\subset U_{i}$. For $\eta'$ small
enough, we have
$$
\Vert \xi_{|\theta_{\rho,i}(\strip[-R,R\strip])}\Vert_{1,p}<\epsilon.
$$

So we restrict attention now to the neighbourhood $U_{i}$ of $x_{i}$~:
$\Vert\xi\Vert_{1,p}$ has to be estimated on $\strip[s_{i}+R,s_{i+1}-R\strip]$,
or after a translation by $-s_{i}-\rho_{i}$, on
$A_{i}^{(\rho_{i}-R)}=\strip[-\rho_{i}+R,\rho_{i}-R\strip]$ (we suppose $\rho_{i}<+\infty$; the
other case is very similar).

Using a local chart, $U_{i}$ can be identified with a ball $B$ of
$\C^{n}$ so that $L$ and $L'$ are identified with two Lagrangian
linear subspaces intersected with $B$ and, moreover, the corresponding
induced almost complex structure $J$ coincides with the standard
almost complex structure at the origin. Indeed, it is standard that
there is a symplectic chart for $U_{i}$ which identifies $L\cap U_{i}$
and $L'\cap U_{i}$ with $B\cap \R^{n}$ and $B\cap i\R^{n}$,
respectively. By composing with a linear symplectic map we insure that
the condition on $J$ is satisfied and $L\cap U_{i}$, $L'\cap U_{i}$
are still identified with linear Lagrangians (obviously, not
orthogonal in general).  In such a chart, the (almost complex) Cauchy
Riemann equation becomes
\begin{equation}
  \label{eq:dbar u+q(u)du=0}
  \dbar u+q(u)\partial u =0.
\end{equation}
where $q=(J+{\i})^{-1}(J-{\i})$ satisfies $q(0)=0$ and $\dbar$ is, of course, associated to the standard complex structure
in $\C^{n}$, $i$.

After possibly rescaling a
smaller neighbourhood to the unit ball in $\C^{n}$, we may assume that
$\Vert q\Vert_{\RegClass^{1}}$ is as small as needed.

Notice that the relation
$$
\exp_{w_{\rho}}\xi=w_{\rho}+\xi'
$$
defines a new map $\xi':\strip[-\rho_i+R,\rho_i-R\strip]\to\C^n$
still satisfying appropriate boundary conditions. Since estimating
$\Vert\xi\Vert_{1,p}$ as a section of $w_{\rho}^*TM$ is equivalent to
estimating $\Vert\xi'\Vert_{1,p}$ (as a $\C^n$-valued function), we still
denote this $\xi'$ by $\xi$ in the sequel.
Thus, $\xi$ is now seen as a map to $\C^n$ instead
of a section of $w_{\rho}^*TM$ and $w_{\rho}+\xi$ is
$J$-holomorphic.
\medskip
Multiplying $\xi$ by appropriate cut-off functions, we obtain
$$
\hat\xi=\cutoff{-\rho+R+1}{-\rho+R}\cutoff{\rho-R-1}{\rho-R}\xi
$$
which is defined on the whole strip $\Sigma_{\rho}$, satisfies the boundary
conditions and belongs to $\L^{1,p}(\strip]-\infty,+\infty\strip[)$.
Lemma \ref{lem:operator P} gives the estimate~:
\begin{equation}
  \label{eq:|^xi|1p<c|dbar ^xi|p}
  \Vert\hat\xi\Vert_{1,p}\leq c\Vert\dbar\hat\xi\Vert_{p}
\end{equation}

Moreover, $\xi$ and $\hat\xi$ coincide on $A_{i}^{(\rho_{i}-R-1)}$ away
from a neighbourhood of the ends, where $\Vert\xi\Vert_{\RegClass^{1}}$
is controlled by $\eta'$~: therefore $\Vert \xi-\hat\xi\Vert_{1,p}\leq
2\Vert\chi\Vert_{\RegClass^{1}}
\eta' \leq 4\eta'$, so that
\begin{equation}
  \label{eq:|xi|1p<c|dbar xi|p+eps}
  \Vert\xi\Vert_{1,p}\leq c\Vert\dbar\xi\Vert_{p}+C\eta'
\end{equation}

To estimate $\Vert\dbar\xi\Vert_{p}$, write \eqref{eq:dbar u+q(u)du=0}
for $u=w_{\rho}+\xi$, and compare with $\alpha=\dbar
w_{\rho}+q(w_{\rho})\partial w_{\rho}$~:
\begin{equation}
  \label{eq:dbar xi complique}
\dbar\xi+
\big(q(w_{\rho}+\xi)-q(w_{\rho})\big)\partial w_{\rho}+
q(w_{\rho}+\xi)\partial\xi=-\alpha
\end{equation}
Recall from \eqref{eq:dbar w_rho} that $\alpha$ is small in $\L^{p}$-norm
for $\rho$ large enough. Developing $q(w_{\rho}+\xi)-q(w_{\rho})$ in the
form $a(z)\xi$, where $\Vert a\Vert_{\infty}$ is controlled by $\Vert
q\Vert_{\RegClass^{1}}$, and observing that $\Vert\partial
w_{\rho}\Vert_{\infty}$ is uniformly bounded, \eqref{eq:dbar xi
  complique} becomes
\begin{equation}
\label{eq:|dbar xi|p<a|xi|1p+eps}
\Vert\dbar \xi\Vert_{p}
\leq
\Vert \alpha \Vert_{p}
+a\Vert \xi\Vert_{1,p}
\end{equation}
where $a$ is a constant as small as desired. Collecting \eqref{eq:|dbar
  xi|p<a|xi|1p+eps} and \eqref{eq:|xi|1p<c|dbar xi|p+eps}, we obtain~:
\begin{equation}
  \big(1-a\,c\big)\Vert\xi\Vert_{1,p}
  \leq\Vert\beta\Vert+C\eta'
\end{equation}
Choosing our neighbourhoods $U_{i}$ small enough so that $a\,c<1$, we
obtain the desired estimate
\begin{equation}
  \Vert\xi\Vert_{1,p}
  \leq\epsilon
\end{equation}
for $\eta,\eta',R$ small/large enough.
\end{proof}

\subsubsection{The Floer equation for perturbations of $w_{\rho}$}\label{subsec:FloerPert}

Let $\Gamma^{1,p}_{\rho}(w_{\rho}^{*}TM)$ be the linear Banach space of $\L^{1,p}$ sections
$\xi$ of $w_{\rho}^{*}TM$ satisfying the boundary conditions and so that for all $i$
\begin{equation}
\label{eq:linear cut}
\int_{0}^{1}\omega(\frac{dw_{\rho}}{dt}(s_{i},t),\xi(s_{i},t))dt=0.
\end{equation}
(this condition is the linear version of
$\A(\exp_{w_{\rho}}\xi(s_{i},\cdot))=a_{i}$).  Let $\mathcal{H}_{i}$
be the (local) exponential image of those $C^{1}$ paths $\xi$
(depending on $t$ only) which verify (\ref{eq:linear cut}).  This is a
smooth hyper-surface in $\mathcal{P}=C^{1}_{L,L'}([0,1],M)$
(independent of $\rho$) and, by using the implicit function theorem in
a $C^{1}$ setting and on a bounded portion of the strip, we see that
each Floer trajectory $C$ in a sufficiently small, fixed neighbourhood
of $C_{\infty}$ crosses these hypersurfaces, thus defining ``linear''
transit times $\tilde{\rho}(C)$ so that $ |\tilde\rho(C)-\rho(C)|\leq
K \Vert\xi\Vert_{C^{1}(\theta_{\rho,i}(\strip[-R,R\strip]))} $ for
constants $K$ and $R$ independent of $C$. Thus, for each such curve
there is one and only one $\tilde\rho(C)>0$ and
$\xi\in\Gamma^{1,p}_{\tilde\rho}(w_{\tilde\rho}^{*}TM)$ such that
$u=\exp_{w_{\tilde\rho}}(\xi)$ is a parametrization of $C$.

Moreovoer, a section $\xi\in\Gamma^{1,p}_{\rho}(w_{\rho}^{*}TM)$ defines a
Floer trajectory if and only if it satisfies a non linear PDE,
$\dbartld \xi=0$, which is the translation in terms of $\xi$ of the usual
Floer equation \eqref{eq:J-holo} for $u=\exp_{w_{\rho}}\xi$~:
$$
\dbar_{J}[\exp_{w_{\rho}}\xi]=0.
$$
This expression takes values in $\Gamma(\Omega^{0,1}u^{*}TM)$.  Using
a $J$-hermitian connection, parallel transport along geodesics of
$g_{L,L'}$ defines an isomorphism
$\Pi_{J}:\Gamma(\Omega^{0,1}u^{*}TM)\to
\Omega^{0,1}(w_{\rho}^{*}TM)$. Finally, $\dbartld_{w_{\rho}}$ is
defined as~:
$$
\dbartld_{w_{\rho}}:
\begin{array}{ccc}
\Gamma_{\rho}(w_{\rho}^{*}TM)&\to&\Omega^{0,1}(w_{\rho}^{*}TM)\\
\xi&\mapsto&\Pi_{J}\Big(\dbar_{J}\big(\exp_{w_{\rho}}\xi\big)\Big)
\end{array}
$$

This map $\dbartld_{w_{\rho}}$ is as smooth as $J$, and one easily checks that all
its derivatives depend only on the derivatives of $J$ and $g$, and are
bounded independently of $\rho$. In particular, there is a constant
$A_{2}$ such that, for all $\rho$ large enough~:
\begin{equation}
  \label{eq:||dbartld||C2<A}
  \Vert\dbartld_{w_{\rho}}\Vert_{\RegClass^{2}}\leq A_{2}
\end{equation}

Moreover, for $\xi=0$, \eqref{eq:dbar w_rho} translates to~:
\begin{equation}
\label{eq:dbartld(0)<Ae^-rho}
\Vert\dbartld_{w_{\rho}} 0\Vert\leq A e^{-\lambda|\rho|}
\end{equation}

Let $D_{\rho}$ denote the linearisation of $\dbartld_{w_{\rho}}$ at $\xi=0$. Then
$D_{\rho}$ is Fredholm, and $\ind D_{\rho}=\ind(C_{\infty})$. Finally, since
the initial curve $C_{\infty}$ is supposed to be regular, for
$\rho=(+\infty,\dots,+\infty)$, $D_{\infty}$ is onto.

\subsubsection{Uniformly bounded right inverse}

\begin{prop}\label{prop:uniformly bounded right inverse}
 For $\rho$ large enough, the operator $D_{\rho}$ has a right inverse
 $P_{\rho}$, uniformly bounded with respect to $\rho$~:
 \begin{equation}
   \label{eq:P_rho unifly bounded}
   \exists C>0,\ \forall\rho>0,\forall\alpha\in\Omega^{0,1}(w_{\rho}^{*}TM)\quad
   \Vert P_{\rho}\alpha\Vert_{1,p}\leq C\Vert \alpha\Vert_{p}
 \end{equation}
\end{prop}

\begin{proof}
To construct a right inverse for $D_{\rho}$, we want to look at it as a
perturbation of $D_{\infty}$. Unfortunately, $D_{\rho}$ acts on
$w_{\rho}^{*}TM$ and the base of these bundles strongly depends on $\rho$.
To work around
this difficulty, $D_{\rho}$ has first to be brought into
$w_{\infty}^{*}TM$, where it will be compared to $D_{\infty}$.

To this end, consider the map~:
$$
u_{\rho}:\Sigma_{\infty}\to M
$$
obtained by multiplying $u_{\infty}$ by appropriate cutoff-functions to
make it constant away from
$\strip[\frac{\rho_{i-1}}{2},\frac{\rho_{i}}{2}\strip]$ on each component~:
$$
u_{\rho,i}(s,t)=
\begin{cases}
  \exp_{x_{i-1}}  (\cutoff{-(\rho_{i-1})/2+1} {-(\rho_{i-1})/2}
     \xi_{i}^{-})&\text{ if }s<-\frac{\rho_{i-1}}{2}+1\\
     u_{\infty,i}(s,t)&\text{ if }-\frac{\rho_{i-1}}{2}+1\leq s\leq \frac{\rho_{i}}{2}-1\\
  \exp_{x_{i+1}}(\cutoff{\rho_{i}/2-1}{\rho_{i}/2}
     \xi_{i}^{+})&\text{ if }s>\frac{\rho_{i}}{2}-1\\
\end{cases}
$$
Notice that $u_{\rho}$ and $w_{\rho}$ coincide where they are not
constant. As before, the Cauchy-Riemann equation for perturbations of
$u_{\rho}$ (satisfying \eqref{eq:linear cut} and the boundary conditions) leads to a PDE~:
$$
\dbartld_{u_{\rho}}:
\Gamma_{\rho}(u_{\rho}^{*}TM)
\to
\Omega^{0,1}(u_{\rho}^{*}TM),
$$
Let $\tilde D_{\rho}$ be its linearization at $0$. For $\rho=\infty$
(i.e. $\rho_{i}=+\infty,\forall i$), $\tilde D_{\infty}$ is just the
usual linearisation of $\dbar_{J}$ at $u_{\infty}=D_{\infty}$. Therefore, it
is onto (for a generic choice of $J$). For $\rho$ large enough,
parallel transport (using $g_{L,L'}$ on the left 
to preserve boundary conditions and a $J$-hermitian connection on the right to preserve
$(0,1)$ forms) induces isomorphisms
\begin{equation} \label{eq:u_rho*TM=u_infty*TM}
\begin{CD}
  \Gamma_{\rho}(u_{\rho}^{*}TM)@>\tilde D_{\rho}>>
  \Omega^{0,1}(u_{\rho}^{*}TM)\\
  @V\Pi_{L,L'} VV @VV\Pi_{J} V\\
  \Gamma_{\infty}(u_{\infty}^{*}TM)@>\tilde D_{\infty}>>
  \Omega^{0,1}(u_{\infty}^{*}TM)\\
\end{CD}
\end{equation}
and $\tilde D_{\rho}$ becomes a continuous family of operators on
$\Gamma(u_{\infty}^{*}TM)$. Thus, for $\rho$ large enough, $\tilde
D_{\rho}$ has a right inverse $\tilde R_{\rho}$, which is uniformly
bounded and continuous in $\rho$.

Now we want to come back to $D_{\rho}$. Consider the map $R_{\rho}$~:
$$
\Omega^{0,1}(w_{\rho}^{*}TM)
\xrightarrow{\mathrm{cut}_{\rho}}
\Omega^{0,1}(u_{\rho}^{*}TM)
\xrightarrow{\tilde R_{\rho}}
\Gamma_{\rho}(u_{\rho}^{*}TM)
\xrightarrow{\mathrm{glue}_{\rho}}
\Gamma_{\rho}(w_{\rho}^{*}TM)
$$
where $\mathrm{cut}_{\rho}$ is the extension by $0$ away from
$\rho_{i-1}<s<\rho_{i}$ on each component~:
$$
\forall(s,t)\in\Sigma_{\infty,i}:\
 \mathrm{cut}_{\rho}(\alpha)(s,t)=
 \begin{cases}
   \alpha(\theta_{\rho}^{-1}(s,t))&\text{ if }\rho_{i-1}\leq s\leq\rho_{i}\\
   0&\text{ otherwise, }
 \end{cases}
$$
and $\mathrm{glue}_{\rho}$ is the following gluing operation, where the
different component overlap on some region~: for
$(s,t)=\theta_{\rho}(z)\in\Sigma_{\infty,i}$, set
$$
 \mathrm{glue}_{\rho}(\xi)(z)=
 \begin{cases}
   \xi_{i}(s,t)+\cutoff{-(\rho_{i-1})/2-1}{-(\rho_{i-1})/2}\ \xi_{i-1}(s',t)
   &\text{ if }-\rho_{i-1}\leq s\leq-\rho_{i-1}/2\\
   \xi_{i}(s,t)&\text{ if }-(\rho_{i-1})/2\leq s\leq\rho_{i}/2\\
   \xi_{i}(s,t)+\cutoff{\rho_{i}/2+1}{\rho_{i}/2}\xi_{i+1}(s'',t)
   &\text{ if }\rho_{i}/2\leq s\leq\rho_{i}
 \end{cases}
$$
where $s'=s+s_{i}-2\rho_{i-1}-s_{i-1}$ and
$s''=s+s_{i}+2\rho_{i}-s_{i+1}$, so that $(s',t)$ and $(s'',t)$ are the
value of $z$ seen in $\Sigma_{\infty,{i-1}}$, and $\Sigma_{\infty,i+1}$ (see figure
\ref{fig:gluing with overlapping}).
\begin{figure}[htbp]
  \begin{center}
\setlength{\unitlength}{1mm}%
\begin{picture}(110,30)(-6,-5)
\def\drawstrip(#1,#2)#3#4{%
  \put(#1,#2){%
    \multiput(0,0)(0,5){2}{
      \multiput(0,0)(-1,0){5}{\line(1,0){.5}}
      \put(0,0){\line(1,0){#3}}
      \multiput(#3,0)(1,0){5}{\line(1,0){.5}}
      }
    \put(-2.5,2){\makebox(0,0)[c]{#4}}
    }
  }
\def\drawstripbar(#1,#2)#3{%
  \put(#1,#2){\line(0,1){5}\put(0,-.5){\makebox(0,0)[t]{#3}}}
  }
\def\drawstripbardot(#1,#2)#3{
  \put(#1,#2){\multiput(0,0.3)(0,1){5}{\line(0,1){.25}}
    \put(0,-.5){\makebox(0,0)[t]{#3}}
    }
  }
\def\cutoff(#1,#2){
  \put(#1,#2){%
    \bezier{0}(0,10)(1,10)(1,5)
    \bezier{0}(1,5)(1,0)(2,0)
    }
  }
\def\cuton(#1,#2){
  \put(#1,#2){%
    \bezier{0}(0,0)(1,0)(1,5)
    \bezier{0}(1,5)(1,10)(2,10)
    }
  }
\drawstrip( 0, 0){100}{$\Sigma_{\rho}$}
\drawstripbar(10,0){$s_{i-1}$}
\drawstripbar(40,0){$s_{i}$}
\drawstripbar(90,0){$s_{i+1}$}
\drawstripbardot(17.5, 0){}\drawstripbardot(17.5, 5){}
\drawstripbardot(25.0, 0){}
\drawstripbardot(32.5, 0){}\drawstripbardot(32.5, 5){}
\drawstripbardot(52.5, 0){}\drawstripbardot(52.5, 5){}
\drawstripbardot(65.0, 0){}
\drawstripbardot(77.5, 0){}\drawstripbardot(77.5, 5){}
\put(0,6){
  \put(-2.5,0){\makebox(0,10)[c]{$\mathrm{glue}_{\rho}(\xi)$}}
  \put(0,10){\line(1,0){30.5}}
  \put(7.75,12){\makebox(0,0)[c]{$\xi_{i-1}$}}
  \cutoff(30.5,0)
  \put(32.5,0){\line(1,0){5}}
  }
\put(0,6){
  \put(52.5,0){\line(-1,0){5}}
  \cuton(52.5,0)
  \put(54.5,10){\line(1,0){45}}
  \put(89.5,12){\makebox(0,0)[c]{$\xi_{i+1}$}}
  }
\put(0,7){
  \put(17.5,0){\line(-1,0){5}\cuton(0,0)}
  \put(19.5,10){\line(1,0){56}}
  \put(42.5,12){\makebox(0,0)[c]{$\xi_{i}$}}
  \cutoff(75.5,0)
  \put(77.5,0){\line(1,0){5}}
  }
\multiput(46.25,0)(12.5,0){4}{
  \put(-1.25,2){\vector(-1,0){5}}
  \put( 1.25,2){\vector( 1,0){5}}
  \put(0,0){\makebox(0,5)[c]{$\scriptscriptstyle\frac{\rho_{i}}{2}$}}
  }
\end{picture}
    \caption{The $\mathrm{glue}_{\rho}$ map.}
    \label{fig:gluing with overlapping}
  \end{center}
\end{figure}

Since the three operators $\mathrm{cut}_{\rho}$, $\tilde R_{\rho}$ and
$\mathrm{glue}_{\rho}$ are uniformly bounded in $\rho$, so is $R_{\rho}$.

The proof of the proposition is finished by the next lemma.

\begin{lem}
The operator $R_{\rho}$ is a quasi inverse for $D_{\rho}$, in the sense that~:
$$
 \lim_{|\rho|\to+\infty}\Vert D_{\rho}\circ R_{\rho}-\Id\Vert=0
$$
\end{lem}
\begin{proof}
Let $\alpha\in\Omega^{01}(w_{\rho}^{*}TM)$ and
$\beta=D_{\rho}R_{\rho}\alpha-\alpha$. We have to estimate
$\Vert\beta\Vert_{p}$. Notice that $\beta$ is supported on the gluing
regions $\strip[s_{i}+\rho_{i}/2,s_{i+1}-\rho_{i}/2\strip]$. Let us focus
on one half of such a region
$\strip[s_{i}+\rho_{i}/2,s_{i}+\rho_{i}\strip]$. Let
$(\alpha_{1},\dots,\alpha_{k})=\mathrm{cut}_{\rho}{\alpha}$ and
$\eta=(\eta_{1},\dots,\eta_{k})=\tilde
R_{\rho}(\mathrm{cut}_{\rho}(\alpha))$. We have~:
\begin{align*}
 \beta(z)
&=D_{\rho}(\mathrm{glue}_{\rho}\eta)(z)-\alpha(z)\\
&=D_{\rho}\big(\eta_{i}(s,t)+\cutoff{}{i}\eta_{i+1}(s'',t)\big)-\alpha(z)
\end{align*}

Notice that on the domain under consideration (and modulo appropriate
translation) $w_{\rho}$, $u_{i,\rho}$, and $u_{i+1,\rho}$ are all
constant, so that $D_{\rho}=\tilde{D}_{\rho}=\dbar$, where $\dbar$ is the
usual Cauchy-Riemann operator associated to the complex vector space
$(T_{x_{i+1}}M,J(x_{i+1}))$. Finally, we obtain
\begin{align*}
 \beta(z)
&=\alpha_{i}(s,t)
+(\dbar\cutoff{}{})\eta_{i+1}(s'',t)
+\cutoff{}{}\alpha_{i+1}(s'',t)
-\alpha(z)
\end{align*}
But $\alpha_{i}(s,t)=\alpha(z)$ and $\alpha_{i+1}(s'',t)=0$. Hence~:
$\beta(z) =(\dbar\cutoff{}{})\eta_{i+1}(s'',t)$ and
\begin{equation}
  \label{eq:|beta|p<A|eta|infty}
 \Vert\beta\Vert_{p}\leq
   A\Vert\eta_{i+1}\Vert_{\infty,\strip[-\rho_{i}/2,-\rho_{i}/2+1\strip]}
\end{equation}

Moreover, we have $\dbar\eta_{i+1}=0$ on $\strip[-\infty,0\strip]$, so
applying lemma \ref{lem:std holo half strips} to $\eta_{i+1}$, we obtain
that $\Vert\beta\Vert_{p}\leq C
\Vert\eta_{i+1}\Vert_{1,p}e^{-\delta\rho_{i}/2}\leq C'
\Vert\alpha\Vert_{p}e^{-\delta\rho_{i}/2}$ with uniform constants $C$ and
$\delta$. Gathering all these inequalities, we obtain uniform constants
$C$ and $\delta$, such that
$$
 \Vert D_{\rho}R_{\rho}\alpha-\alpha\Vert_{p}
\leq Ce^{-\delta|\rho|}\ \Vert\alpha\Vert_{p}.
$$
\end{proof}

This ends the proof of proposition \ref{prop:uniformly bounded right
  inverse}, since, for $\rho$ large enough, $D_{\rho}\circ R_{\rho}$ is
uniformly invertible and we can set $P_{\rho}=R_{\rho}\circ(D_{\rho}\circ
R_{\rho})^{-1}$
\end{proof}

\subsubsection{Isomorphism from $\ker D_{\rho}$ to $\ker D_{\infty}$}

Finally, we need to identify all $\ker D_{\rho}$ to the constant space
$\ker D_{\infty}$. To this end, consider a small neighbourhood $U$ of
the intersection points $x_{i}$, and a compact subset $K =
\bigsqcup_{i=1}^{k} \strip[-R,R\strip]_{i} \subset \Sigma_{\infty}$
with $R$ large enough to have $u_{\infty}(\Sigma_{\infty}\setminus
K)\subset U$.  Then we have (\cite{Siko1})~:

\begin{prop}\label{prop:kerDiso}
  Let $\pi:\Gamma(w_{\rho}^{*}TV)\to\Gamma(w_{\rho}^{*}TV_{|_{K}})\to\ker
  D_{\infty}$ be the $\L^{2}$ orthogonal projection on $\ker D_{\infty}$.
  Then, for $|\rho|$ large enough, the restriction of $\pi$ to $\ker
  D_{\rho}$ is an isomorphism, and there is a constant $C$, uniform with
  respect to $\rho$, such that
  $$
   \Vert \xi\Vert_{1,p}\leq C\Vert\pi(\xi)\Vert_{2}
  $$
\end{prop}

\begin{proof}
  Suppose there is a sequence $(\rho_{n},\xi_{n})$ such that
  $D_{\rho_{n}}\xi_{n}=0$, $\Vert \xi_{n}\Vert_{1,p}=1$ and
  $\lim\pi(\xi_{n})=0$. Then a sub-sequence of $\xi_{n}$ converges on all
  compacts to a section $\xi\in\ker D_{\infty}$. On the other hand,
  $\pi(\xi_{n})$ converges to $\pi(\xi)$, and hence $\xi_{|_{K}}\in(\ker
  D_{\infty})^{\bot}$. Hence, $\xi_{|_{K}}=0$ and $\xi_{n}$ converges to
  $0$ uniformly on $K$ in the $\RegClass^{\infty}$ topology.

  We will derive a contradiction from lemma \ref{lem:std holo bnded
    strips}, which implies that $\Vert\xi_{n}\Vert_{1,p}$ on the
  complement of $K$ is controlled by the behaviour of $\xi$ on the
  boundary and, therefore, should tend to $0$ as $n$ goes to infinity.

  To make this explicit, let us focus on one component
  $[-\rho_{n,i}+R,\rho_{n,i}-R]$ of the complement of $K$. On this piece
  of strip, $w_{\rho}$ takes values in a small neighbourhood of $x_{i}$
   where we choose coordinates in which $L$ 
  and $L'$ are linear spaces and $J(0)=\i$ like in the proof of
  \ref{lem:L1p smallness}. The Cauchy-Riemann equation takes the form
  $$
  \dbar u+q(u)\partial u=0.
  $$
As before, we can replace $\xi_{n}$ by $\xi'_{n}$ with
  $\exp_{w_{\rho_{n}}}\xi_{n}=w_{\rho_{n}}+\xi'_{n}$, keeping control on the $\L^{1,p}$
  norm. Then $\xi_{n}$ satisfies an equation of the form:
  $$
  \dbar\xi_{n}+q_{n}(z)\partial\xi_{n}+a_{n}(z)\xi_{n}=0
  $$
  where $q_{n}(z)=q(w_{\rho_{n}})$ and $a_{n}(z)(\cdot)=D_{w_{\rho_{n}}}q(\cdot)\partial w_{\rho_{n}}$ are uniformly small.

  Using the operator $P$ defined in lemma \ref{lem:operator P}, 
  let $\eta_{n}=\xi_{n}+P(q_{n}(z)\partial\xi_{n}+a_{n}(z)\xi_{n})$ (notice that the operator
  $P$ of the lemma is defined on the full strip, so to be more precise, we
  first extend $q_{n}\partial\xi_{n}+a_{n}\xi_{n}$ by $0$ and then consider the restriction
  of its image under $P$ to our piece of strip). This section $\eta_{n}$ is 
  holomorphic, and takes values in $L$ and $L'$ on the boundary. Moreover,
  we obtain a uniformly small constant $\kappa$ such that, on the relevant
  piece of strip:
  \begin{equation}
    \label{eq:eta_n-xi_n}    
    \Vert \eta_{n}-\xi_{n}\Vert_{1,p}\leq \kappa \Vert\xi_{n}\Vert_{1,p}
  \end{equation}
  As a consequence, near the ends of our piece of strip 
  $\strip[-\rho_{n,i}+R,\rho_{n,i}-R\strip]$, $\Vert \eta_{n}\Vert_{\infty}$ is arbitrarily
  small. According to lemma \ref{lem:std holo bnded strips}, we get
  $\Vert\eta_{n}\Vert_{\infty}\leq C \Vert\eta_{n}\Vert_{\infty,A}e^{-\delta(\rho_{n}-R-|s|)}$, where 
  $A=\strip[-\rho_{n,i}+R,-\rho_{n,i}+R+2\strip]\cup\strip[\rho_{n,i}-R-2,\rho_{n,i}-R\strip]$.
  Integrating this, we get $\Vert \eta_{n}\Vert_{p}\leq C\Vert\eta_{n}\Vert_{\infty,A}$ for
  some uniform constant $C$. In the same way, the Schwartz lemma provides 
  a control of $\partial\eta_{n}$ near the ends by means of $\Vert\eta_{n}\Vert_{\infty}$, 
  from where we derive a similar estimate 
  $\Vert\partial\eta_{n}\Vert_{p}\leq C\Vert\eta_{n}\Vert_{\infty,A}$.
  
  Finally, $\Vert \eta_{n}\Vert_{1,p}$ is arbitrarily small on 
  $\strip[-\rho_{n,i}+R,\rho_{n,i}-R\strip]$, which, in view of \eqref{eq:eta_n-xi_n}, contradicts
  $\Vert\xi_{n}\Vert_{1,p}=1$.

  \end{proof}

\subsubsection{End of the proof}

We will use the following version of the implicit function theorem~:
\begin{prop}
  Let $(f_{\lambda}:E_{\lambda}\to F_{\lambda})_{\lambda\in [0,1[^{m}}$ be
  a family of maps between Banach spaces such that
  \begin{enumerate}
  \item \label{it:f C2-bnded}
    for all $\lambda>0$, $f_{\lambda}$ is of class $\RegClass^{2}$, and
    $\Vert f_{\lambda}\Vert_{\RegClass^{2}}$ is uniformly bounded.

  \item \label{it:f(0)->0}
    $\lim_{\lambda\to0}f_{\lambda}(0)=0$,

  \item \label{it:Df(0) unifly invble}
    $Df_{\lambda}(0)$ is uniformly invertible~: $\exists R_{\lambda}\in
    L(F_{\lambda},E_{\lambda}), Df_{\lambda}(0)\circ R_{\lambda}=\Id$ and
    $\exists C>0,\forall \lambda\ \Vert R_{\lambda}\Vert\leq C$.
Let $H_{\lambda}=R_{\lambda}(F_{\lambda})$.

  \end{enumerate}
Then there exists $\epsilon>0$ such that for all $\lambda$ with 
$|\lambda|<\epsilon$, there are ``uniform'' open subsets 
$U_{\lambda}\subset \ker Df_{\lambda}(0)$ and $V_{\lambda}\subset 
E_{\lambda}$,
and a diffeomorphism $\varphi_{\lambda}:U_{\lambda}\to 
H_{\lambda}$ such that, in the decomposition $E_{\lambda}=\ker 
Df_{\lambda}\oplus H_{\lambda}$, one has:
$$
  f_{\lambda}(x,y)=0 \Leftrightarrow y=\psi_{\lambda}(x).
$$
Here, 
uniform means that the $U_{\lambda}$ all contain a ball whose 
radius is independent of $\lambda$.
\end{prop}

\begin{proof}
  Since $Df_{\lambda}(0)$ has a right inverse, $Df_{\lambda}$ remains
  onto on a neighbourhood $V_{\lambda}$ of $0$ whose size is controlled
  by $\Vert R_{\lambda}\Vert$ and $\Vert
  f_{\lambda}\Vert_{\RegClass^{2}}$. The conditions \ref{it:f C2-bnded}
  and \ref{it:Df(0) unifly invble} imply that this size is uniform in
  $\lambda$. Finally, the Newton algorithm proves that
  $f_{\lambda}^{-1}(0)\cap V_{\lambda}\neq \emptyset$ for $\lambda$ small
  enough, and the usual implicit function theorem gives the result.
\end{proof}

The three conditions in this proposition have been checked for
$\dbartld_{\rho}: \Gamma^{1,p}_{\rho}(w_{\rho}^{*}TM) \to
\Gamma^{p}(\Lambda^{0,1}\Sigma_{\rho}\otimes w_{\rho}^{*}TM)$ in
\eqref{eq:||dbartld||C2<A} and \eqref{eq:dbartld(0)<Ae^-rho} and the
proposition \ref{prop:uniformly bounded right inverse}.

However, to have a better control on the spaces $H_{\rho}$, we will
use a slightly modified right inverse for $D_{\rho}$. Recall the
compact $K$ used in proposition \ref{prop:kerDiso}, and consider the
restriction map $L^{1,p}(\Sigma_{\rho})\xrightarrow{r}L^{1,p}(K)$ and
the projection $L^{1,p}(K)\xrightarrow{\pi_{\rho}}\ker D_{\rho}$ in
the (constant) direction of the $L^{2}$-orthogonal of $r(\ker
D_{\infty})$. The operator $P'_{\rho}=P_{\rho}-\pi_{\rho}(r(P_{\rho}))$
is then still a right inverse of $D_{\rho}$, it is still uniformly
bounded, but has the additional property that the corresponding space
$H_{\rho}=\mathrm{rk}(P'_{\rho})$ is such that
$\pi_{\rho}(r(H_{\rho}))=0$.

We obtain a one to one map
$\varphi_{\rho}=\mathrm{id}\oplus\psi_{\rho}$ from a neighbourhood of
$0$ in $\ker D_{\rho}$ to the space of sections in
$\Gamma^{1,p}_{\rho}(w_{\rho}^{*}TM)$ associated to Floer trajectories
in a neighbourhood of $C_{\infty}$ in $\overline\M(x_{k},x_{0})$.

Composing $\varphi_{\rho}^{-1}$ with the map $\pi:\ker D_{\rho}\to\ker D_{\infty}$,
constructed in the previous section, we obtain a map
$\phi=(\rho,\pi\circ\varphi_{\rho}^{-1})$ between neighbourhoods of $C_{\infty}$ in
$\overline\M(x_{k},x_{0})$ and $(\infty,\dots,\infty,0)$ in $(1,+\infty]^{k-1}\times
T_{C_{\infty}}\M(x_{k},\dots,x_{0})$ with $\rho$ the linear transit times as described
at the beginning of \S\ref{subsec:FloerPert}~:
$$
\phi:\overline\M(x_{k},x_{0}),C_{\infty}
\xrightarrow{(\rho,\pi\circ\varphi_{\rho}^{-1})}
(1,+\infty]^{k-1}\times T_{C_{\infty}}\M(x_{k},\dots,x_{0}),((\infty,\ldots, \infty),0)
$$

\medskip

This map is one to one. For the continuity of $\phi^{-1}$ consider a
converging sequence $(\rho_{n},\xi_{n})\to(\rho,\xi)$. By possibly
extracting a subsequence we may assume that
$C_{n}=(\phi^{-1}(\rho_{n},\xi_{n}))$ converges to a curve $C'$
corresponding to $\xi'\in\ker D_{\rho}$.  The curves $C_{n}$ may be
viewed as sections $x_{n}+y_{n}\in L^{1,p}(w_{\rho_{n}}^{*}(TM)) =
\ker D_{\rho_{n}} \oplus H_{\rho_{n}}$, and, similarly, the curves
$C=\phi^{-1}(\rho,\xi)$ and $C'$ as sections $x+y$ and $x'+y'$.

We now have $\pi(r(y_{n}))=\pi(\pi_{\rho}(r(y_{n})))=0$, so that
$\pi(r(x_{n}+y_{n}))=\pi(r(x_{n}))=\xi_{n}$. On the other hand, this
has to converge to $\pi(r(x'+y'))=\xi'$, so $\xi=\xi'$, and finally
$C=C'$.

Conversely, if $C_{n}\to C$ is a converging sequence of trajectories,
then $\rho(C_{n})$ converges to $\rho(C)$ and $\xi_{n}$ is bounded in
the finite dimensional space $\ker D_{\infty}$. On the other hand, a
converging subsequence of $(\rho_{n},\xi_{n})$ has to converge to
$(\rho,\xi)=\phi(C)$, and hence $\phi(C_{n})$ converges to $\phi(C)$.

Let us turn to the behaviour of $\phi$ with respect to the
stratification~: $\phi$ clearly respects the stratifications since a
trajectory $C$ is broken at $x_{i}$ if and only if $\rho_{i}(C)=+\infty$.
Within each stratum now, all the fiber-bundles $w_{\rho}^{*}TM$ are
topologically equivalent, and we can locally move smoothly from one to
another. The whole family $\pi\circ\varphi_{\rho}^{-1}$ depends then smoothly
on $\rho$ and $\phi$ is a local diffeomorphism.

This ends the proof of theorem \ref{thm:corner structure on Mbar}.

\subsection{Hamiltonian perturbations}

We now shortly discuss the non homogeneous case.  In this case the
usual Cauchy-Riemann equation \eqref{eq:J-holo} is modified by a Hamiltonian
term:
\begin{equation}
  \label{eq:non homogeneous J-holo}
\frac{\partial u}{\partial s}+J(u)\frac{\partial u}{\partial t}=-\nabla H_{t}(u)
\end{equation}
or, equivalently~:
\begin{equation}
  \label{eq:du+Jdui=b}
\dbar_{J}u=du+J(u)du\i=b_{u}
\end{equation}
where $b_{u}$ is the $\C$ anti-linear map defined by $b_{u}(z)\ \xi
=-\bar{\xi}\ \nabla H(u(z))$.

The homogeneous and non homogeneous situations differ in many respects.
In particular, the intersection points are replaced by orbits of the
hamiltonian flow $\phi^{t}$ of $H$ starting on $L$ and reaching $L'$ at
time $1$ (the transversality assumption being replaced by requiring that
$\phi_{1}^{H}(L)$ be transverse to $L'$). Moreover, the ``breaks'' do not arise in the
neighbourhood of a point, but along a curve, making the analysis a bit
more technical. 

\subsubsection{Reduction to the standard non homogeneous equation}

However, using the naturality maps used in \S\ref{subsubsec:naturality}, this case reduces
to the one in which the basic equation is homogenous but the almost complex structure
depends on the variable $t$.  This is the case that we discuss below.
Thus, the model equation of lemmas
\ref{lem:operator P}, \ref{lem:std holo half strips} and \ref{lem:std
  holo bnded strips} is now replaced by
$$
\frac{\partial f}{\partial s}+J(t)\frac{\partial f}{\partial t}=0.
$$

Considering a path $\Phi_{t}$ of symplectic matrices such that
$J(t)=\Phi^{-1}_{t}\i\Phi_{t}$, and letting
$g(s,t)=\Phi_{t}f(s,t)$, we end up with the equation:
$$
\dbar g = b(t) g 
$$
and the boundary conditions become $g(s,0)\in\Lambda$ and $g(s,1)\in\Lambda''=\Phi_{1}\Lambda'$.
Notice that the transversality
assumption on $\Lambda$ and $\Lambda'$ is now replaced by requiring that
the differential equation for $\gamma:[0,1]\to(\C^{n},\Lambda,\Lambda'')$~:
$$
\dot\gamma=\i b(t)\gamma
$$
has no non trivial solution.

Finally, lemma \ref{lem:operator P} in this setting is again a boundary version of Theorem 3.1.13
from \cite{Sch3}, as for lemmas \ref{lem:std
  holo half strips} and \ref{lem:std holo bnded strips}, observe that
$b$ is self adjoint, and replacing the operator
$A=\i\frac{\partial}{\partial t}$ in the proof of these lemmas by
$$
 A=\i\frac{\partial}{\partial t}+b
$$
we get a self adjoint, injective operator, so that
$\alpha(s)=\int_{0}^{1}\Vert f\Vert^{2}dt$ satisfies the same differential
inequality as before. 

This discussion also covers the case when dependence of $s$ is required (for example,
to study Floer chain maps) because $s$-dependence has compact support, where the
convergence of the curves is well controlled.

\end{document}